\documentclass[11pt]{article}
\usepackage{graphicx,epstopdf,subfigure,mathtools,mathrsfs, arydshln} 
\usepackage[font=small,labelfont=bf]{caption}
\usepackage{float, cancel, amssymb}
\usepackage[normalem]{ulem}
\PassOptionsToPackage{usenames,dvipsnames}{xcolor}
\usepackage[usenames,dvipsnames]{xcolor}
 
\usepackage[margin=1 in]{geometry}

\usepackage{amsfonts, soul, authblk}

\newcommand{\rev}[1]{\color{black}#1\normalcolor}

\newcommand{\revp}[1]{\color{black}#1\normalcolor}
\newcommand{\widebar}[1]{%
   \hbox{%
     \vbox{%
       \hrule height 0.5pt 
       \kern0.5ex
       \hbox{%
         \kern-0.1em
         \ensuremath{#1}%
         \kern-0.1em
       }%
     }%
   }%
} 
\graphicspath{{../}{LocalFigs/}{CLinkingFigs/}{NLFigs/}{LudQuadFigs/}{PaperFigs/}}
\makeatletter
\def\input@path{{PaperTex/}}
\makeatother
\newfloat{algorithm}{!t}{loa}
\providecommand{\algorithmname}{Algorithm}
\floatname{algorithm}{\protect\algorithmname}

\newcommand{\V}[1]{\boldsymbol{#1}}                 
\newcommand{\M}[1]{\boldsymbol{#1}}
\newcommand{\Lop}[1]{\boldsymbol{\mathcal{#1}}}
\newcommand{\bm}[1]{\boldsymbol{#1}}
\newcommand{\ind}[2]{{#1}^{(#2)}}
\newcommand{\piX}[2]{\V{X}^{\left(#1\right)}_{ #2}}
\newcommand{\piXs}[2]{\V{\tau}^{(#1)}_{ #2}}
\newcommand{\iXs}[1]{\V{\tau}^{(#1)}}
\newcommand{\pif}[2]{\V{f}^{(#1)}_{ #2}}
\newcommand{\piThree}[3]{{#1}^{(#2)}_{ #3}}
\newcommand{\ddt}[1]{\frac{\partial #1}{\partial t}}

\newcommand{\Knels}[3]{\M{S}_D\left(#1, #2; #3\right)}
\newcommand{\KRPY}[3]{\M{S}_\text{RPY}\left(#1, #2; #3\right)}
\newcommand{\tdisc}[2]{#1_{#2}}
\newcommand{\dco}{\frac{e^3}{24}}

\global\long\def\D#1{\Delta#1}

\global\long\def\norm#1{\left\Vert #1\right\Vert }

\global\long\def\bm#1{\V#1}
\global\long\def\eqd{:=}
\global\long\def\Ubar{\widebar{\V{U}}}

\global\long\def\Slet#1{\M{S}\left(#1\right)}
\global\long\def\Dlet#1{\M{D}\left(#1\right)}
\global\long\def\Knel#1{\M{S}_D\left(#1\right)}
\global\long\def\Xs{\V{\tau}}
\global\long\def\EPMI{\frac{1}{8\pi\mu}}
\global\long\def\ML{ \M{M}^\text{L}}
\global\long\def\MLD{ \M{M}^\text{LD}}
\global\long\def\MNL{ \M{M}^\text{NL}}
\global\long\def\MFP{\Lop{M}^\text{FP}}
\global\long\def\MFPd{\M{M}^\text{FP}}
\global\long\def\MJF{\Lop{M}^\text{c}}
\global\long\def\MJFd{\M{M}^\text{c}}
\global\long\def\Xcol{\V{X}}
\global\long\def\aRPY{b}
\global\long\def\McNot{\left(\MJF\left[\ind{\V{X}}{j}(\cdot)\right]\ind{\V{f}}{j}(\cdot)\right]\left(\ind{\V{X}}{i}(s)\right)}
\global\long\def\McNotd{\V{v}\left(\V{x}\right)}

\global\long\def\eqd{:=}

\usepackage{hyperref}
\hypersetup{
    colorlinks=false,
    pdfborder={0 0 0},
}

\title{An integral-based spectral method for inextensible slender fibers in Stokes flow}
\author[a]{Ondrej Maxian\thanks{Email: om759@nyu.edu}}
\author[a,b]{Alex Mogilner\thanks{Email: mogilner@cims.nyu.edu}}
\author[a]{Aleksandar Donev\thanks{Email: donev@courant.nyu.edu}}
\affil[a]{Courant Institute of Mathematical Sciences, New York University, New York, NY 10012}
\affil[b]{Department of Biology, New York University, New York, NY 10012}
\begin{document}

\maketitle

\begin{abstract}
Every animal cell is filled with a cytoskeleton, a dynamic gel made of inextensible fibers, such as microtubules, actin fibers, and intermediate filaments, all suspended in a viscous fluid. Numerical simulation of this gel is challenging because the fiber aspect ratios can be as large as $10^4$. We describe a new method for rapidly computing the dynamics of inextensible slender filaments in periodically-sheared Stokes flow. The dynamics of the filaments are governed by a nonlocal slender body theory which we partially reformulate in terms of the Rotne-Prager-Yamakawa hydrodynamic tensor. To enforce inextensibility, we parameterize the space of inextensible fiber motions and strictly confine the dynamics to the manifold of inextensible configurations. To do this, we introduce a set of Lagrange multipliers for the tensile force densities on the filaments and impose the constraint of no virtual work in an $L^2$ weak sense. We augment this approach with a spectral discretization of the local and nonlocal slender body theory operators which is linear in the number of unknowns and gives improved spatial accuracy over approaches based on solving a line tension equation. For dynamics, we develop a second-order semi-implicit temporal integrator which requires at most a few evaluations of nonlocal hydrodynamics and a few block diagonal linear solves per time step. After demonstrating the improved accuracy and robustness of our approach through numerical examples, we apply our formulation to a permanently cross-linked actin mesh in a background oscillatory shear flow. We observe a characteristic frequency at which the network transitions from quasi-static, primarily elastic behavior to dynamic, primarily viscous behavior. We find that nonlocal hydrodynamics increases the viscous modulus by as much as 25\%. \revp{Most of this increase, in contrast to the smaller (about 10\%) increase in the elastic modulus, is due to short-ranged intra-fiber interactions. } 
\end{abstract}

\section{Introduction \label{sec:intro}}
Interactions of long, thin, inextensible filaments with a viscous fluid abound in biology, engineering, physics, and medicine. In biology, the swimming mechanisms of flagellated organisms have been of interest for decades, with an initial cluster of studies on how force and torque balances lead to swimming \cite{chwang1971note, berg1973bacteria, brennen1977fluid, lauga2009hydrodynamics}, and a more recent focus on flagellar bundling and propulsion \cite{lim2012fluid, maier2016magnetic, maria2020geometrical}. In physics and engineering, suspensions of high-aspect-ratio fibers have been observed to display non-Newtonian, viscoelastic behavior both experimentally \cite{fibexps} and computationally \cite{mackfibs, fibsusps}. 

Our particular area of interest is the simulation of semi-flexible filaments that make up the cell cytoskeleton. These inextensible filaments, which include microtubules and actin filaments, maintain the cell's structure, control the mechanics of the cell division process, and have aspect ratios from $10^2$ to $10^4$ \cite{alberts}. In vivo, actin filaments are generally bound together into networks by cross-linking proteins, the properties of which determine the viscoelastic behavior of the cytoskeleton \cite{wagner2006cytoskeletal, head2003deformation, ahmed2015dynamic}. While there has been much work recently on microtubule systems \cite{ehssan17, shelley2016dynamics}, there has yet to be, to our knowledge, a systematic study of the influence of hydrodynamic interactions on the mechanics and rheology of cross-linked actin networks. One of the goals of this paper is to develop an efficient numerical technique that can simulate a cross-linked network of thousands of inextensible actin filaments and take into account the filament interactions with a viscous solvent at zero Reynolds number. Our method not only handles the inextensibility and stiffness of the fibers robustly and efficiently, but also accounts for, in near-linear time with respect to the number of fibers, the long-ranged hydrodynamic interactions between fibers, which we show can increase the viscous modulus of the network by as much as 25\%. We will not consider thermal fluctuations or elastic twisting of the filaments in this work. 

Prior to the year 2000, tools for analytical analysis and numerical simulation of filaments in Stokes flow were developed in parallel by several authors. For slender filaments, a useful approach for both analysis and computation is to reduce the problem from three dimensions to one by assuming a certain distribution of singularities along the filament centerline. This approach, referred to as ``slender body theory'' (SBT), was first introduced by Hancock \cite{hancock1953self} and later expanded upon by Batchelor \cite{batchelor1970slender}. By using the method of matched asymptotics, Keller and Rubinow were the first to derive an SBT that is uniformly accurate in the fiber slenderness ratio $\epsilon = $ radius/length \cite{krub}. Johnson further developed the theory by adding higher order corrections and correctly treating a fiber with free ends \cite{johnson}, and G\"otz re-derived the SBT of Keller and Rubinow in a more general context, allowing him to apply the theory to Oseen's and Poisson's equations \cite{gotz2001interactions}. 

Because the SBTs of Keller and Rubinow, Johnson, and G\"otz are uniformly accurate in powers of $\epsilon$, they have formed the basis of most of the more recent analysis. To this end, Mori et al.\ recently showed that these singularity solutions solve a well-posed Stokes problem with non-standard boundary conditions on the filament surface \cite{mori2018theoretical, morifree}. Koens and Lauga also showed that the SBT singularity solution can be recovered by matched asymptotic expansion of the full surface boundary integral formulation of Stokes flow \cite{koens2018boundary}. 

On the numerical side, non-SBT based techniques for the simulation of fibers in Stokes flow have been in use for many decades. The most prevalent among these are regularized singularity methods, in which the fibers are discretized by a series of marker points, each of which is assigned a force according to the fiber physics. The force on each marker is then regularized, and the Stokes equations are solved to obtain a fluid velocity on the marker points due to the collection of regularized forces. The type of regularization determines the particular numerical method. For example, in the immersed boundary (IB) method of Peskin and collaborators, the force is regularized by smearing it onto a background grid on which the fluid equations are solved, and this velocity field is then interpolated back onto the marker points \cite{peskin1972flow, peskin2002acta}. In the special case when the spreading and interpolation are done with a Gaussian kernel, the method is referred to as a force coupling method (FCM) \cite{lomholt2003force, maxey2001localized}. For regularization and interpolation over the surface of a sphere, the force to velocity relationship (mobility matrix) can be computed analytically and is known as the Rotne-Prager-Yamakawa (RPY) tensor \cite{rpyOG, keavRPY, PSRPY}. Finally, the method of regularized Stokeslets describes the case when the Stokes equations are solved analytically for a given regularization function, and the resulting velocity field is evaluated directly on the marker points \rev{without an interpolation kernel}\footnote{\rev{This destroys the symmetry of hydrodynamic interactions, which is otherwise preserved in the IB, FCM, and RPY approaches. }} \cite{cortez2001method, cortez2005method}. 

All of these regularization methods have been used to model immersed rods, but generally with penalty terms to enforce inextensibility \cite{olson2013modeling, lim2008dynamics}. To our knowledge, only the recent approaches of Schoeller et al.\ \cite{keavRPY} and Jabbarzadeh and Fu \cite{inexRS} enforce inextensibility rigorously with Lagrange multipliers. In the case of \cite{inexRS}, the fiber is broken up into segments of regularized point forces, and each segment is updated via a rotation that preserves inextensibility exactly. Yet both \cite{keavRPY} and \cite{inexRS} suffer from the same pitfall as all regularization methods: when modeling slender fibers, the width/radius of the regularization function must be on the order of the fiber radius \cite{ttbring08}. Since the regularization width is also tied to the fiber discretization spacing (and, in the IB method, the fluid grid spacing), semi-flexible slender fibers must be discretized with many more points than would be necessary in a continuum, SBT-based approach. While this limitation can be partially overcome with adaptive mesh refinement \cite{griffith2007adaptive}, grid coarsening with local velocity correction \cite{maxian19}, and kernel-independent fast multipole methods (to accelerate many-body sums) \cite{rostami2016kernel}, the fact remains that to achieve controlled accuracy for dilute suspensions of many fibers, discretizing the fiber by a collection of marker points is much less efficient than treating it as a continuum with SBT. This is especially true for semi-flexible or stiff fibers, where the smooth fiber shapes are well represented in a spectral basis with rapidly decaying coefficients. 

Despite their limitations, regularized singularities are sometimes convenient to work with since they are nonsingular on the fiber centerline, and can therefore be easily evaluated there. A natural workaround to the regularization lengthscale issue is to take a continuum limit of many regularized point forces along the fiber \cite{cortez2012slender} \rev{or along segments of a fiber \cite{hall2019efficient, cortez2018regularized, walker2020regularised}}. For example, Walker et al.\ recently derived an SBT that uses regularized singularities along the fiber centerline and can be used for fibers of non-uniform cross section \cite{walkerregularised}. \rev{In still more recent work \cite{walker2020regularised}, they combined this theory with the regularized Stokeslet segments approach of Cortez \cite{cortez2018regularized} to yield a numerical method which is more efficient than that of regularized point forces, but still quadratic-complexity in the total number of segments. 

Since the shapes of biological filaments are smooth, a better approach is to represent the entire shape using an interpolating polynomial, rather than breaking into segments, and take a continuum limit of regularized singularities along the entire fiber. We show in Appendix \ref{sec:rpyappen} that applying this procedure to the RPY tensor yields a formula for the fiber velocity that is identical to SBT away from the fiber endpoints. For the reasons just listed, efficient simulation of many slender fibers requires a numerical method that can handle an SBT-type formulation for the fiber velocity. }

To our knowledge, Shelley and Ueda were the first to derive such a method and use it to simulate immersed slender fibers. By designing a numerical method around the analytical results of slender body theory, they reduced the complexity of the numerical computations from three dimensions to one \cite{shelley1996nonlocal, shelley2000stokesian}. Their formulation, however, relies on the filament being a closed loop, thus excluding many problems from biology, engineering, and physics where the filament ends are free. 

Tornberg and Shelley treated inextensible filaments with free ends using an SBT-based numerical method \cite{ts04}. In their approach, inextensibility is preserved by deriving an auxiliary \mbox{(integro-)differential} equation for the line tension in the filament, which acts as a Lagrange multiplier to preserve inextensibility. This method has since been used in applications with flexible (and sometimes fluctuating) filaments  \cite{manikantan2013subdiffusive, young2009hydrodynamic}, and was also extended to simulate falling rigid fibers, the novelty there being that many of the SBT-related integrals can be done analytically \cite{tornberg2006numerical}. More recently, Nazockdast et al.\ modified the approach of Tornberg and Shelley to make it feasible to simulate many-body cellular fiber assemblies. By replacing the second-order spatial discretization of Tornberg and Shelley with a spectral spatial discretization and utilizing a kernel-independent FMM to accelerate sums, Nazockdast et al.\ developed a parallel algorithm that makes it possible to simulate $\mathcal{O}(1000)$ fibers in linear time \cite{ehssan17}. 

Despite these recent advances, imposing inextensibility via a tension boundary value problem (BVP) leads to a number of drawbacks which are present in all of the prior SBT-based numerical methods. To begin, the line tension equation of \cite{ts04} involves multiplications of high-order (as high as four) derivatives of the fiber position function. This leads to severe aliasing problems and a loss of spatial accuracy in the spectral formulation \cite{ehssan17}. In addition, the ``inextensibility'' of the filaments is still subject to discretization error and requires inserting a penalty term into the line tension equation that reduces the discrete extensibility \cite{ts04}. For fibers tugged by cross-linkers or strong extensional flows, this penalty parameter will be large, introducing artificial stiffness into the problem. 

The primary focus of this paper is on a new formulation for inextensible filaments. In our approach, the fibers are evolved via a rotation of the tangent vector on the unit sphere, and the fiber positions are then obtained by integration. This approach is similar to that of \cite{inexRS}, but unique because we consider the fiber as a continuum, rather than a collection of discrete line segments. In this way, we maintain strict inextensibility of the fibers without introducing a penalty parameter. To close our formulation, we treat the force due to tension as a Lagrange multiplier and enforce the principle that the constraint forces do no work \cite{varibook}. We couple this advance with recent techniques \cite{barLud, tornquad} for efficient evaluation of nonlocal integrals appearing in SBT to develop a method that is both accurate and robust. \rev{In our spectrally accurate numerical method, we use, as in \cite{ehssan17}, Chebyshev polynomials of relatively low degree ($16-32$) to represent the fiber centerlines \emph{and} the force densities acting on them. This assumes that all of these quantities are smooth enough to be represented in the spectral basis.} \revp{This assumption fails for truly cylindrical fibers near the endpoints, and also when the fibers experience localized forces such as those due to steric repulsion, electrostatics, friction, and molecular motors. In our concluding Section\ \ref{sec:conclusion}, we discuss some possible ways to extend our method to account for these important biophysical forces and fiber shapes.} 

The rest of this paper is laid out as follows. We begin in Section \ref{sec:SBT} by introducing the necessary SBT equations for both local and nonlocal hydrodynamics. \rev{We modify the classical SBT formulation \cite{krub, johnson, gotz2001interactions} to regularize the local drag coefficient for cylindrical fibers and account for inter-fiber interactions through the RPY tensor.} In Section \ref{sec:inexfil}, we parameterize inextensible motions of the fiber as rotations of the unit tangent vector, thus strictly enforcing inextensibility. We then discuss how to determine the Lagrange multiplier forces for inextensibilty by imposing the principle of virtual work in a weak $L^2$ sense. Section \ref{sec:numerics} is devoted to numerical methods. We show how to incorporate a fast method for evaluating far-field hydrodynamic interactions (positively-split Ewald summation \cite{PSRPY}), and how to use specialized quadrature schemes for accurate evaluation of finite part and near-fiber quadratures. In Section \ref{sec:tint}, we design a semi-implicit, second-order temporal integrator that treats bending elasticity implicitly, yet for dilute systems only requires solving a block-diagonal linear system with a single evaluation of the nonlocal hydrodynamics per time step. For more concentrated systems, we use GMRES to solve a dense linear system, but show in Section \ref{sec:gresiters} that at most a few iterations are needed per time step to maintain stability. In the other numerical tests of Section \ref{sec:tests}, we also show how our ``weak formulation'' of inextensibility gives improved spatial accuracy over the traditional ``strong formulation'' of inextensibility in \cite{ts04, ehssan17}. In Section \ref{sec:CLs}, we study the rheology of a cross-linked network of filaments in oscillatory shear by introducing cross-linkers into the SBT formulation. Section \ref{sec:conclusion} gives our conclusions and discusses future work.

\section{Slender body theory \label{sec:SBT}}
We begin here by summarizing the slender body theories of \cite{krub, johnson, gotz2001interactions}, here following in particular Johnson \cite{johnson} and G\"otz \cite{gotz2001interactions}. These SBTs derive a global fluid velocity due to a single slender fiber, then evaluate this velocity asymptotically on the fiber surface to obtain a fiber evolution equation. It remains an open question, however, how to efficiently evaluate the fluid velocity generated by one filament on another filament. Here we formulate a modified treatment of these fiber-fiber interactions that is more physical and motivated by \rev{our observation in Appendix\ \ref{sec:rpyappen} that classical SBT can be formulated in terms of a line integral involving the RPY kernel}. In this section we will not consider any time dependence and look at the velocity of the filament and the Stokes fluid at a specific instant in time. We therefore omit explicit time dependence in our notation for simplicity. 

\subsection{Single filament}
We denote with $\V{X}(s)$ the position of the centerline of a filament, parameterized by arclength $s \in [0,L]$, where $L$ is the fiber length. The tangent vector is $\Xs(s)=\partial \V{X}/\partial s$ and has unit length. The fiber has physical radius $a(s) = r \rho(s)$, where $0 \leq \rho(s) \leq 1$, and slenderness ratio $\epsilon = r/L$. Let the force per unit length on the fiber centerline be denoted by $\V{f}(s)$ and the background flow (e.g. shear flow) at an arbitrary point in the fluid be denoted by $\V{u}_0(\V{x})$.

We recall the Stokeslet and doublet (Laplacian of the Stokeslet) kernels, which are the fundamental solutions to the Stokes equations for a point force and mass source dipole, respectively. If we center the kernels at $\V{x}_0$, take $\V{x}$ to be an arbitrary point in the fluid, and introduce $\V{R} = \V{x}-\V{x}_0$ with $\hat{\V{R}}=\V{R}/\norm{\V{R}}$, we have that 
\begin{equation}
\label{eq:Slet}
\Slet{\V{x},\V{x}_0} = \frac{\M{I}+\hat{\V{R}}\hat{\V{R}}}{\norm{\V{R}}} \qquad \text{ and } \qquad 
\Dlet{\V{x},\V{x}_0} = \frac{\M{I}-3\hat{\V{R}}\hat{\V{R}}}{\norm{\V{R}}^3}. 
\end{equation}

The idea of SBT is to introduce an ansatz for the flow field away from the fiber centerline of the form
\begin{align}
\label{eq:sbtsd}
\V{u}(\V{x}) - \V{u}_0(\V{x}) = & \EPMI \int_0^L \left(\Slet{\V{x},\V{X}(s)}+\beta(s)\Dlet{\V{x},\V{X}(s)}\right)\V{f}(s) \, ds\\[2 pt] 
\label{eq:Kdef}
\eqd & \EPMI \int_0^L \Knels{\V{x}}{\V{X}(s)}{\beta(s)}\V{f}(s) \, ds, 
\end{align}
where $\mu$ is the fluid viscosity. In Eq.\ \eqref{eq:Kdef}, we have defined a kernel $\Knels{\V{x}}{\V{x}_0}{\beta(s)}$ that is a combination of a Stokeslet and a doublet with strength $\beta$. Using the method of matched asymptotic expansions, the velocity integral\ \eqref{eq:sbtsd} can be computed analytically on the surface of the fiber to $\mathcal{O}(\epsilon)$ (see \cite{gotz2001interactions, koens2018boundary} for details on these integrals). The value of $\beta$ comes from imposing the boundary condition that the velocity on the fiber surface be constant to $\mathcal{O}(\epsilon)$; Mori et al.\ \cite{mori2018theoretical, morifree} refer to this as the ``fiber integrity condition.'' For cylindrical \cite{gotz2001interactions} or ellipsoidally-tapered \cite{johnson} filaments, this yields the solution for the velocity in the fluid as
\begin{equation}
\label{eq:sbt2}
\V{u}(\V{x}) - \V{u}_0(\V{x}) =\EPMI \int_0^L  \Knels{\V{x}}{\V{X}(s)}{\frac{(\epsilon L)^2}{2}}\V{f}(s) \, ds. 
\end{equation}
The fluid velocity $\V{u}(\V{x})$ in\ \eqref{eq:sbt2} does not apply inside of the fiber volume; in fact the kernel $\M{S}_D$ in\ \eqref{eq:sbt2} is not even defined on the fiber centerline. Physically, however, the velocity of the fiber centerline, which we denote with $\V{U}(s)$, should be equal to the average of $\V{u}(\V{x})$ around a ring cross section of the fiber with radius $a(s)$. Equivalently, since the function $\V{u}(\V{x})$ is constant on the cross section surface to $\mathcal{O}(\epsilon)$, averaging $\V{u}(\V{x})$ is equivalent to throwing out all terms in its expansion of $\mathcal{O}(\epsilon)$ or higher. \rev{Another approach, which is based on the Rotne-Prager-Yamakawa (RPY) kernel and matched asymptotics, is presented in Appendix\ \ref{sec:rpyappen}. The RPY tensor approximates the hydrodynamic interaction between two spheres of radius $\aRPY$ centered at $\V{x}$ and $\V{y}$ with the kernel \cite{rpyOG, wajnryb2013generalization, PSRPY}
\begin{gather}
\label{eq:rpyknel}
8 \pi \mu \KRPY{\V{x}}{\V{y}}{b}  = 
\begin{cases} 
\Knels{\V{x}}{\V{y}}{2\aRPY^2/3} & r > 2b\\
\left(\frac{4}{3b}-\frac{3r}{8b^2}\right)\M{I}+\frac{1}{8b^2 r}\V{r}\V{r}&  r < 2b
\end{cases},
\end{gather}
where $\V{r}=\V{x}-\V{y}$, $r=\norm{\V{r}}$, and we set the sphere radius (see\ \eqref{eq:rpyeps}) to
\begin{equation}
\label{eq:rpyradius}
b = \epsilon L \frac{e^{3/2}}{4} \approx 1.12 \epsilon L. 
\end{equation}
We express the velocity $\V{U}(s)$ on the fiber as a line integral of the RPY kernel
\begin{equation}
\label{eq:rpylineint}
\V{U}\left(s\right)-\V{u}_0\left(\V{X}(s)\right):= \EPMI \int_0^L  \KRPY{\V{X}(s)}{\V{X}\left(s^\prime \right)}{2b^2/3}\V{f}\left(s'\right) \, ds'.
\end{equation}

Using matched asymptotics to approximate\ \eqref{eq:rpylineint} for slender fibers (see Appendix\ \ref{sec:rpyappen}), we obtain the same result as classical SBT \cite{gotz2001interactions} }
\begin{gather}
\label{eq:onefib}
\V{U}(s)-\V{u}_0\left(\V{X}(s)\right)= \ML\left(\Xs(s); c(s)\right)\V{f}(s) + \left(\MFP\left[\V{X}(\cdot)\right]\V{f}(\cdot)\right)(s), \text{ where} \\[2 pt]
\label{eq:ML}
\ML(\Xs;c)= \EPMI \left(c (\M{I}+\Xs \Xs) +  (\M{I}-3\Xs \Xs)\right), \quad \text{and}\\[2 pt]
\label{eq:Mfp}
\left(\MFP\left[\V{X}\right]\V{f}\right)(s) =  \EPMI \int_{0}^L \left(\Slet{\V{X}(s),\V{X}(s')} \V{f}(s') -\left(\frac{\V{I}+\Xs(s)\Xs(s)}{|s-s'|}\right) \V{f}(s) \right) \, ds'. 
\end{gather}
Here $\ML$ is a $3 \times 3$ local drag matrix that gives the velocity contribution from the force density $\V{f}$ at points $\mathcal{O}(\epsilon)$ away from $\V{X}(s)$.  The integral operator $\MFP\left[\V{X}\right]$ gives the contribution from the rest of the fiber in the form of a finite part integral. The first term in the integrand is the Stokeslet, and the second term is the ``common'' part in the matched asymptotic expansion that comes from expansion of the Stokeslet around $s'=s$. Physically, the finite part integral gives the velocity contribution from forcing at points $\mathcal{O}(1)$ away from $\V{X}(s)$. Thus while both terms in the integrand are singular, their difference is finite \cite{ts04} (see also Section \ref{sec:MFP}). 

In the local drag matrix\ \eqref{eq:ML}, the leading order local drag coefficient is given by \cite{gotz2001interactions}
\begin{equation}
\label{eq:unmodc}
c(s) = \ln{\left(\frac{4s(L-s)}{a(s)^2}\right)}
\end{equation}
and is singular without proper decay of $a(s)$ at $s=0$ and $s=L$. Clearly, if $a(s)$ decays near the fiber endpoints as $2\epsilon\sqrt{s(L-s)}$, then the leading order coefficient\ \eqref{eq:unmodc} is finite at the fiber endpoints \cite{johnson}. This fact is the basis for a general assumption across the SBT literature that the filaments have ellipsoidal shape, so that in most studies $c(s)=-\ln(\epsilon^2)$ is constant for all $s$ \cite{ts04, ehssan17, das2018computing}. 

\begin{figure}
\centering
\includegraphics[width=0.45\textwidth]{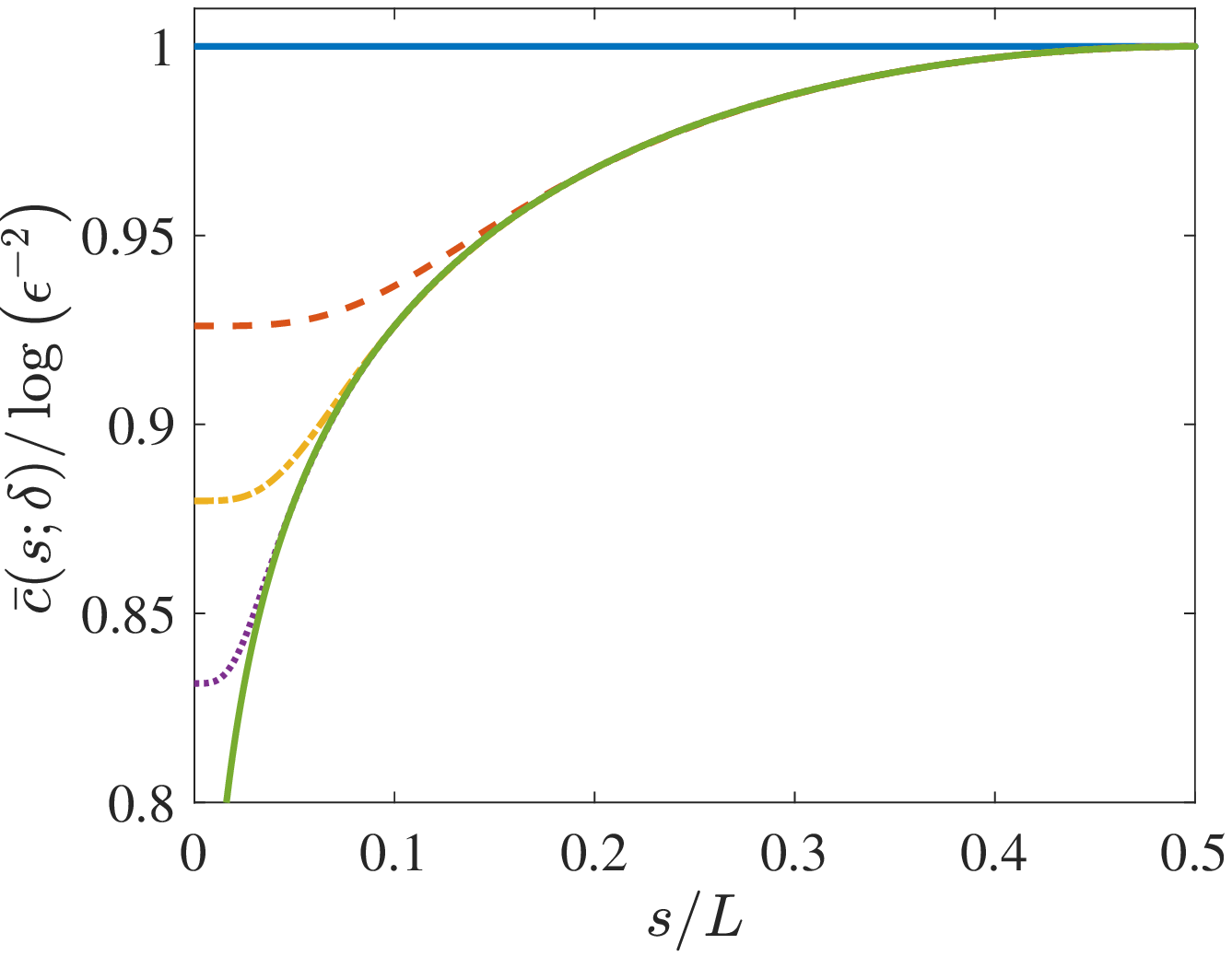}
\includegraphics[width=0.45\textwidth]{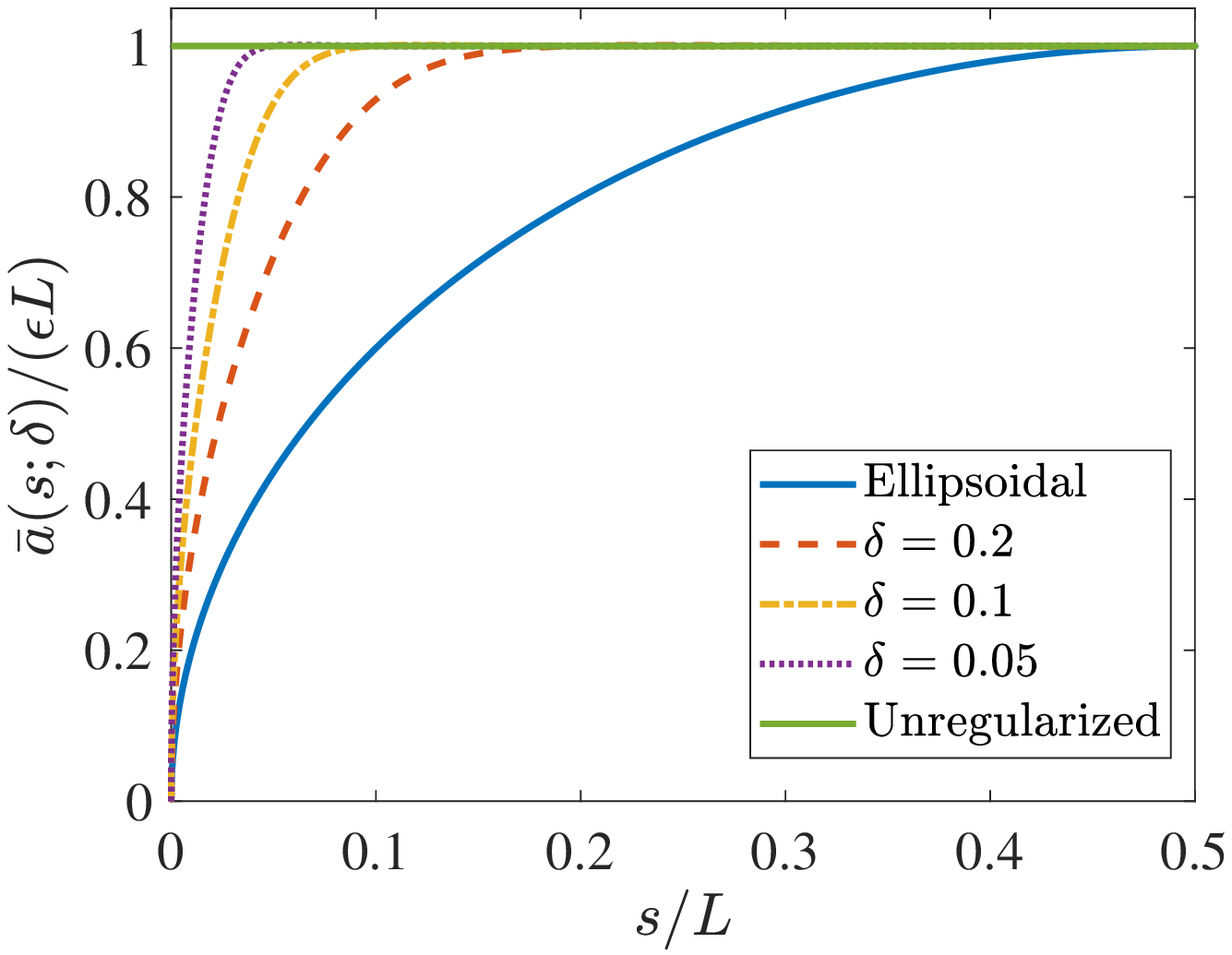}
\caption{\label{fig:regc} \rev{Regularized drag coefficients and effective radius functions. (Left:) The regularized local drag coefficients $\bar{c}(s;\delta)$ for various $\delta$. (Right:) The effective fiber radius $\bar{a}(s;\delta)$ as defined in\ \eqref{eq:creg}. We show only one side of the fiber, with the other being the former's mirror image about $s/L=0.5$.}}
\end{figure}

Actin filaments are best modeled as cylinders with constant radius, so that $a(s) = r = \epsilon L$ on $s \in [0,L]$. In this case, the coefficient\ \eqref{eq:unmodc} becomes singular at the filament ends, and so \rev{we modify the local drag coefficient by effectively tapering the fiber radius over a distance $\sim\delta L$ near the endpoints. Specifically, to regularize $c(s)$, we set $\eta = 2s/L-1$, so that $\eta \in [-1,1]$, and compute a weight function
\begin{equation}
w(s;\delta) = \text{tanh}\left(\frac{\eta(s)+1}{\delta}\right)-\text{tanh}\left(\frac{\eta(s)-1}{\delta}\right)-1, 
\end{equation}
which is 1 near the fiber center ($\eta=0$) and zero at the fiber ends ($\eta = \pm 1$). We then assign to each $s$ a regularized fiber centerline coordinate by
\begin{equation}
\bar{s}(s;\delta) = w(s;\delta) s+ \left(1-w(s;\delta)^2\right)\frac{\delta L}{2}
\end{equation}
on $0 \leq s \leq L/2$, with the corresponding reflection for $s > L/2$. The regularized coefficient for a given $\delta$ is then given by
\begin{equation}
\label{eq:creg}
\bar{c}(s;\delta):=c(\bar{s}(s;\delta)) = \ln{\left(\frac{4\bar{s}(L-\bar{s})}{(\epsilon L)^2}\right)}:= \ln{\left(\frac{4s(L-s)}{\bar{a}(s;\delta)^2}\right)}. 
\end{equation}
Figure\ \ref{fig:regc} shows how the choice of $\delta$ impacts the local drag coefficient $\bar{c}(s)$ and effective radius function $\bar{a}(s)$. We see that the fibers are cylindrical ($\bar{a}(s)=\epsilon L$) for $s/L \gtrsim \delta$, while for $s/L \lesssim \delta$ the effective radius smoothly decays to zero, as it would for an ellipsoidal fiber. Larger values of $\delta$ yield smoother radius functions and smoother local drag coefficients. Throughout this paper, we will use $\delta=0.1$, unless otherwise stated.}

\subsection{Multiple filaments \label{sec:RPYSBT}}
It remains to include in the fiber centerline velocity\ \eqref{eq:onefib} the perturbed flow due to other filaments, i.e., to account for hydrodynamic interactions between fibers. We require more involved notation in this case to distinguish between fibers. In general, we use the symbol $\V{X}$ in an equation whenever it is localized to a single fiber and is the same on all fibers. For example, $\Lop{H}[\V{X}]=0$ implies that $\Lop{H}$ is a functional of a single fiber's position and is zero on every fiber individually. When multiple fibers are involved, we index the $i$th fiber by the superscript $\ind{\V{X}}{i}$. For example, $\Lop{H}\left[\ind{\V{X}}{i},\ind{\V{X}}{j}\right]\ind{\V{f}}{i}$ implies that $\Lop{H}$ is a functional of a pair of fiber positions which acts on the force density on fiber $i$. Whenever an equation applies to one fiber, but involves all other fibers, we will use the notation $\ind{\V{X}}{i}$ for the single fiber and $\Xcol$ to refer to the collection of all fibers. Our use of $\V{X}$ to refer to both a general fiber and collection of fibers is a slight abuse of notation, but the meaning should be clear from the context, specifically whether an equation involves a single fiber or multiple fibers. 

The simplest approach for including hydrodynamic interactions in SBT is the one taken by Tornberg and Shelley \cite{ts04}, in which the fluid velocity due to one fiber\ \eqref{eq:sbt2} is simply evaluated on the centerline of the other fibers. Nazockdast et al.\ also adopted this, except they dropped the doublet term completely and included only the Stokeslet term. \rev{Inspired by the fact that the classical SBT for a single fiber\ \eqref{eq:onefib} can be reformulated as a regularized singularity method using the RPY line integral\ \eqref{eq:rpylineint}, we use the same line integral for interaction with other fibers. In\ \eqref{eq:dco}, we substitute the choice of RPY radius\ \eqref{eq:rpyradius} into the RPY kernel\ \eqref{eq:rpyknel} to obtain the inter-fiber interaction velocity
\begin{align}
\nonumber
\ind{\V{v}}{j}\left(\ind{\V{X}}{i}(s)\right) & \eqd \EPMI \int_0^{L} \Knel{\ind{\V{X}}{i}(s),\ind{\V{X}}{j}(s'), \dco \left(\epsilon L\right)^2}\ind{\V{f}}{j}(s') \, ds' \\[2 pt]
\label{eq:sbtother}
& \eqd \left(\MJF\left[\ind{\V{X}}{j}(\cdot)\right] \ind{\V{f}}{j}(\cdot)\right)\left(\ind{\V{X}}{i}(s)\right)
\end{align}}
Here we have again defined a linear integral operator $\MJF$ which acts on $\ind{\V{f}}{j}$ to give the velocity on the centerline of filament $i$ solely due to filament $j$. More generally, we will denote by $\ind{\V{v}}{j}\left(\V{x}\right)$ the velocity induced by filament $j$ at any point $\V{x}$ on the centerline of any other filament. This ``filament interaction velocity'' differs from the slender body fluid velocity\ \eqref{eq:sbt2} \rev{in that the coefficient of the dipole term is $e^3/24 \approx 0.84$ instead of $1/2$ (see\ \eqref{eq:dco})}.  Since $\ind{\V{v}}{j}\left(\V{x}\right)$ can be evaluated everywhere in the fluid, we use a lowercase letter to denote it and distinguish it from the \textit{fluid} velocity $\ind{\V{u}}{j}$ induced by fiber $j$, which is given in\ \eqref{eq:sbt2}.

\rev{We have used a constant radius $b$ when defining the inter-fiber mobility\ \eqref{eq:sbtother} using the RPY kernel. While our use of the regularized local drag coefficient\ \eqref{eq:creg} implies that the fibers are tapered at the endpoints, and therefore that $\aRPY$ decays to zero at the endpoints, we will use a constant radius for nonlocal hydrodynamics. This assumption, together with the assumption that all of the fibers have the same maximum radius $ \epsilon L$, can both be relaxed, as the RPY kernels derived in \cite{RPY_Shear_Wall} and associated fast methods \cite{SE_Multiblob_SD} can be used to generalize\ \eqref{eq:sbtother} to tapered filaments or fibers with different radii.}

 
Summing the interaction kernel\ \eqref{eq:sbtother} over filaments $j \neq i$ and adding the terms from local drag, we get a slender body theory for the velocity $\ind{\V{U}}{i}(s)$ at position $s$ on filament $i$, 
\begin{gather}
\label{eq:fibevcont}
\ind{\V{U}}{i}(s) -\V{u}_0\left(\ind{\V{X}}{i}(s)\right)=  \ML \left(\iXs{i}(s); c(s)\right)\ind{\V{f}}{i}(s)\\[2 pt] \nonumber + \left(\MFP\left[\ind{\V{X}}{i}(\cdot)\right]\ind{\V{f}}{i}(\cdot)\right)(s) + \sum_{j \neq i}\left(\MJF\left[\ind{\V{X}}{j}(\cdot)\right]\ind{\V{f}}{j}(\cdot)\right)\left(\ind{\V{X}}{i}(s)\right), 
\end{gather}
where $\ML$ is defined in Eq.\ \eqref{eq:ML}, $\MFP$ is defined in Eq.\ \eqref{eq:Mfp}, and $\MJF$ is defined in Eq.\ \eqref{eq:sbtother}. Because $\ind{\V{U}}{i}$ is only defined on the centerline of filament $i$ (and not everywhere in the fluid), we denote it with a capital letter. 

For a single fiber $i$, we will write the mobility\ \eqref{eq:fibevcont} abstractly as
\begin{equation}
\label{eq:ijabstract}
\ind{\V{U}}{i}(s) -\V{u}_0\left(\ind{\V{X}}{i}(s)\right) = \sum_j \left(\Lop{M}_{ij}\left[\V{X}\right] \ind{\V{f}}{j}\right)(s), 
\end{equation}
where the mobility operator $\Lop{M}\left[\V{X}\right]$ is a functional of the positions of all fibers. We will compactly write the velocity\ \eqref{eq:ijabstract} for the entire collection of fibers as 
\begin{equation}
\label{eq:mobeqn}
 \V{U}-\V{u}_0\left(\V{X}\right) \eqd \Lop{M}\left[\V{X}\right]\V{f}. 
\end{equation}
This mobility equation can be closed by defining a constitutive equation for the fiber force densities $\V{f}$, which we do next.

\section{Inextensible filaments \label{sec:inexfil}}
In this paper, we consider inextensible filaments $\ind{\V{X}}{i}(s,t)$ which can bend, but not stretch, as they evolve in time. We assume the fibers are in a constant twist-free equilibrium, since in the absence of externally-applied or internally-generated torques the timescale of twist relaxation is $\mathcal{O}(\epsilon^{-2})$ faster than bending \cite{powers2010dynamics}. 

At every instant in time, each fiber resists bending with bending force density (per unit length) $\V{f}^\kappa\left[\V{X}\right]$. Inextensibility can be enforced by introducing Lagrange multiplier force densities on each fiber $\ind{\V{\lambda}}{i}(s,t)$, where we will again write $\V{\lambda}=\left\{\ind{\V{\lambda}}{i}\right\}$ whenever we refer to the collection of Lagrange multipliers on all fibers. Thus the PDE that we need to solve on every fiber is given by (using the abstract notation of\ \eqref{eq:mobeqn}), 
\begin{equation}
\label{eq:fibPDE}
\frac{\partial \V{X}}{\partial t} - \V{u}_0\left(\V{X},t\right) = \Lop{M}\left[\V{X}\right] \left(\V{f}^\kappa\left[\V{X}\right]+\V{\lambda}\right),
\end{equation}
where the mobility operator $\Lop{M}$ is defined in\ \eqref{eq:fibevcont} and the background flow function $\V{u}_0$ can in general vary in time. 
The fibers are also constrained to be inextensible, so that for every fiber
\begin{equation}
\label{eq:inex}
\Xs(s,t) \cdot \Xs(s,t)=1,
\end{equation}
for all $s$ and $t$. We still need to specify boundary conditions for the evolution equation\ \eqref{eq:fibPDE} and additional conditions on $\V{\lambda}$ to make the solution unique, as we explain shortly. 

\subsection{Bending elasticity}
For fiber mechanics, we use the Euler beam model, in which the bending force density on every fiber is given by
\begin{equation}
\label{eq:bforce}
    \V{f}^\kappa\left[\V{X}\right]=-\kappa\V{X}_{ssss}\eqd \Lop{F}\V{X},
\end{equation}
where the constant linear operator $\Lop{F}$ gives $ \V{f}^\kappa$ taking into account the ``free fiber'' boundary conditions \cite{ts04} 
\begin{gather}
\label{eq:frBCs}
\V{X}_{ss}\left(s=0,t\right)=\V{X}_{sss}\left(s=0,t\right)=\V{0},\\[2 pt]
\nonumber
\V{X}_{ss}\left(s=L,t\right)=\V{X}_{sss}\left(s=L,t\right)=\V{0}.
\end{gather}
Again, because the boundary conditions\ \eqref{eq:frBCs} apply to every fiber without dependence on other fibers, we use the notation $\V{X}_{ss}$ and $\V{X}_{sss}$ to refer to the arclength derivatives along the fiber and drop the superscript $(i)$. 

It is easy to see that the boundary conditions\ \eqref{eq:frBCs} cause the total force and torque on every fiber due to $\V{f}^\kappa$ to be zero, 
\begin{gather}
\label{eq:totfe}
\int_{0}^{L} \V{f}^\kappa \, ds = -\kappa\V{X}_{sss}\Big \rvert^{L}_{0} = \V{0}, \quad \text{and}\\[2 pt]
\left(\int_{0}^{L} \V{f}^\kappa \times \V{X} \, ds\right)^\ell = -\kappa\int_{0}^{L} \left(X^j X^k_{ssss} - X^k X^j_{ssss}\right) \, ds = -\kappa\int_{0}^{L} \left(X^j_{ss}  X^k_{ss} -X^k_{ss} X^j _{ss}\right) = 0. 
\end{gather}
Here the set of superscripts $(j,k,\ell)$ denote vector components and are a cyclic permutation of $(1,2,3)$. In the torque equation, the free fiber boundary conditions lead to the cancellation of boundary terms that arise in integration by parts.

\subsection{Traditional formulation of inextensibility}
In the traditional formulation of inextensibility \cite{ts04}, the inextensibility constraint\ \eqref{eq:inex} is differentiated with respect to time. Then, $s$ and $t$ derivatives are interchanged to yield
\begin{gather}
\label{eq:inexdt}
\left(\ddt{\V{X}}\right)_s \cdot \Xs = 0. 
\end{gather}
In \cite{ts04}, the system was closed by substituting the mobility equation\ \eqref{eq:fibPDE} into the differentiated inextensibility constraint\ \eqref{eq:inexdt}. On each fiber, Tornberg and Shelley then assume that $\ind{\V{\lambda}}{i}=\left(\ind{T}{i}\ind{\Xs}{i}\right)_s$, where $\ind{T}{i}(s,t)$ is an unknown scalar tension \cite{ts04}. This results in the \textit{line tension equation},
\begin{equation}
\label{eq:lineT}
\frac{\partial}{\partial s}\left(\V{u}_0\left(\ind{\V{X}}{i}\right)+\sum_{j}\Lop{M}_{ij}\left[\V{X}\right]\left(\Lop{F}\ind{\V{X}}{j}+\left(\ind{T}{j}\iXs{j}\right)_s\right)\right) \cdot \iXs{i}= 0.
\end{equation}
which holds for each fiber $i$ \cite{ts04}. While the second-order BVP\ \eqref{eq:lineT} is linear in $T$, it is highly nonlinear in $\V{X}$, since the operation $\Lop{F}\V{X}$ gives fourth derivatives of $\V{X}$. Even in the absence of any nonlocal hydrodynamic interactions (i.e. if $\Lop{M}=\ML$) and zero background flow ($\V{u}_0=\V{0}$), the line tension equation still has terms of the form $\V{X}_{sss}\cdot \V{X}_{sss}$ (see \cite[~Eq.\ (13)]{ts04}), which lead to aliasing errors in spectral numerical methods. Because the line tension equation\ \eqref{eq:lineT} enforces inextensibility pointwise along the fiber, we refer to it as a \textit{strong formulation of inextensibility}. 

\subsection{Kinematics of inextensible fibers \label{sec:geo}}
In our approach, we evolve the tangent vector $\Xs(s,t)$, rather than $\V{X}(s,t)=\V{X}(0,t)+\int_0^{s} \Xs\left(s',t\right) \, ds'$. Considering the evolution of $\Xs(s,t)$, the differentiated inextensibility constraint\ \eqref{eq:inexdt} implies that, for every fiber, 
\begin{equation}
\label{eq:omegadef}
\frac{\partial \Xs}{\partial t}(s,t) = \V{\Omega}(s, t) \times \Xs(s,t), 
\end{equation}
i.e., that the fiber evolution can be thought of as rotations of $\Xs$ on the unit sphere. 

At each fiber point, we uniquely define an orthonormal coordinate system using spherical angles $\theta(s,t)$ and $\phi(s,t)$. We represent the unit tangent vector $\Xs(s,t)$ as
\begin{equation}
\label{eq:Xsangle}
\Xs(s,t)= \begin{pmatrix} \cos{\theta} \cos{\phi}\\[2 pt] \sin{\theta} \cos{\phi} \\[2 pt] \sin{\phi} \end{pmatrix}, 
\end{equation}
where we define $\theta$ to be single-valued at $\phi=\pi/2$ by setting $\theta \left(\phi=\pm \pi/2\right)=0$. A choice of normal vectors that are always orthonormal to $\Xs$ on the unit sphere is
\begin{equation}
\label{eq:nangles}
\V{n}_1 =  \begin{pmatrix} -\sin{\theta}\\[2 pt] \cos{\theta}\\[2 pt]0 \end{pmatrix}, \qquad \V{n}_2 =  \begin{pmatrix} -\cos{\theta} \sin{\phi}\\[2 pt] -\sin{\theta} \sin{\phi} \\[2 pt] \cos{\phi} \end{pmatrix}. 
\end{equation}
Because $\V{n}_1$ and $\V{n}_2$ can be determined uniquely from $\Xs$, we denote them henceforth with $\V{n}_j\left(\Xs(s,t)\right)$, for $j=1, 2$. Since $\theta$ is single-valued at $\phi=\pi/2$, each component of the orthonormal coordinate system $(\Xs,\V{n}_1,\V{n}_2)$ is a smooth function of $\V{X}$ when $\Xs$ is smooth. Importantly, our method does not depend on the particular choice of normal vectors\ \eqref{eq:nangles}; any choice that gives smooth $\V{n}_1$ and $\V{n}_2$ for a smooth $\V{X}$ is equally acceptable. For example, the Frenet or Bishop frames could be used \cite{langer1996lagrangian}. 

Because $\Xs \times \Xs=\V{0}$, and since we are not considering twist, $\V{\Omega}(s,t)$ can be restricted to linear combinations of $\V{n}_1$ and $\V{n_2}$. We let 
\begin{equation}
\label{eq:omdef}
\V{\Omega}(s,t) \eqd \V{\Omega}\left(\Xs(s,t),\V{g}(s,t)\right) \eqd g_1(s,t)\V{n}_2\left(\Xs(s,t)\right)-g_2(s,t)\V{n}_1\left(\Xs(s,t)\right),  
\end{equation}
where $g_1(s,t)$ and $g_2(s,t)$ are two specific unknown functions and $\V{g}=\{g_1,g_2\}$. Equation \eqref{eq:omdef} implies that, by the right-handedness of the coordinate system $(\Xs, \V{n}_1, \V{n}_2)$, 
\begin{equation}
\label{eq:Xsupdate}
\ddt{\Xs} = \V{\Omega} \times \Xs = g_1\V{n}_1 + g_2\V{n}_2. 
\end{equation}
Any inextensible velocity of the fiber centerline can now be written in the form 
\begin{equation}
\label{eq:velfib}
\V{U}(s,t)=\ddt{\V{X}}(s,t) = \Ubar(t)+\int_0^s  \sum_{j=1}^2 g_j(s',t)\V{n}_j\left(\Xs(s',t)\right) \,ds' , 
\end{equation}
where $\Ubar(t)=\partial \V{X}/\partial t(s=0,t)$ is a rigid body translation. 


\subsection{Principle of virtual work}
The kinematic formulation of Section \ref{sec:geo} can still be used to solve for the line tensions and fiber velocities. In particular, by substituting the inextensible velocity\ \eqref{eq:velfib} into the left hand side of the evolution equation\ \eqref{eq:fibPDE} and setting $\V{\lambda}=\left(T\Xs\right)_s$, a PDE results with unknowns $g_1$, $g_2$, $\Ubar,$ and $T$. We choose to close our formulation differently, in the process eliminating the need to solve for tension explicitly. 

On every fiber, the principle of virtual work states that the constraint forces $\V{\lambda}$ do no work for any choice of $g_1, g_2,$ and $\Ubar$ \cite{varibook}. Because this constraint holds for all time, for simplicity we drop for the moment the explicit dependence on $t$ in the notation. To impose the principle of virtual work, we use the $L^2$ inner product to compute the total power dissipated in the fluid from $\V{\lambda}$, 
\begin{align}
\label{eq:15}
\mathcal{P}& =\Bigg{\langle} \V{\lambda},\ddt{\V{X}} \Bigg{\rangle} = \int_0^L ds' \left(\Ubar+\int_0^{s'} \left(g_1(s) \V{n}_1\left(\Xs(s)\right) + g_2(s)\V{n}_2\left(\Xs(s)\right)\right)\,ds\right) \cdot \V{\lambda}\left(s'\right).
\end{align}
Changing integration variables, we can rewrite this as 
\begin{align}
\label{eq:16}
\mathcal{P}& = \Ubar \cdot \int_0^L \V{\lambda}\left(s'\right) \, ds' +\int _0^L ds \int_{s}^L \left(g_1(s)\V{n}_1\left(\Xs(s)\right) + g_2(s) \V{n}_2\left(\Xs(s)\right)\right) \cdot \V{\lambda}(s') \, ds'\\[2 pt]
\label{eq:Kstarcont}
& = \Ubar \cdot \int_0^L \V{\lambda}\left(s'\right) \, ds' + \int_0^L \left(g_1(s)\V{n}_1\left(\Xs(s)\right)+ g_2(s) \V{n}_2\left(\Xs(s)\right)\right) \cdot  \left(\int_{s}^L  \V{\lambda}(s') \, ds'\right) \, ds=0. 
\end{align}

\subsubsection{Pointwise formulation}
Since the principle of virtual work\ \eqref{eq:Kstarcont} must hold for any inextensible motion, it must hold for all $\Ubar$ and all sufficiently smooth $g_1$ and $g_2$. Therefore, we must have, for all $s$, 
\begin{equation}
\label{eq:noworkcont}
\begin{pmatrix} \left(\int_s^L \V{\lambda}(s')\, ds'\right) \cdot \V{n}_1\left(\Xs(s)\right)\\[2 pt] \left(\int_s^L \V{\lambda}(s')\, ds'\right) \cdot\V{n}_2\left(\Xs(s)\right)\\[2 pt] \int_0^L \V{\lambda}(s') ds' \end{pmatrix} = \begin{pmatrix} 0 \\[2 pt] 0\\[2 pt] \bm{0}\end{pmatrix}. 
\end{equation}
The first and second components of the constraints\ \eqref{eq:noworkcont} taken together tell us that $\int_{s}^L  \V{\lambda}\left(s'\right) \, ds'$ is orthogonal to both normal vectors. Therefore, $\int_{s}^L  \V{\lambda}\left(s'\right) \, ds'$ is in the direction of $\Xs(s)$ and can be written as
\begin{equation}
\label{eq:tsalm}
\int_{s}^L  \bm{\lambda}\left(s'\right) \, ds' = -T(s)\Xs(s),
\end{equation}
for some scalar function $T(s)$ with $T(s=L)=0$. This gives
\begin{equation}
\label{eq:lamval}
\V{\lambda}(s) = \left(T(s)\Xs(s)\right)_s, 
\end{equation}
which is the form assumed in Tornberg and Shelley \cite{ts04}. Thus our derivation shows that the form of $\V{\lambda}$ taken in \cite{ts04} is equivalent to the principle that the constraint forces perform no virtual work, if the work is given by the standard $L^2$ inner product \cite{varibook}. 

Now, returning to the third of the constraints\ \eqref{eq:noworkcont}, $\int_0^L \V{\lambda}(s) \, ds=\V{0}$, and substituting the derived form of $\V{\lambda}$ in\ \eqref{eq:lamval}, we obtain 
\begin{equation}
\label{eq:BCT}
T(L)\Xs(L) - T(0)\Xs(0)=\V{0}. 
\end{equation}
Since $T(L)=0$, Eq.\ \eqref{eq:BCT} implies that $T(0)=0$ as well, since neither of the tangent vectors is identically 0. So we obtain $T(0)=T(L)=0$, which is exactly the boundary condition for the line tension equation in \cite{ts04}. The form of $\V{\lambda}=\left(T\Xs\right)_s$ and the tension boundary conditions imply that the total torque induced by the constraint forces is zero in continuum, $\int_0^L \V{X}(s) \times \V{\lambda}(s) \, ds = \V{0}$. 

In this sense, the constraint equation \eqref{eq:15} is equivalent to the line tension equation used in prior work \cite{ts04}. Because we showed the equivalence by enforcing constraint\ \eqref{eq:15} for every choice of $g_1(s)$ and $g_2(s)$, we refer to the inextensibility constraint\ \eqref{eq:15} as a \textit{weak formulation of inextensibility}. In the next section, we choose a suitable basis for $g_1(s)$ and $g_2(s)$ to obtain a linear system of equations instead of the pointwise constraint\ \eqref{eq:noworkcont}.

\subsubsection{$L^2$ weak formulation}
\label{sec:numinex}
In this section, we introduce an $L^2$ weak formulation that is suitable for a numerical discretization of the weak inextensibility constraint\ \eqref{eq:15}. The key idea is to expand the unknown functions $g_1(s)$ and $g_2(s)$ as, 
\begin{equation}
\label{eq:basis}
g_j(s) = \sum_k \alpha_{jk} T_k(s), \quad \text{for $j=1, 2$}, 
\end{equation}
where $T_k(s)$ are sufficiently smooth scalar-valued basis functions for $L^2:[0,L]$. Substituting the basis function expansion\ \eqref{eq:basis} into the inextensible velocity\ \eqref{eq:velfib}, we obtain
\begin{equation}
\label{eq:du}
\V{U}(s) = \ddt{\V{X}}(s) =\Ubar + \int_0^s \sum_{j=1}^2\sum_k \alpha_{jk} T_k\left(s'\right) \V{n}_j\left(\Xs(s')\right) \, ds' \eqd (\Lop{K}\left[\V{X}(\cdot)\right]\V{\alpha})(s),
\end{equation}
where we have defined a linear operator $\Lop{K}\left[\V{X}\right]$ on every fiber that acts on $\V{\alpha}=\left(\alpha_{jk},\, \Ubar\right)$ to give an inextensible velocity on the filament centerline (i.e., $\V{\alpha}$ parameterizes the space of inextensible fiber motions). Note the functional dependence of $\Lop{K}$ on $\V{X}$ since $\Lop{K}$ involves the normal vectors $\V{n}_1$ and $\V{n}_2$. Substituting the inextensible velocity\ \eqref{eq:du} into the dynamical equation\ \eqref{eq:fibPDE}, we obtain
\begin{equation}
\label{eq:veleqn}
\Lop{K}\left[\V{X}\right]\V{\alpha} = \V{u}_0\left(\V{X}\right) + \Lop{M}\left[\V{X}\right]\left(\Lop{F}\V{X} +\V{\lambda}\right). 
\end{equation}

This constrained dynamical equation is supplemented by enforcing the principle of virtual work\ \eqref{eq:noworkcont} in an $L^2$ weak sense. We begin by substituting the representation of $g_j$ in\ \eqref{eq:basis} into the power equation \eqref{eq:15} to obtain, for every fiber, 
\begin{align}
\label{eq:Kstardef}
\mathcal{P} & = \bigg{\langle} \V{\lambda},\Lop{K}[\V{X}]\V{\alpha}\bigg{\rangle}\eqd \bigg{\langle} \Lop{K}^*[\V{X}]\V{\lambda}, \V{\alpha}\bigg{\rangle} \\[2 pt]
\nonumber
& =\Ubar \cdot \int_0^L \V{\lambda}(s) \, ds + \int_0^L \left(\int_0^{s} \sum_{j=1}^2\sum_{k} \alpha_{jk} T_k\left(s'\right) \V{n}_j\left(\Xs\left(s'\right)\right)\, ds'\right)\cdot \V{\lambda}\left(s\right) \, ds = 0, 
\end{align}
where we have defined $\Lop{K}^*$ as the $L^2$ adjoint of $\Lop{K}$. Since the power from the constraint forces must be zero for any inextensible motion (any $\V{\alpha}$), each term of the constraint\ \eqref{eq:Kstardef} must be zero. This gives the set of constraints on every fiber
\begin{equation}
\label{eq:noworkcontL2}
\Lop{K}^*[\V{X}]\V{\lambda}:=\begin{pmatrix} \int_0^L \left(\int_0^{s} T_k(s') \V{n}_1\left(\Xs(s')\right) \, ds'\right)\cdot \bm{\lambda}(s) \, ds\\[2 pt] \int_0^L \left(\int_0^{s} T_k(s') \V{n}_2\left(\Xs(s')\right)\, ds'\right)\cdot \bm{\lambda}(s) \, ds\\[2 pt] \int_0^L \bm{\lambda}(s) \, ds \end{pmatrix} = \begin{pmatrix} 0 \\[2 pt] 0\\[2 pt] \V{0}\end{pmatrix}, 
\end{equation}
where the first two constraints hold for all $k$ and the last constraint holds for each of the three Cartesian directions.

\subsection{Summary of dynamical equations}
In our abstract notation, the evolution of the fiber system can be obtained by solving the following system for $\V{\alpha}(t)=\left\{\ind{\alpha}{i}_{jk}(t), \ind{\Ubar \, }{i}(t)\right\}$ and $\V{\lambda} = \left\{\ind{\V{\lambda}}{i}(s,t)\right\}$, 
 \begin{gather}
 \label{eq:abssys}
 \ddt{\V{X}}=\Lop{K}\left[\V{X}\right]\V{\alpha} = \V{u}_0\left(\V{X},t\right) + \Lop{M}\left[\V{X}\right]\left(\Lop{F}\V{X}+\V{\lambda}\right)\\[2 pt]
 \label{eq:nowork}
 \Lop{K}^*\left[\V{X}\right]\V{\lambda} = \V{0}. 
 \end{gather}
 The first equation\ \eqref{eq:abssys} is the mobility equation. The left hand side is the velocity of a fiber centerline, restricted to the space of inextensible motions via the operator $\Lop{K}$ defined in\ \eqref{eq:du}. The right hand side involves all fiber positions and force densities because of hydrodynamic interactions. The second equation \eqref{eq:nowork} is the principle of virtual work and applies on each fiber separately, $\Lop{K}^*\left[\ind{\V{X}}{i}\right]\ind{\V{\lambda}}{i}=\V{0}$.  
 
On a single fiber $i$, the mobility equation\ \eqref{eq:abssys} takes the explicit form
 \begin{gather}
 \label{eq:Kalphexpl}
 \ddt{\ind{\V{X}}{i}}\rev{(s,t)}= \ind{\Ubar \, }{i}(t) + \int_0^s \sum_{j=1}^2\sum_k \ind{\alpha}{i}_{jk}(t) T_k\left(s'\right) \V{n}_j\left(\Xs(s',t)\right) \, ds' =\rev{\V{u}_0\left(\ind{\V{X}}{i}(s),t\right)}  \\[2 pt]
 \label{eq:MLexpl}
 +\EPMI \Bigg{(}\left(\rev{\bar{c}(s;\delta)} \left(\M{I}+\ind{\Xs}{i}(s,t) \ind{\Xs}{i}(s,t)\right) +  \left(\M{I}-3\ind{\Xs}{i}(s,t) \ind{\Xs}{i}(s,t)\right)\right)\ind{\V{f}}{i}(s,t)\\[2 pt]
 \label{eq:MFPexpl}
+ \int_{0}^{L} \left(\Slet{\ind{\V{X}}{i}(s,t),\ind{\V{X}}{i}\left(s',t\right)} \ind{\V{f}}{i}\left(s',t\right) -\left(\frac{\V{I}+\ind{\Xs}{i}(s,t)\ind{\Xs}{i}(s,t)}{|s-s'|}\right) \ind{\V{f}}{i}(s,t) \right) \, ds'\\[2 pt]
\label{eq:MJFexpl}
+ \sum_{j \neq i} \int_0^{L} \left(\Slet{\ind{\V{X}}{i}(s,t),\ind{\V{X}}{j}\left(s',t\right)}+\rev{\dco} \left(\epsilon L\right)^2\Dlet{\ind{\V{X}}{i}(s,t),\ind{\V{X}}{j}\left(s',t\right)}\right) \ind{\V{f}}{j}\left(s',t \right) \, ds'\Bigg{)}, 
 \end{gather}
with $\ind{\V{f}}{i}(s,t) = \Lop{F}\ind{\V{X}}{i}(s,t)+\ind{\V{\lambda}}{i}(s,t)$. The Stokeslet and doublet kernels $\M{S}$ and $\M{D}$ are defined in\ \eqref{eq:Slet}, and the \rev{local drag coefficient $\bar{c}$ is regularized at the endpoints as defined in\ \eqref{eq:creg}}. The mobility equation for fiber $i$ is supplemented by the principle of virtual work\ \eqref{eq:nowork} which is localized to fiber $i$ and takes the explicit form
\begin{gather}
\label{eq:noworkcontL2expl}
 \int_0^{L} \left(\int_0^{s} T_k(s') \V{n}_j\left(\iXs{i}(s',t)\right) \, ds'\right)\cdot \ind{\V{\lambda}}{i}(s,t) \, ds = 0, \quad \forall k \text{ and } j=1, 2,\\[2 pt] \int_0^{L} \ind{\V{\lambda}}{i}(s,t) \, ds =\V{0}. 
\end{gather}

\section{Numerical Methods \label{sec:numerics}}
Our goal in this section is to write the evolution equations\ \eqref{eq:abssys} and \eqref{eq:nowork} in the form of a block-matrix saddle point system. We will replace the operators with matrices and the position functions $\ind{\V{X}}{i}(s)$ with discrete vectors of collocation points $\ind{\V{X}}{i}$. The fiber evolution is then given by
\begin{equation}
\ddt{\ind{\V{X}}{i}}(t) = \M{K}\left(\ind{\V{X}}{i}\right)\ind{\V{\alpha}}{i}(t). 
\end{equation}
The coefficients $\V{\alpha}=\left\{ \ind{\V{\alpha}}{i} \right\}$ can be determined by solving a saddle point system of the form
\begin{equation}
\label{eq:saddlept1}
\begin{pmatrix}
-\M{M}\left(\V{X}\right) & \M{K}\left(\V{X}\right)  \\[2 pt]
\M{K}^*\left(\V{X}\right) & \M{0}
\end{pmatrix}
\begin{pmatrix} \V{\lambda} \\[2 pt] \V{\alpha} \end{pmatrix}=
\begin{pmatrix} \V{u}_0\left(\V{X},t\right)+\M{M}\left(\V{X}\right) \M{F}\V{X} \\[2 pt] \V{0} \end{pmatrix},  
\end{equation}
where as before $\V{X} = \left\{ \ind{\V{X}}{i}\right\}$ and $\V{\lambda} = \left\{ \ind{\V{\lambda}}{i}\right\}$. In a slight abuse of notation, we will write $\M{K}\left(\V{X}\right)$ to represent the block diagonal matrix of kinematic operators for each fiber $i$, $\M{K}\left(\V{X}\right) = \text{Diag}\left\{\M{K}\left(\ind{\V{X}}{i}\right)\right\}$, and likewise for $\M{K}^*\left(\V{X}\right)$. 

Since we expect the fiber shapes to be smooth, we use a spectral spatial discretization, described in Section \ref{sec:discspat}.  In Section\ \ref{sec:Mtot}, we break the discretized mobility matrix $\M{M}$ into three components: the local drag mobility $\ML$ given in\ \eqref{eq:MLexpl}, the finite part mobility\ \eqref{eq:MFPexpl}, and the cross-fiber mobilities\ \eqref{eq:MJFexpl}. The local drag matrix $\ML$ is the $3 \times 3$ matrix whose definition is the same as in continuum. The finite part mobility and cross-fiber mobilities require more specialized quadrature schemes since the integrals involved are near singular or singular and therefore too expensive or impossible to evaluate with direct quadrature. The basic idea of the specialized schemes is to factor out the (near) singularity, expand what remains in a monomial expansion, and compute the integrals involving monomials times the singularity analytically.  In Section\ \ref{sec:MFP}, we discuss this special quadrature scheme for the singular integrals appearing in the finite part mobility\ \eqref{eq:MFPexpl}. 

In Section\ \ref{sec:MJF}, we write a quadratic complexity discretization which uses direct quadrature to compute the SBT interaction kernels\ \eqref{eq:sbtother}. In Section\ \ref{sec:ewald}, we then discuss how to make the complexity linear over a triply periodic, sheared domain using a spectral Ewald method. Since we reformulated the inter-fiber hydrodynamics in terms of the RPY tensor (see Section \ref{sec:RPYSBT}), our Ewald splitting method is exactly the positively split Ewald method of \cite{PSRPY}, with some modifications for a non-orthogonal coordinate system \cite{SpectralSD}. In Section\ \ref{sec:nearfibs} we return to the case when the direct quadrature is insufficiently accurate and corrections are required, for which we use a recently developed monomial-expansion-based special quadrature scheme \cite{barLud} similar to that used for the finite part integral. 

Finally, in Section \ref{sec:tint}, we present a semi-implicit second-order temporal discretization that avoids nonlinear solves and requires a minimum number of evaluations of the nonlocal hydrodynamics for each timestep. For dilute suspensions, our temporal integration strategy is essentially to treat the local drag part of the mobility $\ML \M{F}\V{X}$ implicitly using an implicit trapezoidal method. We treat all of the terms involving the finite part and cross-fiber mobilities explicitly. This leaves a linear system to be solved on each fiber separately. When the suspension becomes more concentrated, this scheme breaks down as the nonlocal hydrodynamics adds stiffness to the problem. When this occurs, we treat the nonlocal and local hydrodynamics implicitly using an implicit trapezoidal method and use GMRES to solve for $\V{\alpha}$ and $\V{\lambda}$. By converting to a residual form based on the solutions of the block diagonal system for dilute suspensions, we are able to use only the minimum number of GMRES iterations necessary to achieve stability without altering accuracy.  

\subsection{Spectral spatial discretization \label{sec:discspat}}
Because the fibers are semi-flexible, their shapes are relatively smooth and can be well represented by a finite number of basis functions. This makes a spectral spatial discretization the logical choice.\rev{\footnote{\rev{Implicit in our choice of spectral discretization is the assumption that the fiber constraint forces $\V{\lambda}(s)$ and tangent vector rotation rates $\V{\Omega}(s)$ are also smooth, which is the case for sufficiently large $\delta$ in the local drag regularization\ \eqref{eq:creg}. If $\delta \ll 1$, $\V{\lambda}$ and $\V{\Omega}$ can become nonsmooth and oscillatory at the fiber endpoints, and a small number of Chebyshev modes can no longer resolve them.}}} We therefore use a first-kind Chebyshev grid for the collocation points on each fiber and Chebyshev polynomials for the basis functions $T_k(s)$, as described in Section \ref{sec:colloc}. For indefinite integration, we use the pseudo-inverse of the Chebyshev differentation matrix, and for definite integration we use Clenshaw-Curtis quadrature. Once these choices are made, the discretization of the kinematic operators $\Lop{K}$ and $\Lop{K^*}$ follows naturally in Section \ref{sec:KKstar}. The discretization of the elastic force operator $\Lop{F}$ is more subtle as the boundary conditions must be treated correctly; for this we use the rectangular spectral collocation approach of \cite{dhale15, tref17} that is described in Section \ref{sec:rsc}. Throughout this section, we consider the discretization on a single fiber.

\subsubsection{Collocation discretization \label{sec:colloc}}
Because we use a collocation discretization, each fiber is discretized as a collection of nodes $s_p$, $p=1, \dots N$, where $s_p$ is a node on a type 1 Chebyshev grid (i.e. a grid that does not include the endpoints). Our notation will shift slightly here to reflect the change from continuous to discrete. We use $\V{X}$ to refer to the $N \times 3$ matrix of fiber positions at the collocation points. The $p$th row of this matrix will be denoted by $\V{X}_p=\V{X}(s_p)$. Likewise, $\Xs$ refers to the $N \times 3$ matrix of tangent vectors at the collocation points with $\Xs_p=\Xs(s_p)$, and $\V{f}$ refers to a matrix of force densities evaluated at the nodes with rows $\V{f}_p=\V{f}(s_p)$. Meanwhile, $\V{X}(s)$ refers to the \textit{Chebyshev interpolant} for $\V{X}$ (this is actually three interpolants, one for each direction), and likewise for $\Xs(s)$. We will not try to distinguish between the unknown ``true'' fiber shape (which could have more than $N$ Chebyshev modes) and its Chebyshev approximation $\V{X}(s)$. 

The tools we use for differentiation and integration are standard \cite{trefethen2000spectral}. For differentiation, we use the Chebyshev differentiation matrix $\M{D}_N$. By $\M{D}_N\V{X}$, we mean the linear operation that takes $\V{X}$, computes the $N-1$ degree Chebyshev polynomial representation $\V{X}(s)$, differentiates it, and returns $\Xs(s)$. We also define $\M{D}_N^\dagger$, the pseudo-inverse of the Chebyshev differentiation matrix, which gives the values of the indefinite integral of a function $f(s)$ modulo an unknown constant, 
\begin{equation}
\left(\M{D}_N^\dagger \V{f}\right)_p \approx \int_0^{s_p} f(s') \, ds' +C. 
\end{equation}
For definite integration, we use Clenshaw-Curtis quadrature with weight $w_p$ associated with each collocation point, 
\begin{equation}
\int_0^L f(s')  \, ds' \approx \V{w}^T \V{f} \eqd \sum_{p=1}^N f_p w_p.
\end{equation}

\subsubsection{Discretization of $\Lop{K}$ and $\Lop{K}^*$ \label{sec:KKstar}}
 To construct a discretization of the kinematic operator $\Lop{K}$ defined in\ \eqref{eq:du}, we first need to choose the basis functions $T_k$ in the representation formula\ \eqref{eq:basis}. We choose $T_k(s)$ to be the Chebyshev polynomial of the first kind of degree $k$ on $[0,L]$. We truncate the sum at $N-2$ basis functions, 
\begin{equation}
\label{eq:basisD}
g_j(s) = \sum_{k=0}^{N-2} \alpha_{jk} T_k(s).
\end{equation}
The choice of $N-2$ for the maximum summation index is a necessary condition to make the representation $\V{U}= \Lop{K}\V{\alpha}$ unique on an $N$ point Chebyshev grid. Increasing the number of basis functions introduces degeneracy without improving the fiber representation. In particular, if the maximum index in the sum were $N-1$, integration of $g_j(s)$ in the inextensible velocity\ \eqref{eq:du} could cause $\V{U}$ to be zero at all $N$ Chebyshev nodes without $\V{\alpha}$ being zero.


\rev{Since the kinematic operators $\Lop{K}$ and $\Lop{K}^*$ act linearly on $\V{\alpha}$ and $\V{\lambda}$, respectively, they can each be discretized as matrices. Because we use a spectral discretization, however, care must be taken to avoid aliasing errors. The key step in doing this is to compute, for each $k=0, \dots, N-2$ and $j=1,2$, the integrals that describe the inextensible motions of the fiber
\begin{equation}
\label{eq:Jints}
\V{J}_q^{(k,j)} \approx \int_0^{s_q} T_k(s') \V{n}_j\left(\Xs\left(s'\right)\right) \, ds'
\end{equation}
on a grid of size $2N$ ($q=1, \dots 2N$).} To do so, we determine the Chebyshev polynomial representation of $\Xs(s)$ on the $N$ point grid and upsample it to a type 1 Chebyshev grid of $2N$ points. We then compute the normal vectors on the $2N$ grid using the polar angle representations\ \eqref{eq:Xsangle} and \eqref{eq:nangles}, and multiply the normal vectors pointwise by $T_k$ evaluated on the $2N$ grid. This gives the integrand in\ \eqref{eq:Jints} on the $2N$ grid. To integrate, we apply the matrix $\M{D}_{2N}^\dagger$ to approximate the integrals\ \eqref{eq:Jints} (modulo a constant) on the $2N$ grid. \rev{Since both $T_k$ and $\V{n}_j$ are Chebyshev polynomials of degree at most $N-1$, these integrals are exact on the grid of size $2N$, modulo a constant.}

Once $\V{J}_p^{(k,j)}$ has been computed, the discretization of $\Lop{K}^*$ is straightforward. The discrete form of the principle of virtual work\ \eqref{eq:noworkcontL2} is \rev{given by inner products of $\V{\lambda}$ with $\V{J}^{(k,j)}$, where $\V{\lambda}$ has been upsampled to a grid of size $2N$ to avoid aliasing errors, }
\begin{equation}
\label{eq:noworkd}
\M{K}^*\left(\V{X}\right)\V{\lambda}:=
\left(\begin{array}{c} 
\sum_{p=1}^{\rev{2}N} \left(\V{J}_p^{(0,1)} \cdot \left(\V{U}\V{\lambda}\right)_p\right) w_p \\[2 pt] 
\vdots \\[2 pt] 
\sum_{p=1}^{\rev{2}N} \left(\V{J}_p^{(N-2,1)} \cdot \left(\V{U}\V{\lambda}\right)_p\right) w_p\\[2 pt] \hdashline
 \sum_{p=1}^{\rev{2}N} \left(\V{J}_p^{(0,2)} \cdot \left(\V{U}\V{\lambda}\right)_p\right) w_p \\[2 pt] 
\vdots \\[2 pt] 
\sum_{p=1}^{\rev{2}N} \left(\V{J}_p^{(N-2,2)} \cdot \left(\V{U}\V{\lambda}\right)_p\right) w_p\\[2 pt] \hdashline 
\M{w}^T\V{\lambda} \end{array}\right) = 
\left(\begin{array}{c} 0 \\[2 pt] \vdots \\[2 pt] 0\\ \hdashline 0 \\[2 pt] \vdots \\[2 pt] 0\\[2 pt] \hdashline \V{0}\end{array}\right).
\end{equation}
This defines $\M{K}^*\left(\V{X}\right)$ as a $(2N+1) \times 3N$ matrix acting on a $3N$ vector $\V{\lambda}$, \rev{which has been upsampled to (i.e., evaluated on) the grid of size $2N$ by applying the upsampling matrix $\M{U}$}. Because $\V{\lambda}$ discretely integrates to zero ($\V{w}^T\V{\lambda}=\V{0}$), adding a constant to $\V{J}_p^{(k,j)}$ does not change the first two rows of $\M{K}^*\left(\V{X}\right)\V{\lambda}$. Thus the fact that $\V{J}_p^{(k,j)}$ gives the integrals modulo a constant is not relevant in the formation of $\M{K}^*$. 

\rev{The discretization of $\M{K}$ is less straightforward. The main issue is that the integrals $\V{J}_p^{(k,j)}$ are exact on the grid of size $2N$, but the fiber velocity is defined on a grid of size $N$. One way around this is to simply double the grid size, i.e., define $\M{K}\V{\alpha}$ on a grid of size $2N$. This would, however, necessitate computing $\M{M}\V{f}$ on a grid of size $2N$ as well, which is unnecessarily expensive. Instead, we incur some aliasing error and downsample the velocity $\V{J}_p^{(k,j)}$ to an $N$ point grid. Specifically, we discretize the matrix $\M{K}$ as
\begin{equation}
\label{eq:Kalph}
\left(\M{K}\left(\V{X}\right)\V{\alpha}\right)_p = \Ubar + \sum_{j=1}^2 \sum_{k=0}^{N-2} \alpha_{jk} \rev{\left(\M{R}\V{J}^{(k,j)}\right)_p}
\end{equation}
so that  $\M{K}=\M{R}\M{J}$ is a $3N \times (2N+1)$ matrix which acts on the $2N+1$ vector $\V{\alpha} = \left( \alpha_{jk}, \Ubar \right)$ to give the three components of the velocity at $N$ points on the Chebyshev grid.   The matrix $\M{R}$ is a downsampling matrix which gives a reduced-order representation of $\V{J}^{(k,j)}$ on the $N$ point grid (the unknown constant in $\V{J}^{(k,j)}$ can be folded into the constant velocity $\Ubar$). If we concisely write $\M{K}^* \V{\lambda}$ as defined in\ \eqref{eq:noworkd} as $\M{K}^* \V{\lambda} = \M{J}^T \M{W}_{2N} \M{U}\V{\lambda}$, where $\M{W}_{2N}$ is the diagonal matrix of Clenshaw-Curtis quadrature weights on the $2N$ grid, and discretize the inner products in\ \eqref{eq:Kstardef} as
\begin{equation}
\label{eq:ip1}
\langle \V{\alpha}, \M{K}^* \V{\lambda}\rangle = \V{\alpha}^T \M{J}^T \M{W}_{2N} \M{U}\V{\lambda} = \langle  \M{K} \V{\alpha},\V{\lambda} \rangle =\left(\M{U}\M{K}\V{\alpha}\right)^T  \M{W}_{2N} \M{U}\V{\lambda}, 
\end{equation}
we obtain the weighted least squares downsampling matrix
\begin{equation}
\M{R} = \left(\M{U}^T \M{W}_{2N} \M{U}\right)^{-1}  \M{U}^T \M{W}_{2N}.
\end{equation}
}

\subsubsection{Discretization of $\Lop{F}$ \label{sec:rsc}}
We use rectangular spectral collocation \cite{tref17, dhale15} to discretize the bending force operator $\Lop{F}$ with the boundary conditions\ \eqref{eq:frBCs}. We recall that the matrix $\V{X}$ gives the positions of the fiber on an $N$ point type 1 Chebyshev grid that does \textit{not} include the boundaries. 

In rectangular spectral collocation, we compute an upsampled representation $\widetilde{\V{X}}$ of $\V{X}$. Since there are four boundary conditions, the upsampled representation is on a $\widetilde{N}=N+4$ point type 2 Chebyshev grid that includes the endpoints. The unique configuration $\widetilde{\V{X}}$ can be obtained by solving
\begin{equation}
\label{eq:deftilde}
\begin{pmatrix} \M{A} \\[2 pt] \M{B} \end{pmatrix} \widetilde{\V{X}} = \begin{pmatrix} \V{X}\\[2 pt] \V{0} \end{pmatrix},
\end{equation}
where $\M{A}$ is an $N \times \widetilde{N}$ resampling matrix and $\M{B}$ is a $4 \times \widetilde{N}$ matrix that encodes the boundary conditions. The linear operation $\M{A}\widetilde{\V{X}}$ has the effect of computing the Chebyshev interpolant of $\widetilde{\V{X}}$ on the $\widetilde{N}=N+4$ point grid and evaluating it at the $N$ original gridpoints. The first block equation simply states that the downsampled $\widetilde{\V{X}}$ has to be the original $\V{X}$. In the next block, the product $\M{B}\widetilde{\V{X}}$ is a vector with 4 entries. The first two entries are the Chebyshev interpolant approximation to $\partial_s^2\widetilde{\V{X}}\left(s=0,L\right)$, respectively, and the second two entries are likewise an approximation to $\partial_s^3\widetilde{\V{X}}\left(s=0,L\right)$. Thus the second block equation simply states that the boundary conditions are satisfied on the type 2 grid, and any modifications to the BCs would modify $\M{B}$ in this formulation. 

The use of a type 1 grid for $\V{X}$ and a type 2 grid for $\widetilde{\V{X}}$ is a sufficient condition for the left-hand side of system\ \eqref{eq:deftilde} to be invertible (see \cite{dhale15} for details). We can therefore write 
\begin{equation}
\label{eq:getX}
\widetilde{\V{X}} = \begin{pmatrix} \M{A} \\[2 pt] \M{B} \end{pmatrix}^{-1} \begin{pmatrix} \V{X} \\[2 pt] \M{0} \end{pmatrix} \eqd \M{E}\V{X}.
\end{equation}
In summary, $\widetilde{\V{X}}$ is the unique upsampled configuration that satisfies the problem boundary conditions and gives $\V{X}$ when downsampled. This is similar to ``ghost cells'' in finite difference schemes which take on unique values so that the boundary stencils satisfy the BCs to some order. The rectangular spectral collocation method can therefore be thought of as a generalization of ghost cell techniques for finite difference methods to collocation-based spectral methods.  

Once the configuration $\widetilde{\V{X}}$ is known, the elastic force density can be computed on it as $\widetilde{\V{f}}^\kappa=-\kappa \M{D}_{\widetilde{N}}^4 \widetilde{\V{X}}$. The elastic force density $\widetilde{\V{f}}^\kappa$ is then downsampled to the original $N$ point type 1 grid to give the final result, 
\begin{equation}
\label{eq:fE}
\V{f}^\kappa=\M{A}\widetilde{\V{f}}^\kappa=-\kappa \left(\M{A}\M{D}_{\widetilde{N}}^4 \M{E}\right)\V{X}:=\M{F}\V{X},
\end{equation}
which defines $\M{F}$, the discrete analogue of $\Lop{F}$. In a slight abuse of notation, we will use the notation $\M{F}\V{X}$ to refer to the bending force calculation on either a single fiber (where $\M{F}$ is as defined in\ \eqref{eq:fE}), or a collection of fibers (where $\M{F}$ is a diagonal block matrix composed of smaller matrices that are defined in\ \eqref{eq:fE}); the meaning should be clear from the context.

\subsection{Discretization of $\Lop{M}$ \label{sec:Mtot}}
We discretize the mobility operator defined in\ \eqref{eq:fibevcont} by computing the relative velocity of point $p$ on fiber $i$ as 
\begin{equation}
\label{eq:discmob}
\ind{\left(\M{M}\V{f}\right)}{i}_p = \ML \left(\piXs{i}{p} ; c_p\right) \pif{i}{p}+ \left(\MFPd\left(\ind{\V{X}}{i}\right)\ind{\V{f}}{i}\right)_p + \sum_{j \neq i} \MJFd\left(\piX{i}{p}, \ind{\V{X}}{j}\right)\ind{\V{f}}{j}. 
\end{equation}
The $3 \times 3$ matrix $\ML$ given in\ \eqref{eq:ML} is unchanged from the continuum, and $c_p=\rev{\bar{c}(s_p;\delta)}$. The matrix $\MFPd$ computes an approximation to the finite part integral\ \eqref{eq:Mfp} (see Section\ \ref{sec:MFP}), and $\MJFd\left(\piX{i}{p},\ind{\V{X}}{j}\right)\ind{\V{f}}{j}$ is the velocity at point $\piX{i}{p}$ induced by fiber $j$ (see Section\ \ref{sec:MJF}). 


We will need notation to separate the local part of the discrete mobility\ \eqref{eq:discmob}, which is easy to invert, from the nonlocal part, which is not. For this we write
\begin{equation}
\label{eq:fibiLNL}
\ind{\left(\M{M}\V{f}\right)}{i}_p = \ML \left(\piXs{i}{p} ; c_p\right) \pif{i}{p} + \ind{\left(\MNL\left(\V{X}\right) \V{f}\right)}{i}_p, 
\end{equation}
where the nonlocal part of the mobility matrix $\MNL\left(\V{X}\right)$ is a function of the collection of fibers $\V{X}$ and acts on the collection of force densities $\V{f}$. For the collection of fibers, we will simply write the splitting\ \eqref{eq:fibiLNL} as 
\begin{equation}
\label{eq:allfibNL}
\M{M}\left(\V{X}\right)=  \MLD\left(\V{X}\right) + \MNL\left(\V{X}\right) ,
\end{equation}
where the block diagonal matrix $\M{M}^{\text{LD}}$ is composed of a collection of $3 \times 3$ local drag matrices on the diagonal. 

\subsubsection{Discretization of $\MFP$ \label{sec:MFP}  \label{sec:MFP}}
Here we discretize the finite part integral\ \eqref{eq:Mfp}. Since the finite part integral involves only a single fiber, we use $\V{X}$ here to denote a matrix of $N \times 3$ positions for a single fiber. Substituting the definition of the Stokeslet\ \eqref{eq:Slet}, we have
\begin{equation} 
\label{eq:MFPfull}
\left(\MFP\left[\V{X}\right]\V{f}\right)(s) =  \EPMI \int_{0}^L \left( \frac{\M{I}+\left(\hat{\V{R}}\hat{\V{R}}\right)\left(\V{X}(s),\V{X}(s')\right)}{\norm{\V{R}\left(\V{X}(s),\V{X}(s')\right)}}\V{f}(s') -\left(\frac{\V{I}+\Xs(s)\Xs(s)}{|s-s'|}\right) \V{f}(s)\right) \, ds'
\end{equation}
for any coordinate on the fiber $s$. We seek to evaluate the finite part integral\ \eqref{eq:MFPfull} at $s=s_p$ on a given fiber. Because each term in the integrand is singular, the integral\ \eqref{eq:MFPfull} cannot directly be evaluated with Clenshaw-Curtis quadrature. One way around this difficulty is to simply skip the singular point in the quadrature, which results in a second-order accurate scheme. Because this destroys the spectral accuracy of our formulation, we seek an improved quadrature that handles the singularity analytically. 

In \cite{ts04}, the integrand was regularized to make it non-singular, and a product integration scheme was used to compute the resulting regularized integral \cite[Section~3.1]{ts04}. The justification for the regularization is that the self mobility operator $\Lop{M}$ is actually not invertible, since its null space contains force densities $\V{f}$ with frequencies higher than $1/\epsilon$. G\"otz \cite{gotz2001interactions} and Tornberg and Shelley \cite[Appendix~B]{ts04} show this analytically by considering a straight fiber and expanding $\V{f}$ as a sum of Legendre polynomials (which diagonalize $\MFP$). They show that the centerline velocity $\V{U}$ for a single fiber\ \eqref{eq:onefib} uniquely gives the Stokeslet strength $\V{f}(s)$ via $\V{U}=\Lop{M}\V{f}$ if the maximum number of polynomials that contribute to $\V{f}$ is less than $\mathcal{O}(1/\epsilon)$. Intuitively, adding polynomials of degree larger than $1/\epsilon$ introduces length scales into $\V{f}$ which are less than $\epsilon$ and cannot be accounted for by SBT. In the discretization of \cite{ts04}, the number of points can exceed $1/\epsilon$ for $\epsilon=10^{-2}$, and so regularization of the integrand is required. 

Because we use a spectral basis with smooth fiber shapes, we never exceed the $\mathcal{O}(1/\epsilon)$ threshold for the number of Chebyshev polynomials. Indeed, having less than $\mathcal{O}(1/\epsilon)$ Chebyshev points (polynomials) is a sensible restriction on the numerics. After all, using such a large number of points (spectral modes) is antithetical to the philosophy of SBT, which, unlike IB methods, eliminates the need to resolve the length scale $\epsilon$. \rev{Mori and Ohm analyze this ``spectral truncation'' approach for infinite (periodic) fibers and find that it yields errors of order at worst $\epsilon$ with respect to the true solution for a ``slender body'' PDE on the fiber surface \cite[Eq.~(20)]{mori2020accuracy}. This accuracy is asymptotically equivalent to that obtained using the regularization of Tornberg and Shelley \cite[Eq.~(22)]{mori2020accuracy}.}

With this in mind, we do not modify the integrand of the finite part integral\ \eqref{eq:MFPfull}. Rather, we use a spectrally accurate method developed in \cite{tornquad} to compute the action of the finite part integral on the fiber force density. The key idea is to isolate the singularity by writing the integrand in\ \eqref{eq:MFPfull} as $\V{g}(s',s) (s'-s)/|s'-s|$ for some function $\V{g}(s',s)$. In particular, we observe that the finite part integral\ \eqref{eq:MFPfull} can be written as 
\begin{equation}
\label{eq:Jrewrite}
\left(\MFP\left[\V{X}\right]\V{f}\right)(s) = \int_0^L \V{g}(s,s') \frac{s'-s}{|s'-s|} \, ds' = \frac{L}{2}\int_{-1}^1 \V{g}(\eta, \eta') \frac{\eta'-\eta}{|\eta'-\eta|} \, d\eta', 
\end{equation}
where $\eta=-1+2s/L$ is a rescaled arclength coordinate on $[-1,1]$ and 
\begin{align}
\label{eq:defg}
\V{g}(s,s') = \EPMI \Bigg{[}\left(\V{I}+\left(\hat{\V{R}}\hat{\V{R}}\right) \left(\V{X}(s),\V{X}(s')\right)\right) &\frac{|s'-s|}{\norm{\V{R}\left(\V{X}(s),\V{X}(s')\right)}} \V{f}(s') \nonumber \\[2 pt] & - \left(\V{I}+\Xs(s) \Xs(s)\right) \V{f}(s)\Bigg{]} \frac{1}{s'-s}. 
\end{align}
The computation is now tractable since $\V{g}$ has a limit as $s' \rightarrow s$. The limit is easily computed by adding and subtracting $\left(\V{I}+\Xs(s)\Xs(s)\right)\V{f}(s')$ inside the square bracket, and Taylor expanding around $s'=s$, to obtain
\begin{equation}
\lim_{s' \to s} \V{g}(s,s') = \EPMI\left(\frac{1}{2}\left(\Xs(s)\V{X}_{ss}(s)+\V{X}_{ss}(s)\Xs(s)\right)\V{f}(s) + \left(\V{I}+\Xs(s)\Xs(s)\right)\V{f}_s(s)\right). 
\end{equation}
Since $\V{g}$ is smooth, we can approximate it by a polynomial expressed in a monomial basis on $[-1,1]$, 
\begin{equation}
\label{eq:gmono}
\frac{L}{2}\V{g}(\eta,\eta') \approx \sum_{k=0}^{N-1} \V{c}_k(\eta) (\eta')^k, 
\end{equation}
where $\V{c}_k$ is a vector of 3 coefficients for each $\eta$. 

We are now ready to discretize $\MFP\left[\V{X}\right]$ with the matrix representation $\MFPd\left(\V{X}\right)$. Substituting the monomial expansion\ \eqref{eq:gmono} into the finite part integrand\ \eqref{eq:Jrewrite} and computing the integrals involving monomials and the singularity analytically, we get
\begin{gather}
\label{eq:expandmono}
\left(\MFPd\left(\V{X}\right)\V{f}\right)_p= \sum_{k=0}^{N-1} \V{c}_k(\eta_p) \int_{-1}^1 (\eta')^k \frac{\eta'-\eta_p}{|\eta'-\eta_p|} \, d\eta' = \sum_{k=0}^{N-1} \V{c}_k(\eta_p) q_k(\eta_p),\\[6 pt]
\label{eq:qints}
\text{where} \qquad 
q_k(\eta_p) = \int_{-1}^1 (\eta')^k \frac{\eta'-\eta_p}{|\eta'-\eta_p|} \, d\eta'= \frac{1+(-1)^{k+1}-2\eta_p^{k+1}}{k+1}.
\end{gather}

An adjoint method can be used to accelerate our computation of the product\ \eqref{eq:expandmono}. Let us introduce the Vandermonde matrix $\M{V}$ with entries $V_{pq}=\eta_p^q$. Let $\V{g}^p\left(\V{X}\right)$ be the $N \times 3$ matrix with rows $\V{g}^p_q = \V{g}(\eta_p,\eta_q)$, with $\V{g}$ defined in\ \eqref{eq:defg}. The $N \times 3$ matrix $\V{c}$ of coefficients of the three polynomial interpolants of the columns of $\V{g}^p$ is $\V{c} = \V{V}^{-1}\V{g}^p$. If $\V{q}$ is an $N$ vector with elements $q_k = q_k(\eta_p)$ as given in\ \eqref{eq:qints}, then the product\ \eqref{eq:expandmono} can be computed efficiently as  
\begin{equation}
\label{eq:specscheme}
\left(\MFPd\left(\V{X}\right)\V{f}\right)_p = \V{c}^T \V{q} = \left(\V{V}^{-1}\V{g}^p\left(\V{X}\right)\right)^T \V{q} = \left(\V{g}^p\left(\V{X}\right)\right)^T \left(\V{V}^{-T}\V{q}\right) \eqd \left(\V{g}^p\left(\V{X}\right)\right)^T \V{b}. 
\end{equation}
Since $\V{b}= \V{V}^{-T}\V{q}$ does not depend on the fiber configuration, it can be precomputed using pivoted $LU$ factorization for each $p=1,2,\dots N$ at the beginning of the simulation. The Vandermonde matrix must be sufficiently well-conditioned to do this calculation accurately; specifically, the fiber discretization can have at most $\sim$40 points in double precision. If higher accuracy is needed, then the fiber must be split into multiple panels or higher precision arithmetic must be used to compute $\V{b}$ in\ \eqref{eq:specscheme}.

\subsubsection{Discretization of $\MJF$ \label{sec:MJF}  \label{sec:MJF}}
In this section, we describe the simplest discretization of inter-fiber hydrodynamic interactions. We recall the definition of the velocity induced by fiber $j$ at point $s$ on fiber $i$ from\ \eqref{eq:sbtother}, 
\begin{align}
\label{eq:dintvel}
\ind{\V{v}}{j}\left(\ind{\V{X}}{i}(s)\right)& = \McNot \\[2 pt] &= \nonumber \EPMI \int_0^{L} \Knel{\ind{\V{X}}{i}(s),\ind{\V{X}}{j}(s'),\rev{\dco}\left(\epsilon L \right)^2}\ind{\V{f}}{j}(s') \, ds'.  
\end{align}
Given a discrete Chebyshev node on fiber $i$, $\piX{i}{p}$, the total disturbance velocity is a sum of the flows generated by all other fibers $j \neq i$. We can therefore restrict our attention to the calculation of the velocity induced by a single ``source'' fiber at a single ``target'' point on another fiber. 

The simplest approach is to discretize the interaction velocity\ \eqref{eq:dintvel} by Clenshaw-Curtis quadrature, 
\begin{align}
\label{eq:dirquad}
\ind{\V{v}}{j}\left(\piX{i}{p}\right) & \approx \EPMI \sum_{q=1}^N w_q \Knels{\piX{i}{p}}{\piX{j}{q}}{\rev{\dco}(\epsilon L)^2}\ind{\V{f}}{j}_q \eqd \MJFd\left(\ind{\V{X}}{i}_p, \ind{\V{X}}{j}\right)\ind{\V{f}}{j}.
\end{align}
The key challenge in evaluating\ \eqref{eq:dirquad} is the quadratic complexity; for each Chebyshev point $\piX{i}{p}$ we must sum over all others. We address this in Section\ \ref{sec:ewald} using the positively split Ewald method \cite{PSRPY}.

While the direct quadrature\ \eqref{eq:dirquad} represents the simplest way to discretize $\MJF$, it becomes inadequate when fibers $i$ and $j$ approach each other. In fact, the Stokeslet-doublet combination kernel $\M{S}_D$ defined in\ \eqref{eq:Kdef} becomes singular if the Chebyshev interpolant $\ind{\V{X}}{j}(s)$ of fiber $j$ approaches the target $\piX{i}{p}$ for some value of $s$. This singularity occurs because the velocity of one fiber due to another\ \eqref{eq:sbtother} only makes sense physically if the two fiber cross sections are not overlapping \rev{(for overlapping fibers the RPY kernel is different, see the second line of\ \eqref{eq:rpyknel})}.  

While it makes little physical sense for fiber cross sections to overlap, it is numerically possible. In this case, we set the velocity at the target point on fiber $i$ to be equal to the centerline velocity at the closest point on fiber $j$. Let us denote the minimum distance from $\piX{i}{p}$ to\footnote{Recall that we cannot and do not distinguish between the unknown ``true'' fiber shape and its Chebyshev interpolant and denote both with $\ind{\V{X}}{j}(s)$.} $\ind{\V{X}}{j}(s)$ by $d$ and denote the closest point on fiber $j$ to $\ind{\V{X}}{i}$ as $\ind{\V{X}}{j}(s^*)$. Then if the two cross sections are (almost) overlapping, we set 
\begin{equation}
\label{eq:CLvel1}
\ind{\V{v}}{j}\left(\piX{i}{p}\right)  = \ML\left(\ind{\Xs}{j}\left(s^*\right); \rev{\bar{c}\left(s^*;\delta \right)}\right)\ind{\V{f}}{j}\left(s^*\right) + \left(\MFPd\left(\ind{\V{X}}{j}\right)\ind{\V{f}}{j}\right)(s^*) \quad \text{for $d \leq  \rev{2\aRPY}$}, 
\end{equation}
so that the influence of fiber $j$ on the target (which is inside the cross section of fiber $j$) is the same as if the target were actually on the centerline of fiber $j$, \rev{with the radius $\aRPY \approx 1.12\epsilon L$ given in\ \eqref{eq:rpyradius}.} By $\ind{\Xs}{j}\left(s^*\right)$, we mean the Chebyshev polynomial $\ind{\Xs}{j}(s)$ evaluated at $s^*$, and likewise for the remaining terms.

For non-overlapping cross sections \rev{($d > 2 \aRPY$)}, the expression for the interaction velocity\ \eqref{eq:sbtother} has to be changed \rev{because evaluating\ \eqref{eq:sbtother} at $d=2\aRPY$ may not be exactly equal to the centerline velocity on fiber $j$\ \eqref{eq:CLvel1}  because of the slenderness approximation used in SBT}. We therefore set the velocity to be the interaction velocity\ \eqref{eq:sbtother} only if $d \geq \rev{4b}$. Between \rev{$2b \leq d \leq 4b$}, we linearly interpolate between the interaction velocity\ \eqref{eq:sbtother} and centerline velocity\ \eqref{eq:CLvel1}. This interpolation procedure is almost identical to that of \cite{ts04}, except we use a different integral kernel $\M{S}_D$ for the interaction velocity\ \eqref{eq:sbtother}, and we estimate $s^*$ using a more robust procedure described in Section\ \ref{sec:nearfibs}. 

Despite our modifications to the interaction velocity for contacting cross sections, the kernel in\ \eqref{eq:sbtother} is still nearly singular, and some quadrature scheme other than direct quadrature\ \eqref{eq:dirquad} is required to accurately determine the interaction velocity for \rev{$d > 2b$}, as we discuss in Section\ \ref{sec:nearfibs}. 

\subsection{Fast summation \label{sec:ewald}}
Putting off for the moment the possible near-singular nature of the interaction integrals\ \eqref{eq:sbtother}, suppose that we discretize every integral using direct quadrature\ \eqref{eq:dirquad}. If $F$ is the number of fibers, each of which is discretized using $N$ Chebyshev points, then the direct evaluation of nonlocal hydrodynamics using\ \eqref{eq:dirquad} requires $\mathcal{O}\left((NF)^2\right)$ operations. In this section, we discuss our choice of fast algorithm to accelerate the evaluation of these many-body sums under triply periodic boundary conditions. We recall that the kernel $\M{S}_D$ in\ \eqref{eq:dirquad} is the RPY kernel for non-overlapping spheres of radius \rev{$\aRPY = \epsilon L e^{3/2}/4$}. Thus\ \eqref{eq:dirquad} reduces to summing the RPY kernel over all pairs of points $\left(\piX{i}{p},\piX{j}{q}\right)$, which is a well-studied problem that can be treated with a number of fast algorithms \cite{PSRPY, liang2013fast, guan2018rpyfmm, yan2018universal}. Alternatively, the kernel $\M{S}_D$ can be viewed as a linear combination of Stokes singularities, and fast algorithms for the individual singularities can be applied \cite{Lindbo:2010ha, af14fast, Klinteberg:2016bj}. In Section\ \ref{sec:ewblob} we describe the Positively Split Ewald (PSE) approach of \cite{PSRPY, SpectralSD}, which assumes a constant value of $b$ (and therefore $\epsilon L$), across all of the fibers; this assumption can be relaxed \cite{SE_Multiblob_SD}. Because we are interested in rheological applications, we use the method of \cite{SpectralSD} to extend the fast Ewald summation technique of \cite{PSRPY} to a parallelepiped sheared unit cell (see Section \ref{sec:shearcoords}). 

For overlapping spheres the RPY kernel $\M{S}_\text{RPY}$ differs from the SBT kernel $\M{S}_D$ in\ \eqref{eq:dirquad}, and is given by \cite{wajnryb2013generalization}
\begin{equation}
\label{eq:rpyclose}
\KRPY{\V{x}}{\V{y}}{b}= \EPMI \left(\left(\frac{4}{3b}-\frac{3\norm{\V{R}}}{8b^2}\right)\M{I}+\frac{\norm{\V{R}}}{8b^2} \hat{\V{R}}\hat{\V{R}}\right) \quad \text{ if } \norm{\V{R}} < 2b, 
\end{equation}
where $\V{R}=\V{x}-\V{y}$. This means that the PSE method mistakenly computes the RPY kernel\ \eqref{eq:rpyclose} between a pair of points separated by a distance less than $2b$ instead of the desired SBT kernel. We need not worry about this, however, since points that are a distance less than $2b$ yield a target-fiber pair for which \rev{we set the velocity at the target to the fiber centerline velocity\ \eqref{eq:CLvel1}.} More generally, even for points farther apart than $2b$ there will be some number of target-fiber pairs for which the direct quadrature\ \eqref{eq:dirquad} fails. Our approach to this is to rely on Ewald splitting to periodize and accelerate the many-body summation, then subtract the free space RPY kernel between the problematic fiber and target from the result. Subtracting the free space kernel leaves the periodic images of the sum, which have been correctly accounted for by Ewald splitting (since they are distant). We then handle the free space kernel $\M{S}_D$ for problematic pairs of fibers and targets using the special quadrature algorithm described in Section\ \ref{sec:nearfibs}. 

\subsubsection{Sheared coordinate system \label{sec:shearcoords}}
\begin{figure}
\centering
\includegraphics[width=0.3\textwidth]{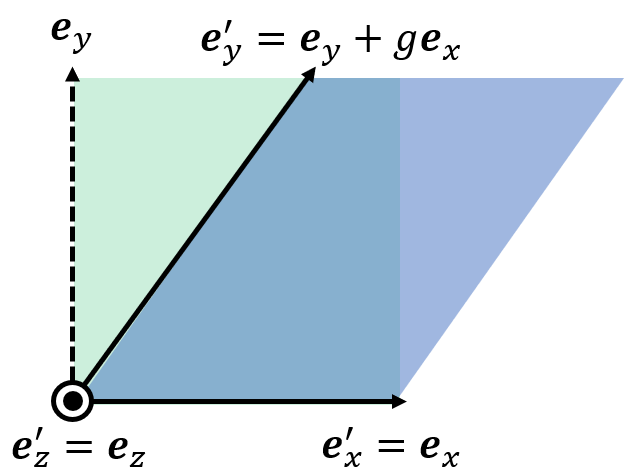}
\caption{The sheared/parallelepiped periodic simulation cell. We denote the dimensionless shear by $g$. The green area shows the periodic cell (2D projection) when $g=0$, and the blue area is the periodic cell for nonzero $g$.}
\label{fig:shearcell}
\end{figure}

In order to implement a shear flow in periodic boundary conditions, a strained coordinate system is necessary. The PSE method was extended to sheared cells in \cite{SpectralSD}, but here we give a more detailed description for completeness. We assume (without loss of generality) that $x$ is the flow direction, $y$ is the gradient direction, and $z$ is the vorticity direction. 
Let the total nondimensional strain be $g(t)$. In Fig.\ \ref{fig:shearcell} we define a strained coordinate system with axes
\begin{equation}
\label{eq:axes}
\V{e}_{x'} = \V{e}_x, \qquad \V{e}_{y'} = \V{e}_y+g(t)\V{e}_x, \qquad \V{e}_{z'} = \V{e}_z,
\end{equation}
and strained wave numbers
\begin{equation}
\label{eq:swnums}
k'_x=k_x, \qquad k'_y = k_y+g(t)k_x, \qquad k'_z = k_z.
\end{equation}
Here $k_x$, $k_y$, and $k_z$ are the wave numbers when the periodicity is over the $x$, $y$, and $z$ directions, while $k_x'$, $k_y'$, and $k_z'$ are the wave numbers when the periodicity is over the $x'$, $y'$, and $z'$ directions. 

The transformation between the two coordinate systems is given by
\begin{gather}
\label{eq:primes}
\V{x}' \eqd \begin{pmatrix} x' \\ y' \\ z' \end{pmatrix} = \begin{pmatrix} 1 & -g(t) & 0 \\ 0 & 1 & 0 \\ 0 & 0 & 1 \end{pmatrix}\begin{pmatrix} x \\ y \\ z \end{pmatrix} \eqd \M{L}\V{x}. 
\end{gather}
In the unsheared to sheared transformation\ \eqref{eq:primes}, the sheared coordinates $x'$, $y'$, and $z'$ are all periodic on $[0,L]^3$ (the blue simulation cell in Fig.\ \ref{fig:shearcell}). 

Now we use the transformation\ \eqref{eq:primes} to transform the derivative operators in the unsheared coordinate system to the sheared one,  
\begin{equation}
\frac{\partial}{\partial x} = \frac{\partial}{\partial x'} \qquad \frac{\partial}{\partial y} = \frac{\partial}{\partial y'}-g(t)\frac{\partial}{\partial x'} \qquad \frac{\partial}{\partial z} = \frac{\partial}{\partial z'}. 
\end{equation}

We therefore have the Laplacian in the transformed space as \cite{moto11}
\begin{equation}
\Delta = \left(\frac{\partial^2}{\partial x'^2}+\left(\frac{\partial}{\partial y'}-g(t)\frac{\partial}{\partial x'}\right)^2+\frac{\partial^2}{\partial z'^2}\right). 
\end{equation}
In Fourier space, $\widehat{\Delta} = \V{k}' \cdot \V{k}'$, where 

\begin{equation}
\label{eq:kprime}
\bm{k}'= (k_x',k_y'-g(t)k_x',k_z'), 
\end{equation}

Using the sheared to unsheared transformation in\ \eqref{eq:swnums}, it is easy to see that $k':=\norm{\bm{k}'}=\norm{(k_x,k_y,k_z)}:=k$. It follows that we can simply replace $k$ in any isotropic Fourier calculations by $k'$ to use a Fourier method in the sheared coordinate system \cite{wheeler97, moto11}. In Appendix \ref{sec:shearBC}, we verify this formulation by considering a set of particles that can be viewed periodically in two ways. This appendix verifies our correct treatment of the sheared periodic boundary conditions.

\subsubsection{Ewald splitting for direct quadratures \label{sec:ewblob}}
Let $\V{x}$ be a target point on the centerline of a fiber. To evaluate the direct quadrature \eqref{eq:dirquad} with periodic boundary conditions, we first compute 
\begin{equation}
\label{eq:velalldirect}
\V{U}^{(\text{PSE})}\left(\V{x}\right) \eqd \sum_P \sum_{i} \KRPY{\V{x}}{\V{y}_i^{(P)}}{b}\V{F}_i \eqd \sum_i \M{S}_\text{RPY}^{(P)}\left(\V{x},\V{y}_i;b\right)\V{F}_i  
\end{equation} 
where $\V{F}_i$ is a force (force density $\times$ weight) assigned to the point $\V{y}_i$. For each target point $\V{x}=\V{y}_j$ for some $j$, the sum is over all discrete fiber points $\V{y}_i$ (including $\V{y}_j$) and over triply periodic images of the points $\V{y}_i$ in the sheared coordinate system $P$. The periodized RPY kernel is denoted by $\M{S}_\text{RPY}^{(P)}$.

We use the Ewald splitting of \cite{PSRPY, SpectralSD} to accelerate the computation of the many-body sum\ \eqref{eq:velalldirect} on a periodic domain. The idea of Ewald splitting or Ewald summation is to split the RPY kernel into a smooth long-ranged part and a remaining short-ranged part. The smooth ``far field'' part has an exponential decay in Fourier space and can be done by standard Fourier methods (namely the non-uniform FFT), and the ``near field'' part decays exponentially in real space and can be truncated so that it is nonzero for $\mathcal{O}(1)$ neighboring points (sources) per target. 

Let $\V{x}'$ and $\V{y}'$ be the coordinates of the points $\V{x}$ and $\V{y}$ in the sheared domain using the coordinate transformation\ \eqref{eq:primes}. The periodic RPY tensor for a sphere with radius $b$ can then be written on the sheared domain as
\begin{equation}
\M{S}_\text{RPY}^{(P)}\left(\V{x}',\V{y}';b\right)=\frac{1}{V\mu}\sum_{\small{\bm{k}'\neq \bm{0}}} e^{i\V{k}' \cdot \left(\V{x}'-\V{y}'\right)} \frac{1}{k'^2}\left(\bm{I}-\hat{\bm{k}'}\hat{\bm{k}'}^T\right)\text{sinc}^2\left(k'b\right), 
\end{equation} 
where $V$ is the domain volume and $\V{k}' = 2\pi \V{m}/L_d$, where $\V{m}$ is a vector of three integers and $L_d$ is the periodic domain length. Using the screening function of Hasimoto \cite{Hsplit},
\begin{equation}
\label{eq:HspE}
H(k',\xi)=\left(1+\frac{k'^2}{4\xi^2}\right)e^{-k'^2/4\xi^2}, 
\end{equation}
we split the periodic kernel $\M{S}_\text{RPY}^{(P)}$ into a far field and near field component, $\M{S}_\text{RPY}^{(P)}=\M{S}_\text{RPY}^{\text{(FF)}}+\M{S}_\text{RPY}^{\text{(NF)}}$, where the far field is given in Fourier space by 
\begin{equation}
\label{eq:ufar}
\M{S}_\text{RPY}^{\text{(FF)}}\left(\V{x}',\V{y}';b\right)=\frac{1}{V\mu}\sum_{\small{\bm{k}'\neq \bm{0}}} e^{i\V{k}' \cdot \left(\V{x}'-\V{y}'\right)} \frac{1}{k'^2}\left(\bm{I}-\hat{\bm{k}'}\hat{\bm{k}'}^T\right)\text{sinc}^2\left(k'b\right) H(k',\xi). 
\end{equation}
Here $\xi$ is a splitting parameter that controls the decay of the far field kernel in Fourier space and of the near field kernel in real space, and is chosen to optimize performance. The total far field sum is obtained by summing the far field kernel\ \eqref{eq:ufar} over all points $\V{y}'$, with the $k'=0$ mode set to zero since in continuum the total force on the system is zero. We use standard NUFFT algorithms (in particular the Flatiron NUFFT library \cite{barnettES}) to compute these sums at all points $\V{x}'$ in log-linear time in the number of points (see \cite{PSRPY, SpectralSD} for more details). The Flatiron NUFFT relies on a new ``exponential of a semicircle'' kernel to do spreading and interpolation \cite{barnettES}, which is more efficient than the traditional Gaussian featured in \cite{PSRPY, SpectralSD}. See \cite[Eq.~(31)]{tornbergfft} for error estimates using the exponential of a semicircle kernel, although these have yet to be extended to sheared domains in the manner of \cite[Eq.~(55)]{SpectralSD}.

Assuming that the near field decays rapidly enough that Fourier series can be replaced by Fourier integrals, the near field mobility can be computed in real space by inverse transforming its Fourier space representation,
\begin{equation}
\label{eq:nearRPY}
\M{S}_\text{RPY}^{\text{(NF)}}\left(\V{x},\V{y}; b \right)=F(r,\xi,b)\left(\bm{I}-\hat{\bm{r}}\hat{\bm{r}}^T\right)+G(r,\xi,b)\hat{\bm{r}}\hat{\bm{r}}^T, 
\end{equation}
where $\V{r}=\left(\V{x}-\V{y}\right)^*$, $r=\norm{\V{r}}$, and the $*$ denotes the nearest periodic image in the sheared domain (blue in Fig.\ \ref{fig:shearcell}). The exact forms of $F$ and $G$ are given in \cite[Appendix~A]{PSRPY}. The total near field sum is computed at $\V{x}$ by summing the near field kernel\ \eqref{eq:nearRPY} over neighboring points whose minimum image distance from $\V{x}$ is less than a precomputed value $r^*$. We choose $r^*$ so that the velocities are computed to a relative tolerance of $10^{-3}$, and set $\xi$ such that $r^*$ is small enough that only the nearest periodic image contributes to the near field sum for each pair of points.

While the nearest image for near field calculations is over the sheared domain, $\V{r}$ and $r$ are computed using the Euclidean metric. We search for pairs of points closer than $r^*$ apart in log-linear time using a kD tree implemented in SciPy for a rectangular periodic cell. To adjust for the fact that the $\V{x}'$ coordinates are given on a non-orthogonal coordinate system, we bound the Euclidean distance between points in the unsheared coordinates by their ``distance'' in sheared coordinates, 
\begin{equation}
\label{eq:safety}
\norm{\V{r}} = \sqrt{\V{r}^T \V{r}} = \sqrt{\left(\V{r'}\right)^T \M{L}^{-T} \M{L}^{-1} \V{r'} }\leq \left(1+\frac{1}{2}\left(g^2+\sqrt{g^2(g^2+4)}\right)\right) \norm{\V{r}'}:=\psi \norm{\V{r}'}, 
\end{equation}
where we have used the maximum eigenvalue of $\M{L}^{-T} \M{L}^{-1}$ to bound the norm \cite{SpectralSD}. The factor $\psi$ defined in\ \eqref{eq:safety} can be thought of as a ``safety factor'' in the sense that points that are $r^*$ apart using the Euclidean metric in primed coordinates are at most $\psi r^* $ apart in physical space.  

For multiple interacting fibers, the velocity $\V{U}^{(\text{RPY})}\left(\piX{i}{p}\right)$ is obtained by summing both the far field and near field over fibers $j$ and points $q$. One inconvenience is that the PSE sum\ \eqref{eq:velalldirect} is over all pairs of points, including a point with itself ($\V{y}_i=\V{x}$ is automatically included in the sum). Specifically, the PSE sum\ \eqref{eq:velalldirect} will include interactions of a fiber with itself using the RPY kernel; this is not correct since those interactions should be computed by the SBT formula\ \eqref{eq:onefib}. Since we know \textit{a priori} that the velocity at target $\piX{i}{p}$ incorrectly includes the direct quadrature sum due to fiber $i$, we subtract the free space RPY kernel, defined in\ \eqref{eq:rpyknel}, from the RPY sum $\V{U}^{(\text{PSE})}\left(\piX{i}{p}\right)$ defined in\ \eqref{eq:velalldirect} to obtain the final Ewald sum 
\begin{equation}
\label{eq:velallsph}
\V{U}^{(\text{RPY})}\left(\piX{i}{p}\right) = \V{U}^{(\text{PSE})}\left(\piX{i}{p}\right) - \sum_q \KRPY{\piX{i}{p}}{\piX{i}{q}}{b} \ind{\V{f}}{i}_q w_q. 
\end{equation}

\subsection{Near fiber quadrature \label{sec:nearfibs}}
The Ewald splitting scheme\ \eqref{eq:velallsph} gives the velocity at all points $\piX{i}{p}$ due to all other points $\piX{j}{q}$, including target-fiber pairs where the direct quadrature scheme\ \eqref{eq:dirquad} is inaccurate. If the box size is sufficiently large, these inaccuracies only happen for the periodic image of the target $\piX{i}{p}$ that is closest to the fiber $\ind{\V{X}}{j}$. The other periodic images are handled correctly using the Ewald splitting scheme\ \eqref{eq:velallsph}, since they are sufficiently far from the fiber. In this section, we describe the special quadrature scheme we use to correct the velocity for fibers $j \neq i$ that are close to target $\piX{i}{p}$.

We need to compute the interaction velocity
\begin{equation}
\label{eq:intvel}
\V{v}\left(\V{x}\right)= \EPMI \int_0^{L} \Knel{\V{x},\V{X}(s),\rev{\dco}\left(\epsilon L \right)^2}\V{f}(s) \, ds,   
\end{equation}
where $\V{x}$ is a target point and $\V{X}(s)$ is the Chebyshev interpolant of the centerline of a fiber. 
There are several components in our scheme to compute the interaction velocities $\McNotd$ to a guaranteed tolerance regardless of how close the target point $\V{x}$ is to the centerline of fiber $\V{X}(s)$. We need to:
\begin{enumerate}
\item Understand how far $\V{x}$ can be from $\V{X}(s)$ for the direct quadrature\ \eqref{eq:dirquad} to remain sufficiently accurate. 
\item Obtain a reliable metric to compute or bound the minimum distance between the target $\V{x}$ and fiber $\V{X}(s)$,
\begin{equation}
d \eqd \min_s \norm{\V{x}-\V{X}(s)}. 
\end{equation}
Our procedure to compute $d$ will be different for $d/L = \mathcal{O}(1)$, when fibers are far apart and minimizing over a discrete set of nodes is sufficiently accurate, than for $d/L=\mathcal{O}(\epsilon)$, where it is more efficient to actually solve the continuous minimization problem.   
\item Use a special quadrature scheme to compute the integral for $d=\mathcal{O}(\epsilon)$. The special quadrature scheme will of course be more expensive than direct quadrature\ \eqref{eq:dirquad}, but less expensive than actually using the requisite number of direct quadrature points required to get the same accuracy. The scheme we use here is taken directly from \cite{barLud} and is based on extending the ideas used for the finite part integration in Section \ref{sec:MFP} to near-singular quadrature.
\item If \rev{$d < 4\aRPY$}, compute the closest point on fiber $\V{X}(s)$ to the target $\V{x}$ and denote it by $\V{X}\left(s^*\right)$. Use the distance metric for small $d$ to correctly blend the interaction velocity\ \eqref{eq:intvel} with the centerline velocity\ \eqref{eq:CLvel1} of the fiber at $\V{X}\left(s^*\right)$. 
\end{enumerate}

\begin{algorithm}
\centering
\includegraphics[width=\textwidth]{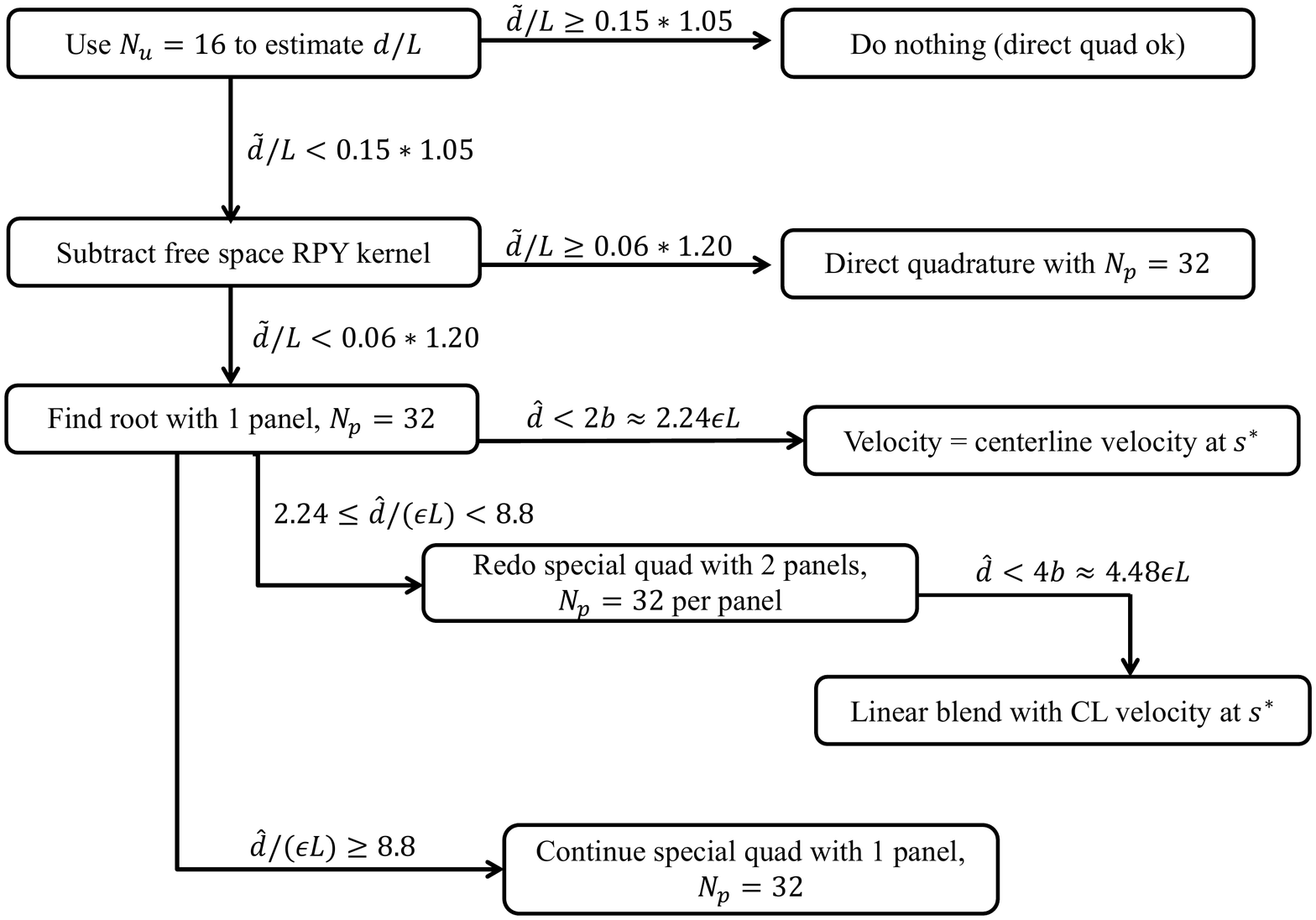}
\caption{\label{fig:algflow}The discrete operation $\McNotd=\MJFd\left(\V{x}, \V{X}\right)\V{f}$ for a fiber $\V{X}(s)$ and target $\V{x}$. Beginning with a target $\V{x}$ and $N$ Chebyshev points on the fiber centerline $\V{X}(s)$, we first estimate the distance $d/L$ by $\tilde{d}/L$, where $\tilde{d}$ is computed from the discrete minimization over $16$ uniform points. This determines whether the result from direct quadrature with $N$ (assumed to be at least 16) points is sufficiently accurate. If direct quadrature is not sufficiently accurate, we subtract the free space RPY kernel\ \eqref{eq:vRPY} from the velocity\ \eqref{eq:velallsph} and, if $N < 32$, redo the calculation with direct quadrature with $N_p=32$ points or special quadrature. If special quadrature is needed, we calculate the complex root $\eta^*$ of $\norm{\V{x}-\V{X}(\eta)}=0$, use the projection\ \eqref{eq:sstar} to find the closest point on the fiber, and compute the non-dimensional distance $\hat{d}/(\epsilon L)$ in\ \eqref{eq:dhat}. The distance $\hat{d}/(\epsilon L)$ determines whether to use special quadrature with 1 or 2 panels of 32 points each, and whether to use a linear combination with the fiber centerline velocity at $s^*$, given in\ \eqref{eq:CLvel1}. The algorithm is designed to ensure the calculation of $\McNotd$ to 3 digits of accuracy most of the time, and is specific to $\epsilon \approx 10^{-3}$. }
\end{algorithm}

To begin, we define an acceptable tolerance for the integrals. Since slender body theory itself is only accurate to $\mathcal{O}(\epsilon)$, it does not make sense to set a tolerance less than $\epsilon$. For actin filaments $\epsilon \approx 10^{-3}$, so we define the tolerance as $10^{-3}$ and set $\epsilon=10^{-3}$ in our accuracy tests. That is, our goal is to compute the interaction velocities $\McNotd$ to three digits of accuracy regardless of the distance between a target $\V{x}$ and fiber $\V{X}(s)$. Algorithm\ \ref{fig:algflow} gives the schematic flowchart of our method, the details of which are discussed next. While this method cannot guarantee 3 digits of accuracy, it does so for most target-fiber pairs of interest to us; see Appendix\ \ref{sec:nearfibacc} for numerical results.

\subsubsection{Distance where direct quadrature breaks down \label{sec:dbrkdwn}}
Our first goal is to determine when direct quadrature breaks down. To do this, we simply measured the accuracy of direct quadrature\ \eqref{eq:dirquad} with $N=16$ and $N=32$ Chebyshev nodes for randomly-generated pairs of fibers and targets. Specifically (more detail is given in Appendix\ \ref{sec:nearfibacc}), we generated 100 inextensible fibers with 16 nonzero Chebyshev coefficients decaying exponentially (in expectation) by four orders of magnitude. We placed 100 targets around each fiber a distance $d$ in the normal direction and computed the quadratures\ \eqref{eq:dirquad}. Measuring error with respect to a refined direct quadrature, we found that $N=16$ gives 3 digits of accuracy for all test cases when the non-dimensional distance $d/L \geq 0.15$. Likewise, direct quadrature with $N=32$ points gives 3 digits of accuracy when $d/L \geq 0.06$. 

\subsubsection{Estimating $d$ for $d/L = \mathcal{O}(1)$}
Since direct quadrature breaks down for $d/L < 0.15$ when $N=16$, we need to determine whether a target $\V{x}$ is indeed a distance less than $d/L=0.15$ from the fiber centerline $\V{X}(s)$ (the analogous statement holds for $d/L=0.06$, $N=32$). To do this quickly, we resample the fiber centerline at $N_u=16$ uniformly spaced points (in the arclength coordinate $s$) and perform a discrete minimization over the uniform fiber points to estimate the distance $d$. We denote this approximation by $\tilde{d}$. Using the same random set of fibers and targets as in Section\ \ref{sec:dbrkdwn}, we found the estimate of $d$ by $\tilde{d}$ has a relative error in the distance $|\tilde{d}-d|/d$ of at most $5\%$ for $d/L=0.15$ and $20\%$ for $d/L=0.06$.

Since direct quadrature breaks down at $d/L = 0.15$ for $N=16$, after accounting for the error in estimating $d$, we have that if $\tilde{d}/L \geq 0.15 \times 1.05$, direct quadrature with $N=16$ gives the integral to 3 digits. If $\tilde{d}/L < 0.15 \times 1.05$, it is possible that the direct quadrature\ \eqref{eq:dirquad} is not accurate to 3 digits, and so the many-body Ewald sum\ \eqref{eq:velallsph} contains a direct quadrature between $\V{x}$ and $\V{X}$ that is not accurate enough. We therefore subtract, for the target point $\V{x}$, this incorrect part of the sum, specifically the free space RPY kernel between the fiber $\V{X}(s)$ and target $\V{x}$, 
\begin{equation}
\label{eq:vRPY}
\V{v}^{\left(\text{RPY}\right)}(\V{x}) = \sum_{p=1}^N \KRPY{\V{x}}{\V{X}}{b} \, w_p . 
\end{equation}
We then recompute the integral\ \eqref{eq:intvel} to 3 digits using some other method. 

For $N=16$, our first resort to compute the integral\ \eqref{eq:intvel} is to sample the Chebyshev polynomial $\V{X}(s)$ at $32$ points and use direct quadrature\ \eqref{eq:dirquad} with $32$ points. Since direct quadrature with $32$ points gives 3 digits of accuracy when $d/L \geq 0.06$, including the error bounds we have that if $\tilde{d}/L \geq 0.06 \times 1.20$, the direct quadrature gives the integral\ \eqref{eq:intvel} to 3 digits. For $\tilde{d}/L < 0.06 \times 1.20$, we abandon direct quadrature and use a special quadrature routine. 

This initial step to upsample and integrate directly is designed to take care of most of the near target and fiber pairs without incurring a significant computational cost. Our empirical correlation between the distance $\tilde{d}$ and the direct quadrature error is less rigorous than the direct quadrature error estimates of \cite{barLud}. These estimates require information about the near singularity which, as discussed in the next section, must be computed using more expensive root finding.  

\subsubsection{Special quadrature.} The special quadrature routine is taken directly from \cite{barLud}. As in Section \ref{sec:MFP}, the underlying idea is to find the near singularity in the integrand, factor it out, expand what remains in a monomial basis, and compute integrals with the singularity and monomials analytically. This time, however, the integrand is not actually singular on the fiber centerline. By expanding the fiber representation to the complex plane, the nearby singularity can be found in the complex plane and the entire procedure of Section \ref{sec:MFP} can be repeated.

In more detail, the interaction velocity integral\ \eqref{eq:sbtother} can be written so that the numerator is smooth as $\V{x}$ approaches the centerline of $\V{X}(s)$. Starting from
\begin{gather}
\label{eq:MJFrwr}
\McNotd=\EPMI \int_0^{L}\left( \frac{\V{f}(s)}{\norm{\V{R}}}+
\frac{\left(\left(\V{R}\V{R}\right)+\rev{\dco}\left(\epsilon L\right)^2 \M{I}\right)\V{f}(s)}{\norm{\V{R}}^3} -3\rev{\dco}\left(\epsilon L \right)^2 \frac{\left(\V{R}\V{R}\right)\V{f}(s)}{\norm{\V{R}}^5}\right) \, ds, 
\end{gather}
where $\V{R}=\V{x}-\V{X}(s)$, we rescale $s$ by $\eta = -1+2s/L$. Then each of the terms in the integral\ \eqref{eq:MJFrwr} can be written in the form 
\begin{equation}
\label{eq:nsing}
\int_{-1}^1 \frac{\V{h}_m\left(\V{x},\eta \right)}{\norm{\V{x} -\V{X}(\eta)}^m}\, d\eta, 
\end{equation}
where $m=1,3,5$. For each $m$, $\V{h}_m\left(\V{x},\eta \right)$ is a density that depends on the target point $\V{x}$ and varies smoothly along the fiber $j$ arclength coordinate $\eta$. 

Now, the idea of \cite{barLud} is to extend the representation of $\V{X}(\eta)$ from $\eta \in [-1,1] \subset \mathbb{R}$ to the \textit{complex} plane $\mathbb{C}$ and compute the complex root of $\norm{\V{x} -\V{X}(\eta)}=0$. Because the centerline representation $\V{X}(\eta)$ is available as a Chebyshev series, it is simple to solve for the root via Newton iteration. We denote this root by $\eta^*$, i.e., $\norm{\V{x} -\V{X}(\eta^*)}=0$ with $\eta^* \in \mathbb{C}$. 

Once the root is known, the algorithm of \cite{barLud} finds the radius of the Bernstein ellipse associated with $\eta^*$. This radius then bounds the direct quadrature error for the integral\ \eqref{eq:MJFrwr}. If the upper bound on the direct quadrature error with $N=32$ points is less than $10^{-3}$, we proceed with direct quadrature. Otherwise, we use the special quadrature scheme of \cite{barLud}, which is the same as in the finite part integration in Section \ref{sec:MFP}. The complex singularity is factored out of the integrand, leaving a ``smooth'' function which can be expanded in a monomial basis. Integrals involving monomials multipled by the singularity are computed analytically, and an inner product of the monomial coefficients with the analytical integrands yields the approximation to $\McNotd$. 

The only difference from Section \ref{sec:MFP} is that the location of the singularity (complex root $\eta^*$) now depends on the fiber position $\V{X}$. This means that the roots and monomial coefficients must be computed at every time step using an $LU$ factorization of the Vandermonde matrix. Since the Vandermone matrix is a function of the nodes $s_p$ on the fiber, its $LU$ factorization is the same for all fibers and can be precomputed once at the start of a dynamic simulation. We refer the reader to \cite[Section~3]{barLud} for more details on this quadrature scheme. 

While error bounds exist for direct quadrature, the special quadrature scheme of \cite{barLud} does not provide error bounds or a method for selecting the fiber discretization (number of panels, points per panel, etc.). Because of this, we performed an empirical study on the same set of 100 fibers and targets as in Section\ \ref{sec:dbrkdwn}. The randomized testing described in Appendix\ \ref{sec:nearfibacc} showed that, for most fiber configurations of interest to us, 1 panel of 32 points is sufficient to compute the integral to 3 digits using special quadrature as long as the non-dimensional distance between the target and fiber is $d/(\epsilon L) > 8$. Otherwise, 2 panels of 32 points are required. Fibers with high curvature typically give the largest errors for a given discretization.  

\subsubsection{Estimating $d$ for $d/L = \mathcal{O}(\epsilon)$ \label{sec:closeds}} We still require a robust numerical procedure to determine when the the target point is too close to the cross section of the fiber (i.e., when \rev{$d < 4b =e^{3/2}\epsilon L$}). Our idea is to use the real part of the root $\eta^*$ as the closest arclength coordinate on the fiber to the target. We know the root $\eta^*$ solves $\norm{\V{x}-\V{X}(\eta^*)}=0$. It seems sensible, therefore, for the real part of the root to approximately minimize (over real $\eta$) $\norm{\V{x}-\V{X}(\eta)}$ when the root is close to the real line. 
We therefore define $s^*$, the closest point on the fiber to the target, from the complex root $\eta^*$ by removing the imaginary part of the root and rescaling, 
\begin{equation}
\label{eq:sstar}
s^*=\begin{cases} \frac{L}{2}(\text{Re}(\eta^*) +1) & -1 \leq \text{Re}(\eta^*) \leq 1\\[2 pt] 0 & \text{Re}(\eta^*) < -1 \\[2 pt] L &  \text{Re}(\eta^*) > 1 \end{cases}. 
\end{equation} 
The shortest distance from $\V{x}$ to the fiber can then be estimated as 
\begin{equation}
\label{eq:dhat}
\hat{d}\eqd \norm{\V{x}-\V{X}(s^*)}, 
\end{equation}
where the position $\V{X}(s^*)$ is computed by evaluating the Chebyshev interpolant at $s^*$. Our randomized tests showed that this estimate gives an error of at most $10\%$ for $d/(\epsilon L) \leq 8$, which is sufficiently accurate for our purpose. Because we use a point $s^*$ on the fiber centerline to estimate $\hat{d}$, this $10\%$ error is always an overestimation. For this reason we use 2 panels for special quadrature when $\hat{d}/(\epsilon L) \leq 8.8$.

Combining all of our steps, we obtain an algorithm to compute $\McNotd$ for all targets and fibers that is presented as a flowchart in Algorithm\ \ref{fig:algflow}. In Appendix \ref{sec:nearfibacc} we show that this quadrature scheme gives 3 digits of accuracy in the integrals\ \eqref{eq:intvel} with high probability.

\subsection{Temporal discretization \label{sec:tint}}
In this section, we discuss how we discretize the evolution equation\ \eqref{eq:abssys} in time. We use the notation of\ \eqref{eq:allfibNL} to split the mobility $\M{M}$ into an $\mathcal{O}(\log{\epsilon})$ local part $\MLD$ and $\mathcal{O}(1)$ non-local part $\MNL$. Our goal in this section is to develop a second-order temporal integrator with the properties that:
\begin{enumerate}
\item A minimum number of evaluations of the nonlocal hydrodynamics are needed per time step. 
\item Bending elasticity is treated implicitly. 
\item Any linear solves are block-diagonal, or fiber by fiber, so that the complexity of solving evolution equation\ \eqref{eq:abssys} scales linearly with the number of fibers and can be trivially parallelized. 
\end{enumerate}
Since we have separated the mobility matrix into the dominant $\mathcal{O}(\log{\epsilon})$ block diagonal local drag matrix $\MLD$ and the $\mathcal{O}(1)$ term that remains, we will begin by treating terms associated with $\MLD$ implicitly, thereby alleviating some of the stability restrictions associated with the bending force $\M{F}\V{X}$. This kind of approach gives almost unconditional stability as long as the density of fibers is small enough for the local drag term to dominate the fiber's motion. When the density of fibers is larger, we will treat the bending force implicitly in the $\MNL$ term as well; the resulting linear system can be solved approximately by a few iterations of GMRES \cite{gmres}.

To avoid nonlinear solves and still achieve second-order accuracy, we extrapolate values from previous time steps to the midpoint of the next time step and use these extrapolated values as the arguments for nonlinear functions (e.g. $\M{M}(\V{X})$ and $\M{K}(\V{X})$). More precisely, we define the extrapolated midpoint fiber positions as
\begin{equation}
\label{eq:Xstar}
\tdisc{\V{X}}{n+1/2, *} = \frac{3}{2}\tdisc{\V{X}}{n}-\frac{1}{2}\tdisc{\V{X}}{n-1}, 
\end{equation}
where we have used the notation $\tdisc{\V{X}}{n}$ to denote the fiber positions $\V{X}$ at the $n$th time step. As is our usual convention, we have not used fiber indices in the extrapolation\ \eqref{eq:Xstar} since it applies to every fiber independent of the others.

\subsubsection{Semi-implicit method for dilute suspensions \label{sec:bdiag}}
To discretize the evolution equation\ \eqref{eq:abssys} and principle of virtual work\ \eqref{eq:nowork} in time, we split the mobility matrix into $\M{M} = \MLD+\MNL$. Since the elastic force density $\M{F}\V{X}$ involves fourth derivatives, it must be treated implicitly to maintain stability as the number of Chebyshev grid points $N$ increases. For dilute suspensions, we assume that, since $\MLD$ is the dominant $\mathcal{O}(\log{\epsilon})$ contribution at each point, only treating the term $\MLD\M{F}\V{X}$ implicitly will still give improved stability. 

A second-order, semi-implicit discretization of the evolution equation\ \eqref{eq:abssys} begins by solving the saddle-point system
\begin{gather}
\label{eq:CNnL}
\tdisc{\MLD}{n+1/2,*} \left(\tdisc{\V{\lambda}}{n+1/2}+\frac{1}{2}\M{F}\left(\tdisc{\V{X}}{n} + \tdisc{\V{X}}{n+1, *}\right)\right)\\[2 pt] \nonumber 
+\tdisc{\MNL}{n+1/2,*} \left(\tdisc{\V{\lambda}}{n+1/2, *}+\M{F}\tdisc{\V{X}}{n+1/2, *}\right)+\V{u}_0\left(\tdisc{\V{X}}{n+1/2,*}\right)= \tdisc{\M{K}}{n+1/2, *} \tdisc{\V{\alpha}}{n+1/2},\\[2 pt]
\nonumber
\tdisc{\M{K}}{n+1/2, *}^* \tdisc{\V{\lambda}}{n+1/2} =\rev{\V{0}} 
\end{gather}
for $\tdisc{\V{\lambda}}{n+1/2}$ and $\tdisc{\V{\alpha}}{n+1/2}$, where the notation $\tdisc{\MLD}{n+1/2,*}$ means $\MLD\left(\tdisc{\V{X}}{n+1/2,*}\right)$ (and likewise for $\MNL, \M{K}$, and $\M{K}^*$). To obtain a second-order block-diagonal system, we extrapolate previous $\V{\lambda}$ values to the midpoint of the next time step, 
\begin{equation}
\tdisc{\V{\lambda}}{n+1/2, *} = 2\tdisc{\V{\lambda}}{n-1/2}-\tdisc{\V{\lambda}}{n-3/2}. 
\end{equation} 
\rev{When the nonlocal mobiility in the block diagonal discretization\ \eqref{eq:CNnL} is applied to $\tdisc{\V{\lambda}}{n+1/2, *}$ and $\M{F}\tdisc{\V{X}}{n+1/2, *}$, the total force might be nonzero because of discretization errors in $\M{F}$. In this case, we manually set the total force to zero in the PSE method of Section\ \ref{sec:ewald}.} We also introduce the approximation
\begin{gather}
\label{eq:Xn1star}
\tdisc{\V{X}}{n+1, *} = \tdisc{\V{X}}n + \D t \tdisc{\M{K}}{n+1/2, *}\tdisc{\V{\alpha}}{n+1/2}
\end{gather}
to make\ \eqref{eq:CNnL} a linear system in $\tdisc{\V{\alpha}}{n+1/2}$ and $\tdisc{\V{\lambda}}{n+1/2}$. Since $\tdisc{\V{X}}{n+1}$ will actually be computed by rotating the tangent vectors and integrating the result (see Section \ref{sec:Xsupdate}), the update\ \eqref{eq:Xn1star} is a second-order approximation to the actual $\tdisc{\V{X}}{n+1}$. Nevertheless, it is a sufficient approximation to give the same stability properties as if the actual $\tdisc{\V{X}}{n+1}$ were included in a nonlinear saddle-point system.


By substituting the approximation $\tdisc{\V{X}}{n+1, *}$ in\ \eqref{eq:Xn1star} into saddle-point system\ \eqref{eq:CNnL}, we obtain the following saddle-point linear system to solve at every time step, 
\begin{gather}
\label{eq:spttemp}
    \begin{pmatrix}
    -\MLD & \M{K}-\frac{\D t}{2}\MLD \M{F}\M{K} \\[2 pt]
    \M{K}^* & \rev{\M{0}}
    \end{pmatrix}_{n+1/2,*}
    \begin{pmatrix} 
    \tdisc{\V{\lambda}}{n+1/2}\\[2 pt]
    \tdisc{\V{\alpha}}{n+1/2}
    \end{pmatrix} =  \\[2 pt]
    \nonumber 
    \begin{pmatrix} 
    \tdisc{\MLD}{n+1/2,*}\M{F}\tdisc{\V{X}}{n} + \tdisc{\MNL}{n+1/2,*}\left(\tdisc{\V{\lambda}}{n+1/2,*}+\M{F}\tdisc{\V{X}}{n+1/2, *}\right) + \V{u}_0\left(\tdisc{\V{X}}{n+1/2, *}\right)\\[2 pt]
   \rev{\V{0}}
   \end{pmatrix}.
\end{gather}
This system can be solved fiber by fiber, since all of the matrices $\MLD, \M{K}$, and $\M{K}^*$ on the left hand side of\ \eqref{eq:spttemp} are block diagonal. System\ \eqref{eq:spttemp} is not invertible in general because the representation $\bm{K}\bm{\alpha}$ is not necessarily unique. To see this, suppose that $\V{n}_j$ is a degree $N-1$ polynomial. Then the inextensible velocity\ \eqref{eq:du} could be zero at all the nodes without $\V{\alpha}$ being identically zero. We therefore use the least squares solution for $\V{\alpha}$ while emphasizing that $\V{\alpha}$ itself has no significance; only $\M{K}\V{\alpha}$ has physical meaning. 

\subsubsection{Implicit method for denser suspensions \label{sec:gmres}}
In the case when the fibers are packed densely enough to make the temporal discretization\ \eqref{eq:spttemp} unstable, we treat the bending force in the nonlocal hydrodynamics implicitly and use GMRES to solve for $\tdisc{\V{\lambda}}{n+1/2}$ and $\tdisc{\V{\alpha}}{n+1/2}$. The new linear system of equations is 
\begin{gather}
\label{eq:gressys}
\tdisc{\MLD}{n+1/2,*} \left(\tdisc{\V{\lambda}}{n+1/2}+\frac{1}{2}\M{F}\left(\tdisc{\V{X}}{n} + \tdisc{\V{X}}{n+1, *}\right)\right)\\[2 pt] \nonumber 
+\tdisc{\MNL}{n+1/2,*} \left(\tdisc{\V{\lambda}}{n+1/2}+\frac{1}{2}\M{F}\left(\tdisc{\V{X}}{n} + \tdisc{\V{X}}{n+1, *}\right)\right)+\V{u}_0\left(\tdisc{\V{X}}{n+1/2,*}\right)= \tdisc{\M{K}}{n+1/2, *} \tdisc{\V{\alpha}}{n+1/2},\\[2 pt]
\nonumber
\tdisc{\M{K}}{n+1/2, *}^* \tdisc{\V{\lambda}}{n+1/2} =\rev{\V{0}}.
\end{gather}

Now, let us denote the solutions of the block diagonal system\ \eqref{eq:spttemp} by $\tdisc{\widetilde{\V{\lambda}}}{n+1/2}$ and $\tdisc{\widetilde{\V{\alpha}}}{n+1/2}$. By subtracting the fully implicit system\ \eqref{eq:gressys} from the locally implicit system\ \eqref{eq:CNnL}, we obtain the residual form of the saddle-point system
\begin{gather}
\label{eq:resid}
    \begin{pmatrix}
    -\left(\MLD +\MNL\right) & \M{K}-\frac{\D t}{2}\left(\MLD+\MNL\right) \M{F}\M{K} \\[2 pt]
    \M{K}^* & \rev{\M{0}}
    \end{pmatrix}_{n+1/2,*}
    \begin{pmatrix} 
    \Delta \tdisc{\V{\lambda}}{n+1/2}\\[2 pt]
    \Delta \tdisc{\V{\alpha}}{n+1/2}
    \end{pmatrix} =  \\[2 pt]
    \nonumber 
    \begin{pmatrix} 
    \tdisc{\MNL}{n+1/2,*}\left(\M{F}\left(\tdisc{\V{X}}{n}+\frac{\Delta t}{2}\tdisc{\M{K}}{n+1/2, *}\tdisc{\widetilde{\V{\alpha}}}{n+1/2}-\tdisc{\V{X}}{n+1/2,*}\right)+\tdisc{\widetilde{\V{\lambda}}}{n+1/2}-\tdisc{\V{\lambda}}{n+1/2,*}\right)\\[2 pt]
   \V{0}
   \end{pmatrix}
\end{gather}
to be solved using GMRES for the perturbations
\begin{gather}
\label{eq:perturb}
\Delta \tdisc{\V{\lambda}}{n+1/2} = \tdisc{\V{\lambda}}{n+1/2} - \tdisc{\widetilde{\V{\lambda}}}{n+1/2} \quad \text{and}\\[2 pt] \nonumber
\Delta \tdisc{\V{\alpha}}{n+1/2} = \tdisc{\V{\alpha}}{n+1/2} - \tdisc{\widetilde{\V{\alpha}}}{n+1/2}.
\end{gather}
The right hand side of system\ \eqref{eq:resid} is zero to second order in $\Delta t$. We therefore expect the perturbations to be $\mathcal{O}(\Delta t^2)$. While these perturbations have no impact on the temporal accuracy of our scheme, obtaining a good approximation to them is vital for stability. In Section \ref{sec:gresiters}, we quantify empirically (for $\epsilon=10^{-3}$) how many GMRES iterations are enough to obtain unconditional stability. Note that smaller values of $\epsilon$ require fewer GMRES iterations for stability since local drag is more dominant for smaller $\epsilon$, and vice versa for larger $\epsilon$.

To solve system\ \eqref{eq:resid} rapidly with GMRES, we use the block diagonal preconditioner 
\begin{equation}
\label{eq:prec}
\tdisc{\M{P}}{n+1/2}=\begin{pmatrix} -\MLD & \M{K}-\frac{\D t}{2}\MLD \M{F}\M{K} \\[2 pt]
    \M{K}^* & \rev{\M{0}}
    \end{pmatrix}^{-1}_{n+1/2,*}
\end{equation}
that appears in\ \eqref{eq:spttemp}. Since $\MLD$ dominates over $\MNL$, this preconditioner is effective, with increased effectiveness for smaller $\epsilon$. 

Our overall scheme to solve for $\tdisc{\V{\lambda}}{n+1/2}$ and $\tdisc{\V{\alpha}}{n+1/2}$ can be summarized as follows:
\begin{enumerate}
\item Solve the block diagonal system\ \eqref{eq:spttemp} for $\tdisc{\V{\lambda}}{n+1/2}$ and $\tdisc{\V{\alpha}}{n+1/2}$. If the fiber suspension is sufficiently dilute (see Section\ \ref{sec:gresiters}), continue to the next time step.
\item Otherwise, set $\tdisc{\widetilde{\V{\lambda}}}{n+1/2}$ and $\tdisc{\widetilde{\V{\alpha}}}{n+1/2}$ to be the solutions of the block diagonal system\ \eqref{eq:spttemp}, run a few iterations of GMRES to solve the residual system\ \eqref{eq:resid}, and update $\tdisc{\V{\lambda}}{n+1/2} $ and $\tdisc{\V{\alpha}}{n+1/2} $ using\ \eqref{eq:perturb}.
\end{enumerate}
To initialize, in the first and second time steps (for $n=0,1$), we solve system\ \eqref{eq:gressys} by converging GMRES with a relative residual tolerance of $10^{-6}$. When $n=0$, we set $\tdisc{\V{X}}{n+1/2, *} = \tdisc{\V{X}}{n}$.

\subsubsection{Updating $\Xs$ and $\V{X}$ \label{sec:Xsupdate}}
Once we have computed $\tdisc{\V{\alpha}}{n+1/2}$, we use a discrete form of the tangent vector rotation\ \eqref{eq:omegadef} to update the tangent vectors. This is done fiber by fiber, and so here we use $\V{X}$ to stand for a single fiber, rather than the entire collection of positions. To avoid double subscripts, in a slight abuse of notation we superscript the time step index $n$ in this section. 
\renewcommand{\tdisc}[2]{#1^{#2}} 

Our goal is to rotate the set of tangent vectors $\tdisc{\Xs}{n}$ on the unit sphere by the (axis-angle) rotations $\tdisc{\V{\Omega}}{n+1/2} =  \V{\Omega}\left(\tdisc{\Xs}{n+1/2,*},\tdisc{\V{\alpha}}{n+1/2}\right)$. To do this in a stable way, we cannot use $\tdisc{\V{\alpha}}{n+1/2}$ directly in the computation, since the kinematic coefficients $\tdisc{\V{\alpha}}{n+1/2}$ do not have physical meaning \rev{and take values which are sensitive to discretization and ill-conditioning}. Since \rev{the velocity} $\tdisc{\M{K}}{n+1/2}\tdisc{\V{\alpha}}{n+1/2}$ has a physical meaning and is \rev{less sensitive to numerical artifacts}, we compute
\begin{equation}
\label{eq:Omega}
\tdisc{\V{\Omega}}{n+1/2} = \tdisc{\Xs}{n+1/2, *} \times \M{D}_N\tdisc{\M{K}}{n+1/2, *}\tdisc{\V{\alpha}}{n+1/2}
\end{equation} 
on an upsampled grid, where $\M{D}_N$ is the Chebyshev differentiation matrix on the $N$ point grid. To do this, we upsample $\tdisc{\Xs}{n+1/2,*}$ and the derivative $\M{D}_N\tdisc{\M{K}}{n+1/2, *}\tdisc{\V{\alpha}}{n+1/2}$ to a grid of size $2N$ and do the cross product. We then downsample the result to obtain $\tdisc{\V{\Omega}}{n+1/2}_p$ for $p=1, \dots N$. 

Once $\tdisc{\V{\Omega}}{n+1/2}$ is known, we use the Rodrigues rotation formula \cite{rodrigues1840lois} to compute the rotated tangent vectors. Letting $\Omega = \norm{\V{\Omega}}$ and $\widehat{\V{\Omega}}=\V{\Omega}/\Omega$, we compute the rotated tangent $\Xs$ at each node $p$ by
\begin{gather}
\label{eq:updateXs}
\tdisc{\Xs_p}{n+1} = \Xs_p^n \cos{\left(\tdisc{\Omega_p}{n+1/2}\D t\right)} + \left(\tdisc{\widehat{\V{\Omega}}_p}{n+1/2} \times \tdisc{\Xs_p}{n} \right)  \sin{\left(\tdisc{\Omega_p}{n+1/2}\D t\right)}+\\[2 pt] \nonumber \tdisc{\widehat{\V{\Omega}}_p}{n+1/2}\left(\tdisc{\widehat{\V{\Omega}}_p}{n+1/2} \cdot \tdisc{\Xs_p}{n}\right)\left(1-\cos{\left(\tdisc{\Omega_p}{n+1/2}\D t\right)}\right).
\end{gather}
We then compute $\tdisc{\V{X}}{n+1}$ from $\tdisc{\Xs}{n+1}$ via Chebyshev integration. Specifically, we compute the Chebyshev series coefficients of $\tdisc{\Xs}{n+1}$, apply the spectral integration matrix \cite{greengard1991spectral} to compute the Chebyshev series of $\tdisc{\V{X}}{n+1}$, then evaluate this series at the nodes on the $N$ point grid. To fix the integration constant, on each fiber we set the position at the first node
\begin{equation}
\tdisc{\V{X}_1}{n+1} = \tdisc{\V{X}_1}{n+1, *},  
\end{equation}
where $\tdisc{\V{X}}{n+1,*}$ is defined in\ \eqref{eq:Xn1star}.

\section{Numerical tests \label{sec:tests}}
In this section, we validate each component of our method and demonstrate the method's spatial and temporal accuracy. We study spatio-temporal accuracy in Section\ \ref{sec:spttemp} with simple examples of four falling fibers in free space and three fibers in periodic shear flow. Section \ref{sec:gresiters} gives our most important result for computational complexity: the number of hydrodynamic evaluations per time step required to stably evolve the dynamics of a fiber suspension is at most five. By varying the fiber number density and bending modulus, we show that one block diagonal solve combined with at most three iterations of GMRES per time step are needed to maintain stability.\footnote{The extra hydrodynamic evaluation to give a total of five comes in the conversion to residual form\ \eqref{eq:resid}.}

In general, we will use an $L^2$ function norm to compute the differences between configurations throughout this section. Given two fiber configurations, we evaluate the Chebyshev interpolant of each on a 1000 point type 2 Chebyshev grid and calculate the discrete $L^2$ error using Clenshaw-Curtis quadrature. Whenever there are multiple fibers, we compute the error on the first fiber $\ind{\V{X}}{1}$ unless otherwise specified. 

Since fibers in shear flow are our primary interest, some of our examples will use shear flows. The general form of a time-oscillatory shear flow is given by
\begin{equation}
\label{eq:shear}
\V{u}_0(\V{x},t)=\dot{\gamma_0} \cos{(\omega t)}(y,0,0). 
\end{equation}
The corresponding strain is given by $g(t)=\left(\dot{\gamma}_0/\omega\right) \sin{(\omega t)}$ for $\omega > 0$ and $g(t)=\dot{\gamma}_0 t$ for $\omega = 0$. The code and input files for the examples of fibers in triply periodic shear flow presented here are available for download at \url{https://github.com/stochasticHydroTools/SlenderBody}. 

\subsection{Spatio-temporal accuracy \label{sec:spttemp}}
In this section, we study spatio-temporal convergence for two examples where nonlocal hydrodynamics has a nontrivial impact on the fiber trajectories. We first verify our method for fibers in gravity by comparing the results to those obtained using the method prescribed in \cite{ehssan17} and in the process show improved robustness and temporal accuracy. We then demonstrate second-order temporal convergence and spectral spatial accuracy for a periodic three fiber system in shear flow. We also show that, if stable, our block diagonal solver\ \eqref{eq:spttemp} with one hydrodynamic evaluation per time step is the most efficient way to resolve the dynamics to a given tolerance. 

\subsubsection{Comparison to strong formulation \label{sec:falling}}
Our first goal is to validate our weak formulation of inextensibility by numerically comparing to the strong formulation. To do this, we consider four fibers centered around a circle of radius $d=0.2$. The fibers have initial tangent vector $\Xs(t=0)=(0,0,1)$ (they are aligned in the $z$ direction) and positions $\ind{\bm{X}}{1}(t=0)=(d,0,s-1)$, $\ind{\bm{X}}{2}(t=0)=(0,d,s-1)$, $\ind{\bm{X}}{3}(t=0)=(-d,0,s-1)$, and $\ind{\bm{X}}{4}(t=0)=(0,-d,s-1)$, where $0 \leq s \leq L=2$. For simplicity, we set $\mu=\kappa=1$, and $\epsilon=10^{-3}$. For this test only, we use fibers with ellipsodial cross sections and set the local drag coefficient $c(s) = -\ln{\left(\epsilon^2\right)}$. We simulate until $t=0.25$.

\begin{figure}
\centering
\subfigure[Trajectory]{\label{fig:fallingfibs}
\includegraphics[width=0.3\textwidth]{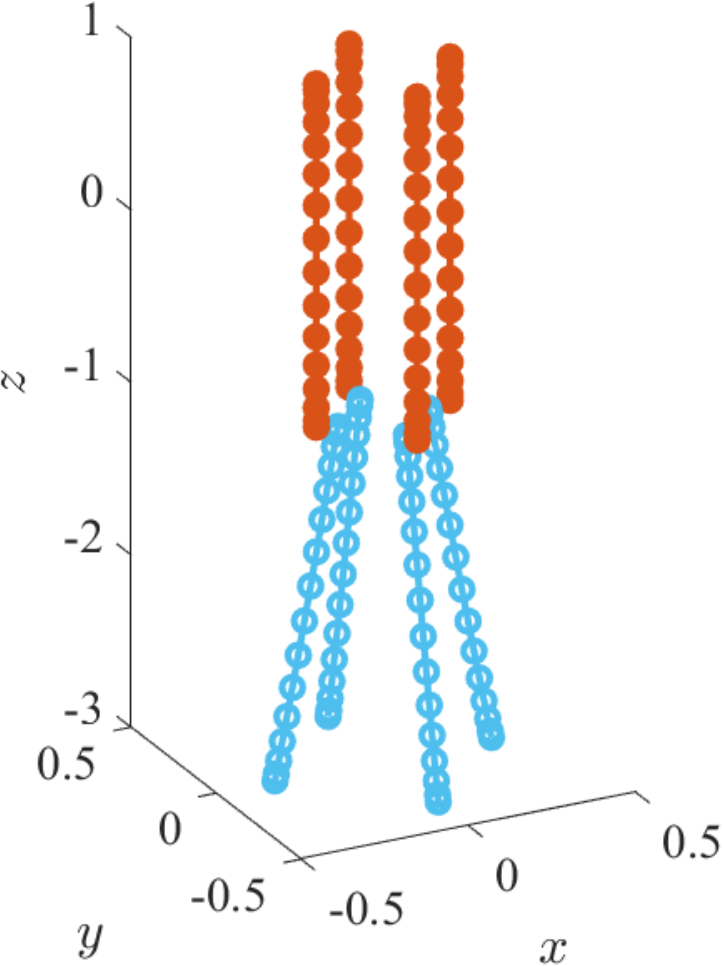}}
\subfigure[Convergence]{\label{fig:stfall}
\includegraphics[width=0.5\textwidth]{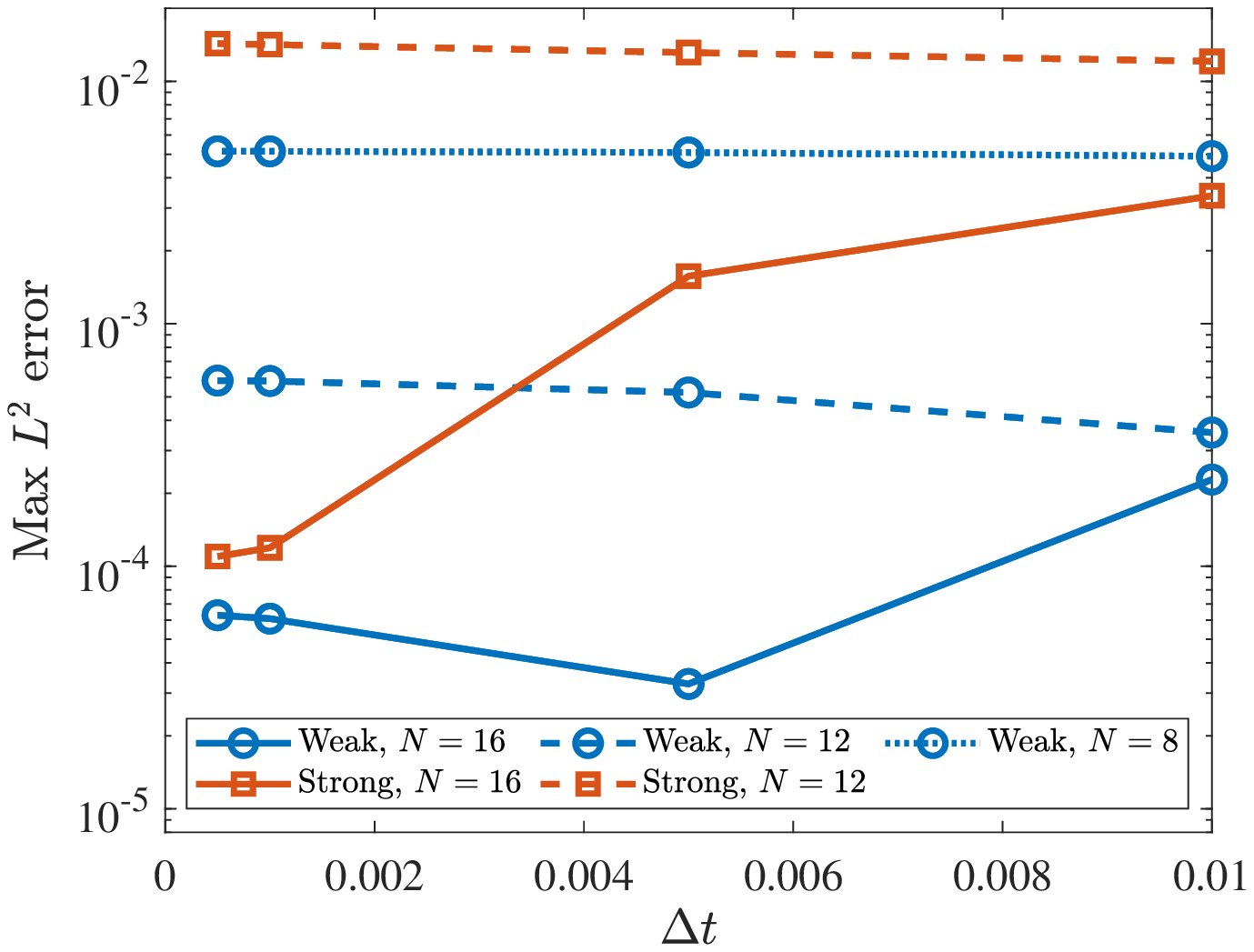}}
\caption{Four fibers in gravity. (a) Initial (filled orange) and final (light blue) configurations of the fibers. There are $N=16$ points on each fiber. (b) Spatio-temporal convergence of our weak formulation (circles) compared to the strong formulation of \cite{ehssan17} (squares). In both cases, the exact solution is a trajectory with $N=24$, $\Delta t = 5 \times 10^{-4}$. For small $\Delta t$, the spatial error dominates and the spatio-temporal error saturates.}
\end{figure}

Each fiber has a uniform gravitational force density $\V{f}^g = (0,0,-5)$ placed on it. In the absence of nonlocal interactions (i.e., if $\M{M}=\MLD$), the fibers fall straight downward. When the interactions between fibers are included, however, the fibers influence each other and have an $x$ and $y$ direction to their motion. Figure\ \ref{fig:fallingfibs} shows the initial and final ($t=0.25$) configurations of the fibers in this test.

Our goal here is to verify our results by comparing to results obtained using the method of \cite{ehssan17}. Because the method of \cite{ehssan17} uses regularization for the finite part integral (and we do not), we drop the finite part integral in this calculation and only include local drag and cross fiber interactions in the mobility. The fibers are sufficiently far from each other that the dipole term in the kernel $\V{S}_D$, which has coefficient of 0 in \cite{ehssan17}, has a very small effect on the result (two orders of magnitude less than our smallest spatio-temporal error). We use free space boundary conditions (no periodicity) and compute all nonlocal integrals by direct quadrature\ \eqref{eq:dirquad}, without any upsampling. We use the block diagonal solves\ \eqref{eq:spttemp} for temporal integration and do \textit{not} perform any GMRES iterations. 

We first verify that our results match those of the strong formulation \cite{ehssan17} when the spatial and temporal discretizations are well-refined. Considering $N=24$ and $\Delta t = 5 \times 10^{-4}$ in both algorithms, we obtain a maximal $L^2$ difference of $1.8 \times 10^{-4}$ in the position of the first fiber, which, as we show in Fig.\ \ref{fig:stfall}, is on the order of magnitude of the discretization error. 

It is instructive to compare the spatio-temporal error between the two algorithms. We define the ``exact'' solution for both algorithms to be a trajectory with $N=24$ and  $\Delta t = 5 \times 10^{-4}$. Figure\ \ref{fig:stfall} shows the maximum $L^2$ errors over the time interval $[0,0.25]$ for both algorithms with different spatial and temporal discretizations. For small $\Delta t$, the spatial error dominates and the spatio-temporal error saturates. We observe that our weak formulation outperforms the strong formulation of \cite{ehssan17} in two ways. First, for coarse discretizations (e.g. $N=12$, dashed lines in Fig.\ \ref{fig:stfall}), our saturated spatial error is more than an order of magnitude lower than the saturated spatial error of \cite{ehssan17}. This is likely because the line tension equation of \cite{ehssan17} has larger aliasing errors for coarser spatial discretizations, and because of our improved treatment of the free fiber boundary conditions using rectangular spectral collocation. Secondly, our errors saturate at a much larger time step size than those of \cite{ehssan17}. For example, when $N=16$ our error saturates for $\Delta t = 5 \times 10^{-3}$, whereas the error from \cite{ehssan17} does not saturate until $\Delta t = 1 \times 10^{-3}$. This occurs because our temporal integrator is second-order accurate. This simple example demonstrates the improved spatial accuracy of our new weak formulation over the strong one, and the improved accuracy of our temporal discretization, even in the absence of GMRES iterations. 

\subsubsection{Spatio-temporal convergence}
\begin{figure}
\centering
\subfigure[Problem setup]{\label{fig:threeshear}
\centering
\includegraphics[width=\textwidth, trim = {1cm 0cm 1cm 1cm}, clip]{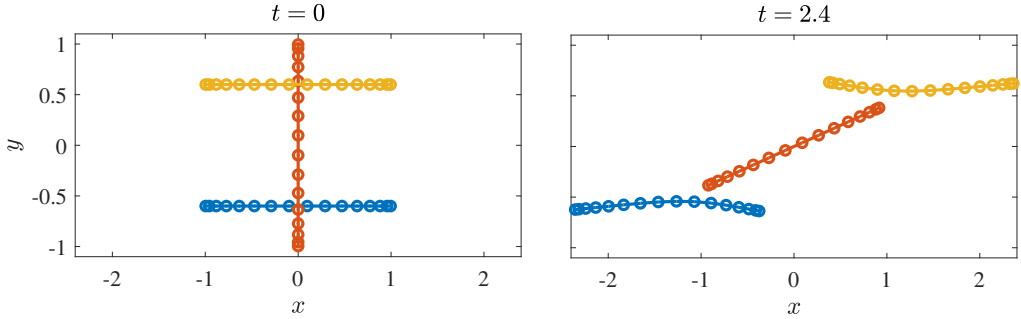}}
\subfigure[Temporal error]{\label{fig:threeconvtemp}
\includegraphics[width=0.4\textwidth]{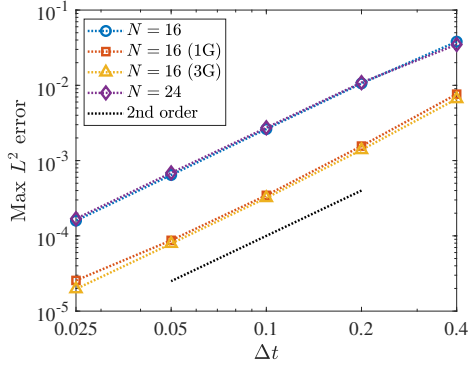}}
\subfigure[Spatio-temporal error]{\label{fig:threeconv}
\includegraphics[width=0.4\textwidth]{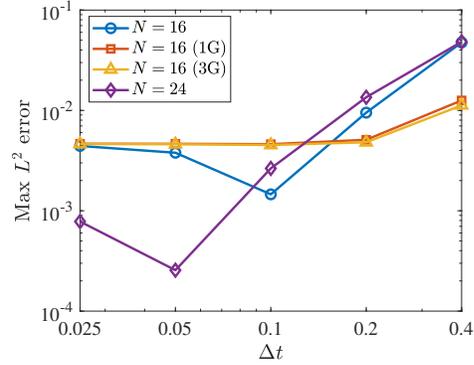}}
\caption{\label{fig:three} Three fibers in shear flow. (a) Fiber configurations at $t=0$ and $t=2.4$. (b) Second-order temporal convergence. We measure the maximum $L^2$ error over time in the first fiber position using successive refinements and observe second-order temporal convergence for block diagonal solves with $N=16$ (blue circles) and $N=24$ (purple diamonds). We also see second-order convergence, and reduced temporal errors, for both 1 (red squares) and 3 (yellow triangles) GMRES iterations. (c) The spatio-temporal errors (measured against a more accurate solution with $N=32$, $\Delta t=0.00125$) are shown for $N=16$ and $N=24$.}
\end{figure}

We next verify the temporal convergence of our algorithm for periodic sheared domains by choosing an example where interactions between the fibers contribute significantly to the dynamics. We consider three sheared fibers with $L=2$ and $\ind{\V{X}}{1}(s)=(s-1,-0.6,\rev{-0.04})$, $\ind{\V{X}}{2}(s)=(0,s-1,0)$, and $\ind{\V{X}}{3}(s)=(s-1, 0.6,\rev{0.06})$. As shown in Figure\ \ref{fig:threeshear}, this corresponds to an ``I'' shaped initial configuration of the fibers, with the fibers stacked in the $z$ direction.

We set the periodic domain length $L_d=2.4$, the Ewald parameter $\xi=3$, and set $\mu=1$, $\epsilon=10^{-3}$, and $\kappa=0.01$. We use a constant shear flow\ \eqref{eq:shear} with $\gamma_0=1$ and $\omega=0$. Because of the small bending rigidity, the fibers deform from their straight configurations in a shear flow. In this example, two of the fibers are initially aligned with the $x$ direction. Without nonlocal hydrodynamics, they would stay aligned with the $x$ direction and simply translate. When nonlocal interactions are included, the flows generated by the middle fiber $\ind{\V{X}}{2}$ induce deformations of the top and bottom fibers. This is evident in Fig.\ \ref{fig:threeshear}, which shows the final fiber positions at $t=2.4$. 

To quantify the temporal convergence, we fix $N=16$ and simulate from $t=0$ to $t=2.4$ with $\Delta t = 0.4,0.2,0.1,0.05,$ and $0.025$ using the block diagonal solver\ \eqref{eq:spttemp} (no GMRES iterations) and successive refinements to measure the error. Figure\ \ref{fig:threeconvtemp} shows that we obtain second-order temporal convergence for the block diagonal solver\ \eqref{eq:spttemp} (blue circles) and for the GMRES system\ \eqref{eq:resid} with 1 iteration (red squares) and 3 iterations (yellow triangles). The temporal error is about an order of magnitude smaller when we perform one GMRES iteration in addition to the block diagonal solvers. Note, however, that this comes at the cost of two additional nonlocal hydrodynamic evaluations (one to convert to residual form, one in the GMRES iteration). 

For spatio-temporal convergence, we simulate a refined trajectory with $N=32$ and $\Delta t = 0.00125$ and compute the maximum $L^2$ errors in trajectories with $N=16$ and $24$. As shown in Fig.\ \ref{fig:threeconv}, increasing the number of points by 8 decreases the spatio-temporal error by a factor of 4 (compare blue circles with purple diamonds), which is consistent with spectral spatial accuracy. 

Performing one GMRES iteration with $\Delta t = 0.2$ approximately matches the spatial and temporal errors and costs a total of three hydrodynamic evaluations per $0.2$ units of time. Running the block diagonal solver with $\Delta t = 0.1$ matches the temporal error with the spatial error and costs two hydrodynamic evaluations per $0.2$ units of time. We therefore conclude from Figure\ \ref{fig:threeconv} that the most efficient way to obtain the maximum accuracy for a given spatial resolution is to run the block diagonal solver with the smaller time step size, assuming it is stable. 

\subsection{Stability \label{sec:gresiters}}
Because our block-diagonal semi-implicit temporal discretization (BDSI) described in Section \ref{sec:bdiag} treats the bending force explicitly in the nonlocal term, it will become unstable when the fiber suspension is too concentrated and the cumulative effect of nonlocal hydrodynamics is comparable to that of local drag. In this case, we switch to the GMRES solver described in Section\ \ref{sec:gmres}. As discussed there, since the block diagonal solver is already second-order accurate, the perturbations to $\V{\alpha}$ and $\V{\lambda}$ that come from the residual GMRES solve\ \eqref{eq:resid} do not impact the overall temporal accuracy (see Fig.\ \ref{fig:three}), but do impact stability. For this reason, we run only a fixed number of GMRES iterations until we obtain stability. Our goal in this section is to determine an upper bound on the number of required GMRES iterations. 

To do this, we consider a suspension of $F=1000$ fibers and vary the density of fibers by changing the periodic domain length $L_d$. If $f=F/L_d^3$ is the number density of fibers and $L$ is the length of a fiber, a dimensionless density $fL^3 < 1$ is considered a dilute fiber suspension, while a semi-dilute suspension is one with $r fL^2 = \epsilon fL^3 \ll 1$, and a semi-concentrated one has $ \epsilon fL^3 = \mathcal{O}(1)$ \cite{fibsusps}. Here we explore the semi-dilute and semi-concentrated regimes and derive empirical bounds on how many GMRES iterations are required to maintain stability for a variety of bending moduli. Our conclusion is that at most five nonlocal hydrodynamic evaluations are sufficient to maintain stability, even for semi-concentrated suspensions.  

We simulate $F=1000$ initially straight fibers of length $L=2$ and radius $r=2\times 10^{-3}$ (so that $\epsilon = 1 \times 10^{-3}$), and use $N=16$ points per fiber. In the oscillatory shear flow\ \eqref{eq:shear}, we set $\omega=2\pi$ and $\dot{\gamma}_0=\omega/10$, so that the maximum strain is $g=0.1$ and the time for one cycle is 1. We expect to need at least 20 time steps per cycle to obtain reasonable accuracy, so we set $\Delta t = 0.05$, although in reality we find that smaller time step sizes are needed to accurately resolve the dynamics of dense suspensions. Since we find that changing the frequency $\omega$ has a negligible impact on the results, we non-dimensionalize $\Delta t$ by the bending timescale $\tau = 8\pi \mu L^4/\left(\ln{\left(\epsilon^{-2}\right)}\kappa\right)$, where here we use $\mu=1$. 

\begin{figure}
\centering
\includegraphics[width=0.75\textwidth]{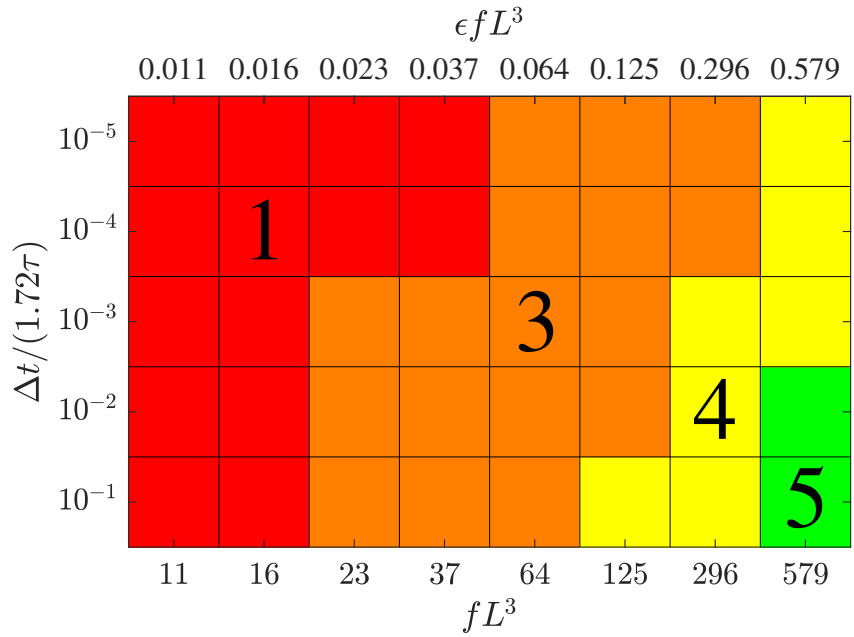}
\caption{\label{fig:itercounts}Number of hydrodynamic evaluations needed for stability, indicated by color. We consider a system of $F=1000$ fibers of length $L=2$, vary the periodic domain length $L_d$ and fiber stiffness $\kappa$, and show the number of hydrodynamic evaluations needed for stability for each set of parameters. One hydrodynamic evaluation occurs in the block diagonal solver\ \eqref{eq:spttemp}, another in the conversion to residual form\ \eqref{eq:resid}, and one evaluation occurs per GMRES iteration. We report fiber density in units of $fL^3=FL^3/L_d^3$ on the $x$ axis and fiber stiffness in units of $\Delta t / (1.72 \tau)$ on the $y$ axis, where the elastic timescale $\tau = 8\pi \mu L^4/\left(\ln{\left(\epsilon^{-2}\right)}\kappa\right)$.}
\end{figure}

We simulate 5 cycles of motion, until $t=5$. Figure\ \ref{fig:itercounts} shows the number of hydrodynamic evaluations required for stability for a given bending modulus and fiber density. The fiber number density $fL^3$ is reported on the bottom $x$ axis and on the top $x$ axis we give $\epsilon f L^3$. For semi-dilute suspensions ($\epsilon fL^3 \approx 0.01$), we see that BDSI is stable, i.e., only a single hydrodynamic evaluation is needed for stability (red squares in Fig.\ \ref{fig:itercounts}). For a fixed number density, we first see instabilities for BDSI for smaller $\tau$, so that stiffer fibers require more GMRES iterations for stability. As the fiber suspension becomes semi-concentrated ($\epsilon fL^3 \geq 0.1$), we see that we need at least three hydrodynamic evaluations regardless of the fiber stiffness. For our stiffest and densest suspensions ($\epsilon fL^3 \approx 0.5$), we need at most five hydrodynamic evaluations per time step to obtain stable dynamics. For comparison, Nazockdast et al.\ report $9$ to $16$ GMRES iterations for a system with a similar number of fibers \cite[Table~1]{ehssan17}. 

We caution that these evaluation counts are the minimum number needed for \textit{stability}. For dense fiber suspensions, $\Delta t =0.05$ might be too large to obtain reasonable accuracy, since fibers that are in close contact are subject to nearly nonsmooth velocity fields that come from the near singular velocity kernel\ \eqref{eq:dintvel} and the possible combination with the fiber centerline velocity\ \eqref{eq:CLvel1}. In practice, we find this causes fibers that are in close contact to oscillate around each other when $\Delta t$ is too large. A smaller $\Delta t$ can resolve these issues, and we know from Fig.\ \ref{fig:itercounts} that the number of iterations required for stability drops with the time step size.

\section{Application: cross-linked actin mesh \label{sec:CLs}}
The cell cytoskeleton is a dynamic network of cross-linked actin filaments and myosin motors that allows cells to migrate, divide, and adapt to new environments \cite{alberts, ahmed2015dynamic}. A number of experimental \cite{gardel2004scaling, janmey1994mechanical, kasza2010actin, lieleg2009cytoskeletal} and computational \cite{ma2018structural, head2003deformation} studies have shown that the viscoelastic rheology of actin networks comes from specialized cross-linking proteins dynamically binding and unbinding to actin fibers, with the rates of binding and unbinding determining the ratio of viscous to elastic behavior. When the cross-linkers (CLs) are permanently attached, the network has traditionally been viewed as purely elastic, while instant unbinding of CLs has been seen as pure viscous behavior \cite{ahmed2015dynamic}. The reality is more nuanced than this, since the network is embedded in an underlying fluid which contributes to the viscous modulus of the network and affects the movement of the filaments between binding and unbinding events.

To our knowledge, there has been no systematic study of the contribution of the underlying fluid to the dynamics and viscoelastic properties of cross-linked actin networks. We defer a full study of this for the future; here we develop the cross linking model and provide some initial results. In particular, we will leave transient cross-linker dynamics for a future study and consider the special case of a permanently cross-linked network, which could model, for example, a network of actin fibers cross-linked by scruin proteins \cite{gardel2004scaling}. Since our fibers are represented by Chebyshev interpolants $\V{X}(s)$, we seek a continuum force density on the fiber due to cross linking. This force density must be smooth relative to the discretization to preserve the spectral accuracy of our algorithm. 

We begin in Section\ \ref{sec:CLmod} by presenting our model of a cross-linker as an elastic spring between the fibers. Although an elastic spring model might not be appropriate for some cross-linking proteins, for example short and stiff rods like $\alpha$-actinin \cite{meyer1990bundling}, it can model longer and more flexible CLs like filamin \cite{weihing1985filamins}. In Section\ \ref{sec:rheo}, we then discuss how we compute rheological information from our model to facilitate comparison with experiments. In Section\ \ref{sec:results}, we also report the sensitivity of this information to the number of collocation points $N$, width $\sigma$ of the CL Gaussian smoothing function, \rev{and local drag regularization length scale $\delta L$}. We conclude this section with our results for a permanently cross-linked actin network. We consider a fixed ratio of twelve CLs per fiber and study the viscous and elastic behavior of the network using both local drag and fully nonlocal hydrodynamics to evolve the system. Our conclusion is that there exists a critical time scale $\tau_c$ on which the network relaxes to a dynamic steady state under oscillatory shear flow. For shear frequencies $\omega \ll \tau_c^{-1}$, the behavior is primarily elastic and dominated by the quasi-steady state of the network. For $\omega \gg \tau_c^{-1}$, dynamics, including hydrodynamics, matter and we see more viscous behavior.

\subsection{Cross-linker model \label{sec:CLmod}}
Suppose we have two fibers, $\ind{\V{X}}{i}$ and $\ind{\V{X}}{j}$, and that a CL connects two fibers by attaching to arclength coordinate $s^*_i$ on fiber $i$ and $s^*_j$ on fiber $j$, where these coordinates are not necessarily Chebyshev points. We define the force density due to the CL at arclength coordinate $s_p$ on fiber $i$ as
\begin{equation}
\label{eq:fcl1}
\ind{\V{f}}{\text{CL}, i}_p\left( \V{X}\right)=  -K_c  \left(1-\frac{\ell}{\norm{\ind{\V{X}}{i}(s^*_i)-\ind{\V{X}}{j}(s^*_j)}}\right)\delta_h(s_p-s_{i}^*) \sum_{q=1}^N \left(\ind{\V{X}}{i}\left(s_p\right)-\ind{\V{X}}{j}(s_q)\right) \delta_h(s_q-s_{j}^*) w_q,  
\end{equation}
where $K_c$ is the spring constant for the CL (units force/length), $\ell$ is the rest length, and $\delta_h$ is a Gaussian smoothing function replacing a Dirac delta function. The CL force density\ \eqref{eq:fcl1} links the point $s_p$ on fiber $i$ to every point on fiber $j$, with a weight related to the distance on fiber $j$ between the anchor coordinate $s_j^*$ and Chebsyshev point $s_q$ by the Gaussian function $\delta_h$. The prefactor outside of the sum is zero when the anchor points are exactly length $\ell$ apart. If the two anchor points are farther than $\ell$ apart, the force between them is attractive; otherwise it is repulsive. 

The force density\ \eqref{eq:fcl1} exerts no net force or torque on the system both in continuum (replace the sum in\ \eqref{eq:fcl1} by an integral) and discretely. Specifically,
\begin{gather}
\label{eq:nof}
\sum_{p=1}^N \ind{\V{f}}{i, \text{CL}}_p w_p + \sum_{q=1}^N \ind{\V{f}}{j, \text{CL}}_q w_q = 0 \\[2 pt]
\label{eq:notorq}
 \sum_{p=1}^N \left(\piX{i}{p} \times \ind{\V{f}}{i, \text{CL}}_p\right) w_p + \sum_{q=1}^N \left(\piX{j}{q} \times \ind{\V{f}}{j, \text{CL}}_q\right) w_q = 0.
\end{gather}
These identities, which imply that each pair of cross-linked fibers is force-and-torque-free, hold regardles of the form of $\delta_h$. 

For a given fiber discretization, we choose $\delta_h$ so that the forcing is smooth in the Chebyshev basis. We consider $\delta_h$ to be a Gaussian density of the form
\begin{equation}
\label{eq:dh}
\delta_h(r) = \frac{1}{Z} \exp{\left(-\frac{r^2}{2\sigma^2}\right)}, 
\end{equation}
where $\sigma$ is a parameter that controls the smoothness and spread of $\delta_h$. The factor $Z$ is a normalization factor that ensures $\delta_h$ discretely integrates to 1 along the fiber length, i.e., $\sum_{p} \delta_h\left(s_p-s_i^*\right)w_p=1$. Far from the endpoints, $Z=\sqrt{2\pi \sigma^2}$, but if the CL is bound to the end of the fiber some of the Gaussian weight might be truncated. For a given $N$, we choose $\sigma$ to be the minimum value that gives 3 digits of accuracy in $\V{f}$, where the error is measured relative to a refined $\V{f}$ computed on a 1000 point Chebyshev grid. \revp{Note that small values of $\sigma$ lead to instabilities if $N$ is not large enough to resolve $\delta_h$. }For $N=16$, we use $\sigma/L=0.1$; we will study the influence of $\sigma$ on physical observables numerically in Section\ \ref{sec:Nsig}.

\renewcommand{\tdisc}[2]{#1_{#2}} 
We modify the BDSI temporal integrator\ \eqref{eq:CNnL} to treat the cross linker forces in an \textit{explicit} second-order fashion, 
\begin{gather}
\label{eq:CNnLCLs}
\tdisc{\MLD}{n+1/2,*} \left(\tdisc{\V{\lambda}}{n+1/2}+\frac{1}{2}\M{F}\left(\tdisc{\V{X}}{n} + \tdisc{\V{X}}{n+1, *}\right)+\tdisc{\ind{\V{f}}{\text{CL}}}{n+1/2,*}\right)\\[2 pt] \nonumber 
+\tdisc{\MNL}{n+1/2,*} \left(\tdisc{\V{\lambda}}{n+1/2, *}+\M{F}\tdisc{\V{X}}{n+1/2, *}+\tdisc{\ind{\V{f}}{\text{CL}}}{n+1/2,*}\right)+\V{u}_0\left(\tdisc{\V{X}}{n+1/2,*}\right)= \tdisc{\M{K}}{n+1/2, *} \tdisc{\V{\alpha}}{n+1/2}, 
\end{gather}
where $\tdisc{\ind{\V{f}}{\text{CL}}}{n+1/2,*}=\ind{\V{f}}{\text{CL}}\left(\tdisc{\V{X}}{n+1/2,*}\right)$ and $\tdisc{\V{X}}{n+1/2,*}$ is the extrapolation\ \eqref{eq:Xstar}. As before, when the block diagonal semi-implicit temporal integrator\ \eqref{eq:spttemp} is unstable, we switch to the residual system for GMRES\ \eqref{eq:resid}, which is unchanged by the CLs. The time step size $\Delta t $ for the CL discretization\ \eqref{eq:CNnLCLs} is limited by stability for cases when the fibers are more flexible than the CLs; we leave implicit treatment of CLs for future work.

\subsection{Rheological experiments \label{sec:rheo}}
Experimental studies on actin networks \cite{janmey1994mechanical, kasza2010actin} typically report the viscous and elastic moduli. These quantities can be computed from the system stress tensor $\M{\sigma}$, which in turn can be computed from the force densities on each fiber's centerline. In this section, we briefly lay out the minimum details needed for the calculation of the viscous and elastic moduli. 

Recalling the time-oscillatory shear flow\ \eqref{eq:shear}, the only nonzero component of the rate of strain tensor is constant in space and is given by 
\begin{equation}
\dot{\gamma}_{21}(t)= \frac{\partial u^x_0}{\partial y} = \dot{\gamma}_0\cos{(\omega t)}, 
\end{equation}
and the relevant component of the \textit{strain} tensor is therefore
\begin{equation}
\gamma_{21}(t) = \int_0^t \dot{\gamma}_{21}(t') \, dt' =  \frac{\dot{\gamma}_0}{\omega}\sin{(\omega t)}:=\gamma_0\sin{(\omega t)}, 
\end{equation}
where we have defined $\gamma_0$ as the maximum strain in the system. We define the bulk elastic ($G'$) and viscous modulus ($G''$) from $\gamma_0$ by \cite{morrison2001understanding}
\begin{equation}
\label{eq:stress}
\frac{\sigma_{21}}{\gamma_0} = G' \sin{(\omega t)}+ G'' \cos{(\omega t)}. 
\end{equation}
Notice that the elastic modulus $G'$ gives the part of the stress that is in phase with the strain, and the viscous modulus $G''$ gives the part of the stress that is in phase with the rate of strain. 

The stress tensor itself can be decomposed into a part coming from the background fluid and a part coming from the internal fiber stresses, 
\begin{equation}
\label{eq:strdef}
\sigma_{21}= \sigma_{21}^{(\mu)} + \sigma_{21}^{(f)}= \mu  \frac{\partial u^x_0}{\partial y}  + \sigma_{21}^{(f)}. 
\end{equation}
For pure viscous fluid, the fiber contribution to the stress is zero and the stress is given entirely by the viscous stress tensor, 
\begin{equation}
\frac{\sigma_{21}^{(\mu)}}{\gamma_0} = \frac{\omega \mu}{\dot{\gamma_0}}\left(\dot{\gamma}_0\cos{(\omega t)}\right) = \omega \mu \cos{(\omega t)}.
\end{equation}
Thus the viscous modulus due to the fluid is $G''=\omega \mu$. 

The stress due to the fibers depends on the force the fibers exert on the fluid. Because\ \eqref{eq:nof} shows that the total force exerted by the cross-linker on the pair of fibers it connects is zero, we use Batchelor's formula \cite{batch70} for the volume-averaged stress due to the fibers and CLs. In the slender limit (i.e., the case when the surface area force density is constant on fiber cross sections), the bulk stress due to the fibers and CLs in a volume $V$ is given by $\M{\sigma}^{(f)}=\M{\sigma}^{(i)} +  \M{\sigma}^{\text{(CL)}}$, where
\begin{gather}
\label{eq:Batch}
\M{\sigma}^{(i)} = -\frac{1}{V}\left(\sum_{i=1}^F \sum_{p=1}^N \left(\piX{i}{p}\left(\piThree{\V{\lambda}}{i}{p} + \left(\M{F}\ind{\V{X}}{i}\right)_p \right)w_p\right)\right),\\[2 pt]
\label{eq:sigCL}
\M{\sigma}^{\text{(CL)}} = -\frac{1}{V} \left(\sum_{c=1}^{C} \sum_{p=1}^N \left(\piX{c_1}{p}\piThree{\V{f}}{\text{CL}, c_1}{p} w_p +  \piX{c_2}{p}\piThree{\V{f}}{\text{CL}, c_2}{p}w_p\right)\right),
\end{gather}
where the double sum for $\M{\sigma}^{(i)}$ is over points on fibers and the double sum for $\M{\sigma}^{\text{(CL)}}$ is over points on pairs of cross-linked fibers. In the cross-linker stress, $C$ is the number of CLs and the notation $c_1$ and $c_2$ means that cross-linker $c$ links fibers $c_1$ and $c_2$. The sums must be separated since different periodic images of $\ind{\V{X}}{i}$ could be involved for different cross linkers. In $\M{\sigma}^{(CL)}$, the positions $\piX{c_1}{p}$ and $\piX{c_2}{p}$ must be the periodic images of the two fibers that are connected by the CL. Because\ \eqref{eq:notorq} shows that the total torque exerted by the cross-linker on the pair of fibers it connects is zero, we expect the stress tensor\ \eqref{eq:Batch} to be symmetric to spectral accuracy, since the constraint force $\V{\lambda}$ and elastic forces $\M{F}\V{X}$ exert exactly zero torque in continuum but not discretely \cite{batch70}. 

\renewcommand{\tdisc}[2]{#1^{#2}} 
We discretize the stress tensor at the midpoint of each time step. Specifically, we substitute $\tdisc{\V{X}}{n+1/2,*}$ for $\V{X}$ and $\tdisc{\V{\lambda}}{n+1/2}$ for $\V{\lambda}$ in the internal fiber stress $\M{\sigma}^{(i)}$, and evaluate the cross-linking forces at the extrapolated midpoint $\tdisc{\V{X}}{n+1/2,*}$ in $\M{\sigma}^{\text{(CL)}}$. Assuming that the final time $T$ is an integer multiple of the period $2\pi/\omega$, we then compute the bulk moduli by discretizing the integrals
\begin{equation}
\label{eq:Gpdp}
G'  = \frac{2}{\gamma_0 T} \int_0^{T}  \sigma_{21}  \sin{(\omega t)} \, dt \qquad  G'' = \frac{2}{\gamma_0 T} \int_0^{T} \sigma_{21} \cos{(\omega t)} \, dt. 
\end{equation}
by the midpoint rule.

\subsection{Results \label{sec:results}}
We now consider a permanently cross-linked network of $F$ filaments, fix the physical parameters of the network, and analyze how the viscoelastic behavior depends on the frequency of oscillations. We set the fiber length $L=2$ $\mu$m with aspect ratio $\epsilon =10^{-3}$ and bending stiffness $\kappa=0.01$ $\text{pN}\cdot\mu$m$^2$, and use a box of size $L_d=4$ $\mu$m. For the fluid viscosity, we use $\mu =1$ Pa$\cdot$s. The CLs have rest length $\ell=0.5$ $\mu$m and spring constant $K_c=1$ pN/$\mu$m. To bind the CLs, at $t=0$ we resample each fiber centerline to 16 uniformly separated points. We then randomly select a pair of these points, $\ind{\V{X}}{i}\left(s_i^*\right)$ on filament $i$ and $\ind{\V{X}}{j}(s_j^*)$ on filament $j$, where $i \neq j$. If the two selected points are initially separated by a distance less than $\ell$, we bind a CL connecting fibers $i$ and $j$, with one end centered on $\ind{\V{X}}{i}\left(s_i^*\right)$ and the other end centered on $\ind{\V{X}}{i}(s_j^*)$. We continue this process until $12F$ CLs have been attached.\footnote{Since there are only 16 sites on each fiber and one CL takes up two sites, this requires that we allow more than one CL to bind to a specific site.} Our use of such a large number of CLs effectively makes the network into a single interconnected cluster, so that our periodic domain can be viewed as a sample of a bulk interconnected fiber gel.

\subsubsection{Effect of changing $N$, $\sigma$, \rev{and $\delta$} \label{sec:Nsig}}
For our experiments in the rest of this section, we will typically use $N=16$ Chebyshev points per fiber, \rev{cross-linker standard deviation $\sigma/L=0.1$, and local drag regularization parameter $\delta=0.1$}. Here we test the effect of changing $N$, $\sigma$, \rev{and $\delta$} from these baseline parameters by considering a set of $F=100$ straight fibers and 1200 CLs. In this section, we are concerned only with the spatio-temporal accuracy of our CL formulation, and not necessarily the actual values of the viscous and elastic moduli, and so for this test we start measuring stress at $t=0$, despite the fact that the network could be far from a steady state. Since we have already tested the spatio-temporal accuracy of our hydrodynamics in Section\ \ref{sec:tests}, we run here with local drag only. We use $\omega=2\pi$ and $\dot{\gamma}_0=0.2\pi$ and run until $T=6$ seconds (6 periods). To obtain a set of refined trajectories, we use $\Delta t = 0.005$ for $N=16$, $\Delta t = 0.0025$ for $N=24$, and $\Delta t = 0.001$ for $N=32$. In addition to $\sigma/L=0.10$ for all discretizations, we also measure stress for $N=24$ with $\sigma/L=0.07$ and $N=32$ with $\sigma/L=0.05$. 

\begin{figure}
\centering
\subfigure[Constant $\sigma/L=0.1$ and constant $\delta=0.1$]{\label{fig:Nsigstress}
\includegraphics[width=0.45\textwidth]{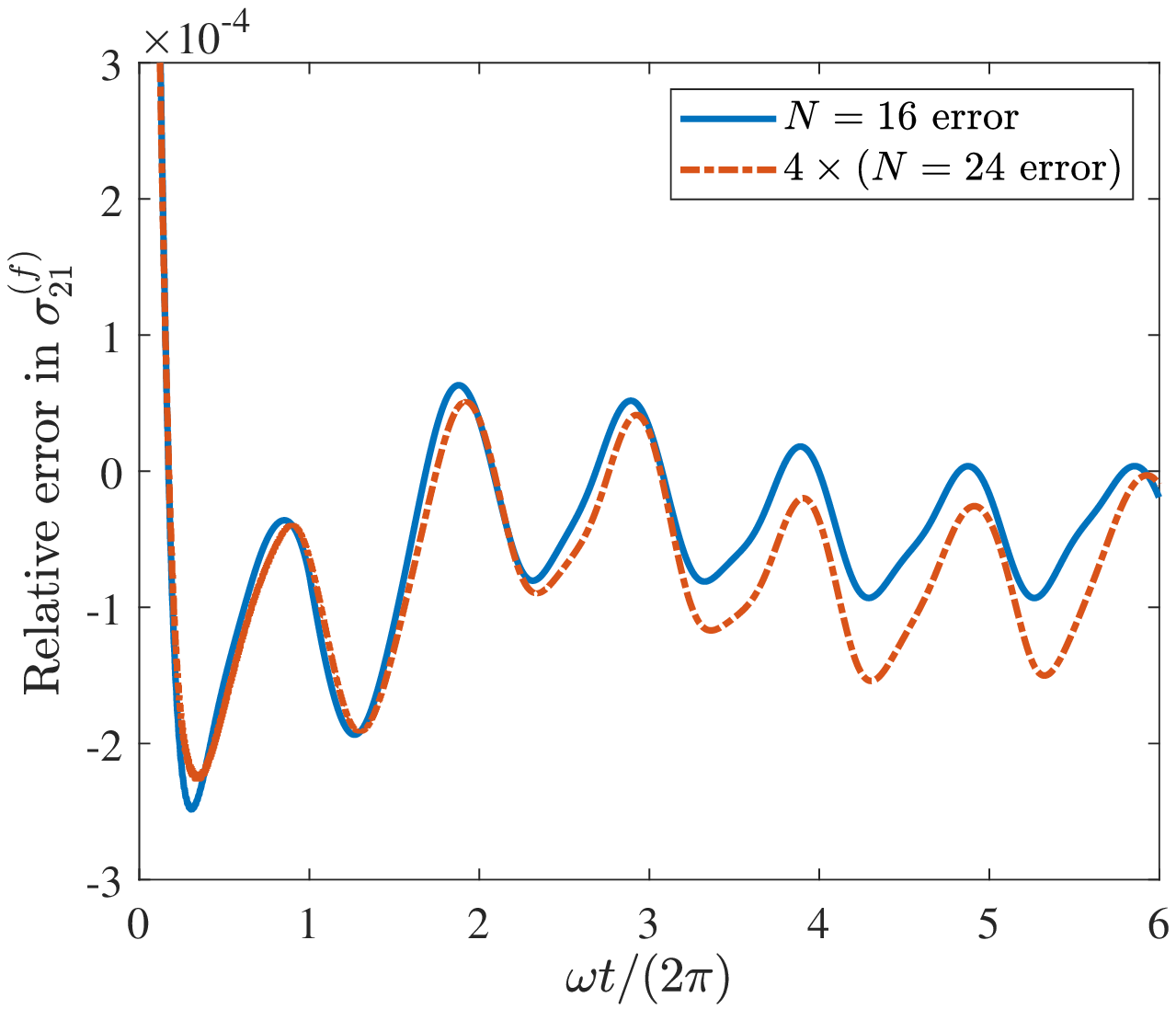}}
\subfigure[Varying $\sigma/L$, constant $\delta=0.1$]{\label{fig:NsigstressVar}
\includegraphics[width=0.45\textwidth]{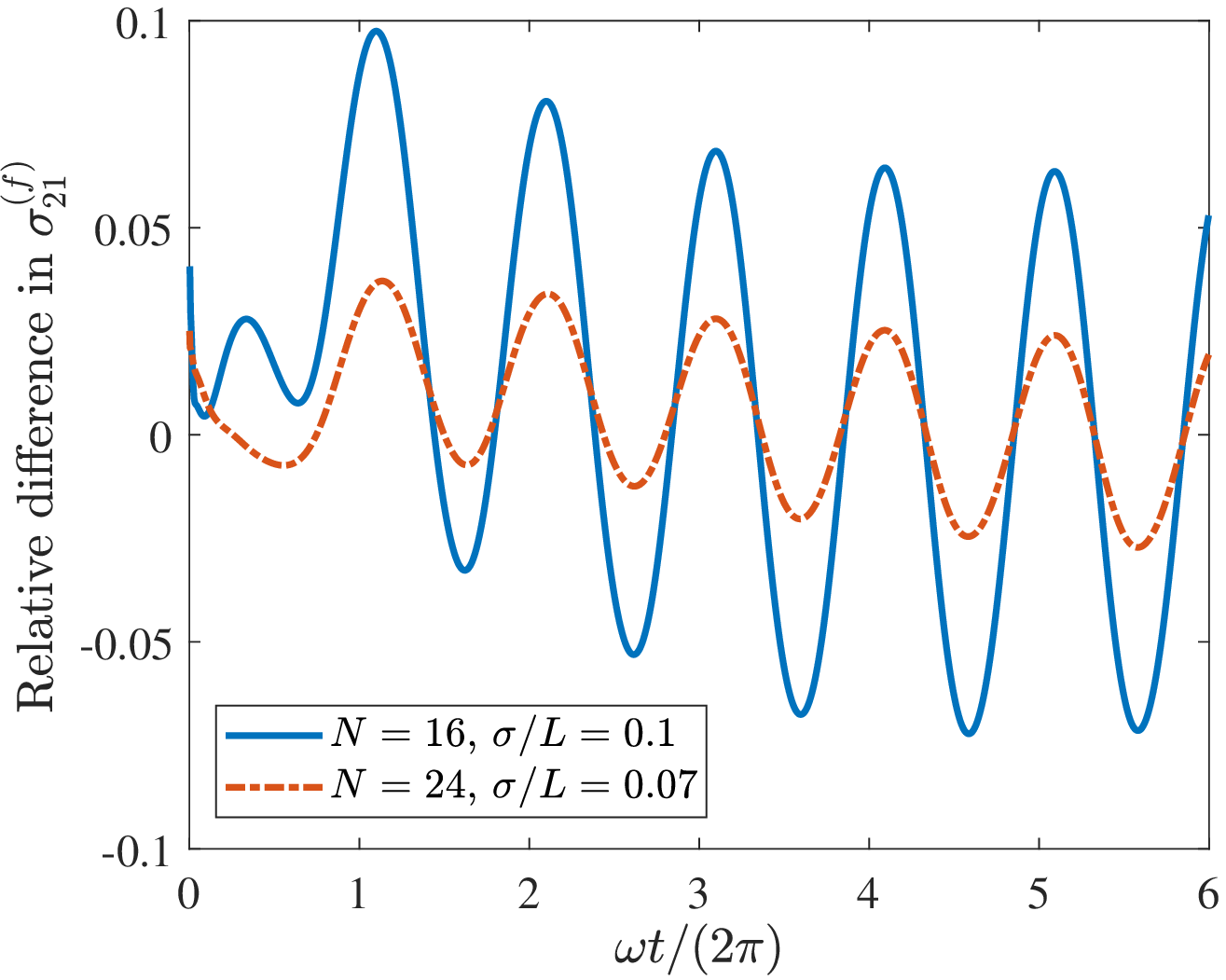}}
\subfigure[Varying $\delta$, constant $\sigma/L=0.1$]{\label{fig:DeltaStressVar}
\includegraphics[width=0.45\textwidth]{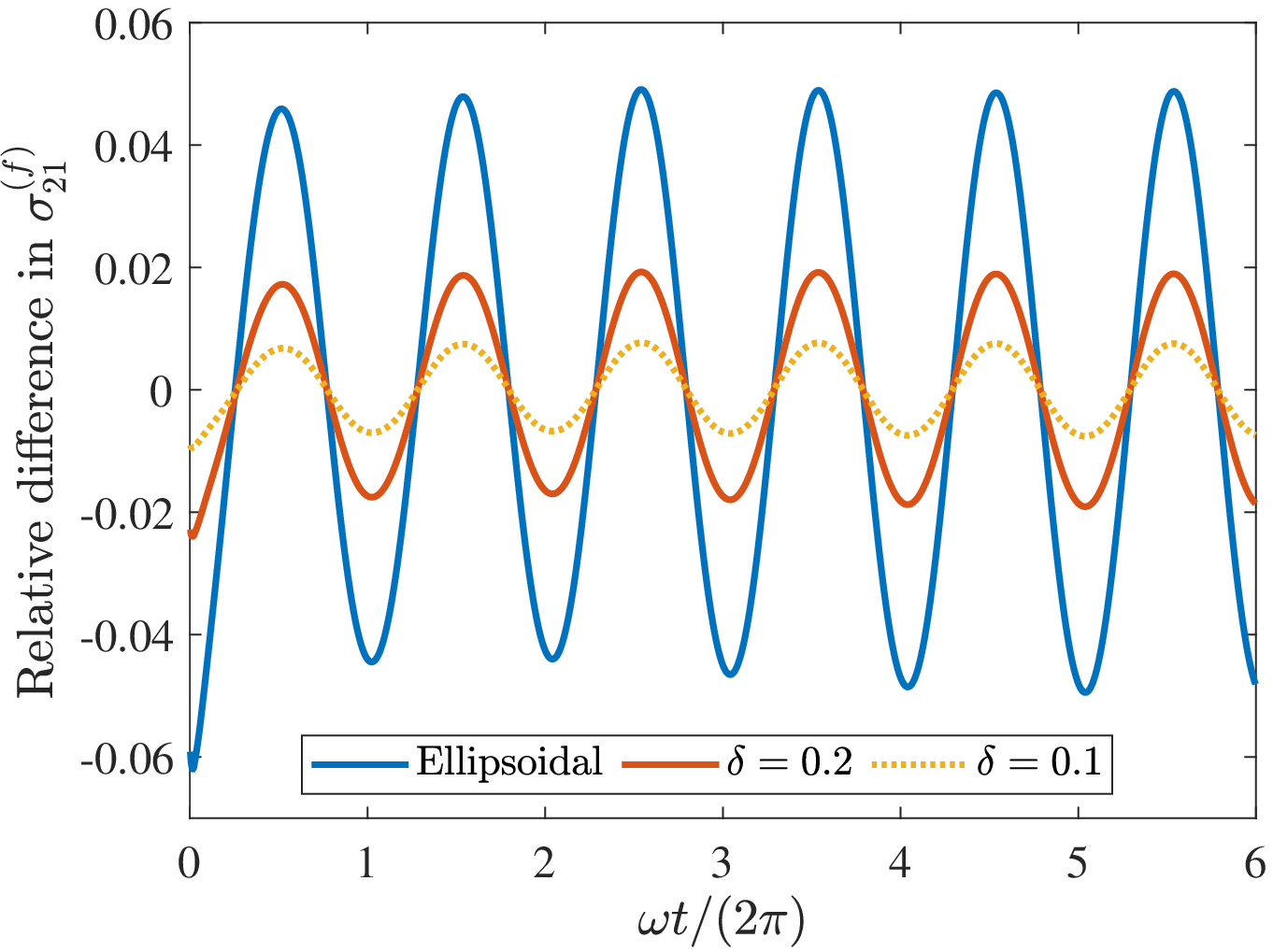}}
\caption{Differences in the fiber stress $\sigma_{21}^{(f)}$ in a suspension of $F=100$ fibers and 1200 CLs with \rev{(a) constant cross-linker standard deviation $\sigma/L$ and local drag regularization parameter $\delta$, (b) varying $\sigma/L$, and (c) varying $\delta$}. Normalization is the maximum absolute value of the stress for the reference solution in all cases. (a) For a reference solution with $N=32$ and $\sigma/L=0.1$, we see (after initial transients) rapid convergence of the stress consistent with second-order convergence in time. (b) For a reference solution with $N=32$ and $\sigma/L=0.05$, we see a $10\%$ difference in stress using $N=16$ and $\sigma/L=0.1$, and an approximately $3\%$ difference for $N=24$ and $\sigma/L=0.07$. \rev{(c) For a reference solution with $N=16$, $\sigma/L=0.1$, and $\delta=0.05$ in the regularized local drag coefficient\ \eqref{eq:creg}, we see a $5\%$ difference in stress using ellipsoidal fibers and a difference of $1$\% and $2$\% for $\delta=0.1$ and $\delta=0.2$, respectively. } }
\end{figure}

We first study the spatio-temporal convergence of stress by fixing $\sigma/L=0.10$.  The errors in the stress under spatio-temporal refinement are shown in Figure\ \ref{fig:Nsigstress}, where we observe rapid convergence consistent with second-order convergence in time. The errors in the stress near $t=0$, immediately after the flow is turned on, are more chaotic since the straight fibers are initially pulled by CLs into a curved shape, and our time step sizes are too large to accurately capture these dynamics.

Perhaps a more important question is the influence of the Gaussian regularization parameter $\sigma$. In Fig.\ \ref{fig:NsigstressVar}, we measure the differences in stress for $N=16$, $\sigma/L=0.1$ and $N=24$, $\sigma/L=0.07$ relative to our (most resolved) reference solution with $N=32$, $\sigma/L=0.05$. We see a relative difference of at most 10\% in the stress tensor for $\sigma/L=0.1$ and 5\% for $\sigma/L=0.07$. 

\rev{We also study the impact of the regularization parameter $\delta$ in the local drag coefficient\ \eqref{eq:creg} on the stress tensor. In Fig.\ \ref{fig:DeltaStressVar}, we plot the relative errors in stress using $\delta=0.05$ as a reference solution. We observe differences of about 5\% in the stress tensor when we use ellipsoidal fibers instead of $\delta=0.05$, and we see the differences in stress decrease as $\delta$ decreases, with a difference of only about 1\% between $\delta=0.1$ and $\delta=0.05$. This shows that the precise value of $\delta$ has a small effect on the macroscopic rheology.}

\subsubsection{Viscoelastic behavior}
We now turn to the measurement of the viscous and elastic moduli for a network of $F=700$ fibers and $8400$ CLs. In order to avoid transient behavior, we first find a steady state configuration of the network by initializing straight fibers with CLs and running the system forward in time using local drag and $\Delta t = 0.005$ without any background flow. \rev{After $t = 2500$ seconds, the maximum $L^2$ norm of the fiber velocity is approximately $4 \times 10^{-6}$ $\mu$m/s, which indicates a near steady state.}

In order to measure the steady state viscous and elastic moduli, we must wait for some intrinsic time on which the network reaches a new steady state in the shear flow. \revp{This relaxation time scale, which we denote by $\tau_c$, combines the cross-linker and fiber relaxation time scales. For cross-linkers, a characteristic time scale of link relaxation is $\tau_\text{CL} = \mu \ell /K_c = 0.5$ seconds for our parameters. For fibers, we assume that the network is sufficiently constrained that the length scale on which the fibers can relax is the mesh size $\ell_m$, or characteristic distance between filaments.}\footnote{\revp{Another possibility is to assume that the fiber relaxation length scale is governed by the distance between two CLs on each fiber. Since there are an average of 24 CL connections on each fiber ($12F$ CLs, each of which binds to 2 fibers), this gives a distance $\ell_c = 2/24 \approx 0.08 \, \mu$m, and a time scale $\mu \ell_c^4/\kappa$ of about $0.005$ seconds, which is much faster than that measured in our numerical experiments. We thank an anonymous referee for suggesting these time scale estimates.}} \revp{We estimate $\ell_m$ by assuming that $L_d/L$ filaments can fit in one direction while $L_d/\ell_m$ can fit in the other two, so that $\ell_m \sim \sqrt{L_d^3/(FL)} \approx 0.21 \, \mu$m, and a characteristic time scale of fiber relaxation is $\tau_F = \mu \ell_m^4/\kappa \approx 0.2$ seconds. Thus our expectation is for $\tau_c$ to be on the order $0.1-1$ second. }

\begin{figure}
\centering
\includegraphics[width=0.45\textwidth]{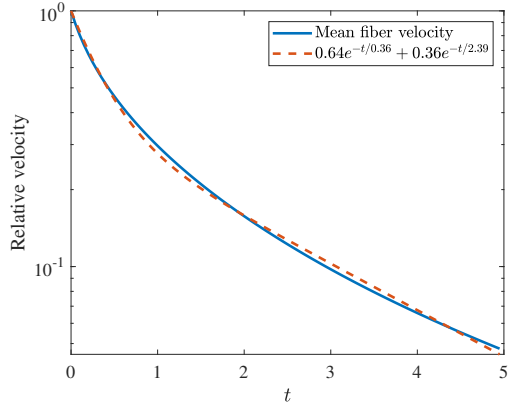}
\caption{\label{fig:VelDecay} Decay of the mean fiber velocity to a steady state after being sheared for one cycle with $\omega=0.2\pi$. We plot the mean fiber velocity, given by\ \eqref{eq:avgfibvel} and normalized by $\bar{v}(t=0)$, over five seconds and compare the result to the double-exponential fit $0.64e^{-t/0.36}+0.36e^{-t/2.39}$ to estimate the network relaxation time scale $\tau_c$.}
\end{figure}

To measure the time scale $\tau_c$ more precisely, we start with the steady state configurations, turn on a shear flow\ \eqref{eq:shear} with $\omega=0.2\pi$ rad/s and $\dot{\gamma}_0=0.02\pi$ 1/s, and run for one cycle (until $T=10$ seconds) with $\Delta t = 0.005$ s. We then turn off the shear flow and measure the velocity of the fibers for another $5$ seconds. We track the mean $L^2$ fiber velocity, given by 
\begin{gather}
\label{eq:avgfibvel}
\bar{v}(t) = \frac{1}{F}\sum_{i=1}^F \ind{v}{i}(t), \quad \text{where} \\ \nonumber
\ind{v}{i}(t) = \left(\sum_{p=1}^N \norm{\ind{\V{X}}{i}_{p}(t)- \ind{\V{X}}{i}_{p}(t+0.05)}^2 \, w_p \right)^{1/2}, 
\end{gather}
and normalize by $\bar{v}(0)$ to obtain the exponential-like decay shown in Fig.\ \ref{fig:VelDecay}. The time scale of relaxation to steady state is $\tau_c \approx 0.5-2$ s, \revp{with the best fit being a sum of two exponentials $\tau_1=0.4$ s and $\tau_2 = 2.4$ s. These relaxation time scales are in line with our physical estimates of $\tau_F = 0.2$ s and $\tau_\text{CL}=0.5$ s.}

\begin{figure}
\centering
\includegraphics[width=0.32\textwidth]{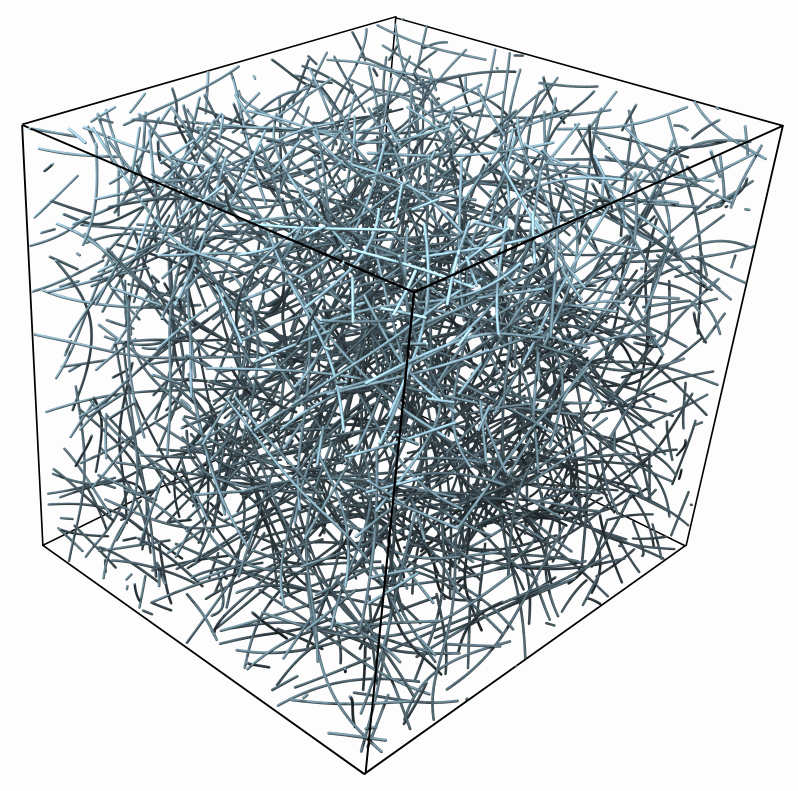}
\includegraphics[width=0.32\textwidth]{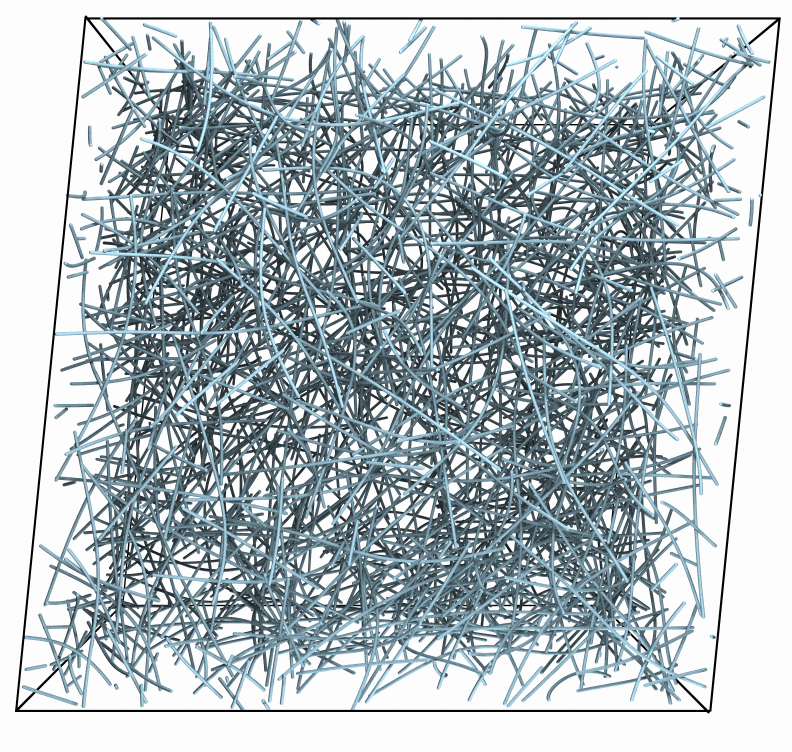}
\includegraphics[width=0.32\textwidth]{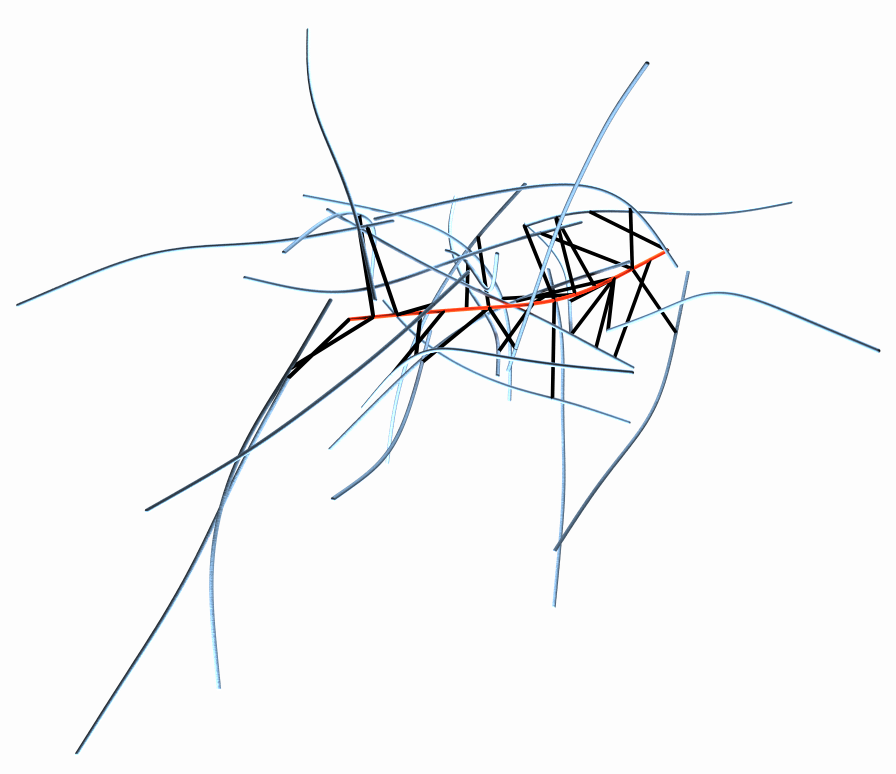}
\caption{\label{fig:network} Steady state fiber configurations for the network of 700 fibers and 8400 CLs  with $\omega=1$ Hz and $g=0.1$. We show (left) a three-dimensional snapshot of all the fibers \rev{in the unit cell}, (middle) a view along the $z$ axis, and (right) a snapshot of all the fibers \rev{(white/blue)} and links \rev{(black)} bound to a single fiber \rev{(orange)} located near the center of the simulation cell. }
\end{figure}

Thus in order to measure the steady state moduli, we wait one second or one cycle (whichever is longer) prior to measuring the stress (and moduli) over three cycles of shear flow. We use a maximum strain $\gamma_0=\dot{\gamma}_0/\omega=0.1$ to stay in the linear regime (data not shown), and we give frequencies in Hz. To initialize for a given $\omega$, we use the final network configuration from the previous (next smallest) frequency. A sample configuration of the fibers in the network, taken with $\omega=1$ Hz at the point of maximum strain, is shown in Fig.\ \ref{fig:network}. \rev{Dynamic movies of the simulations for varying frequencies are available for viewing at \url{https://cims.nyu.edu/~om759/FiberVideos.}}

\begin{figure}
\centering
\subfigure[Moduli with full hydrodynamics]{\label{fig:Moduli700}
\includegraphics[width=0.48\textwidth]{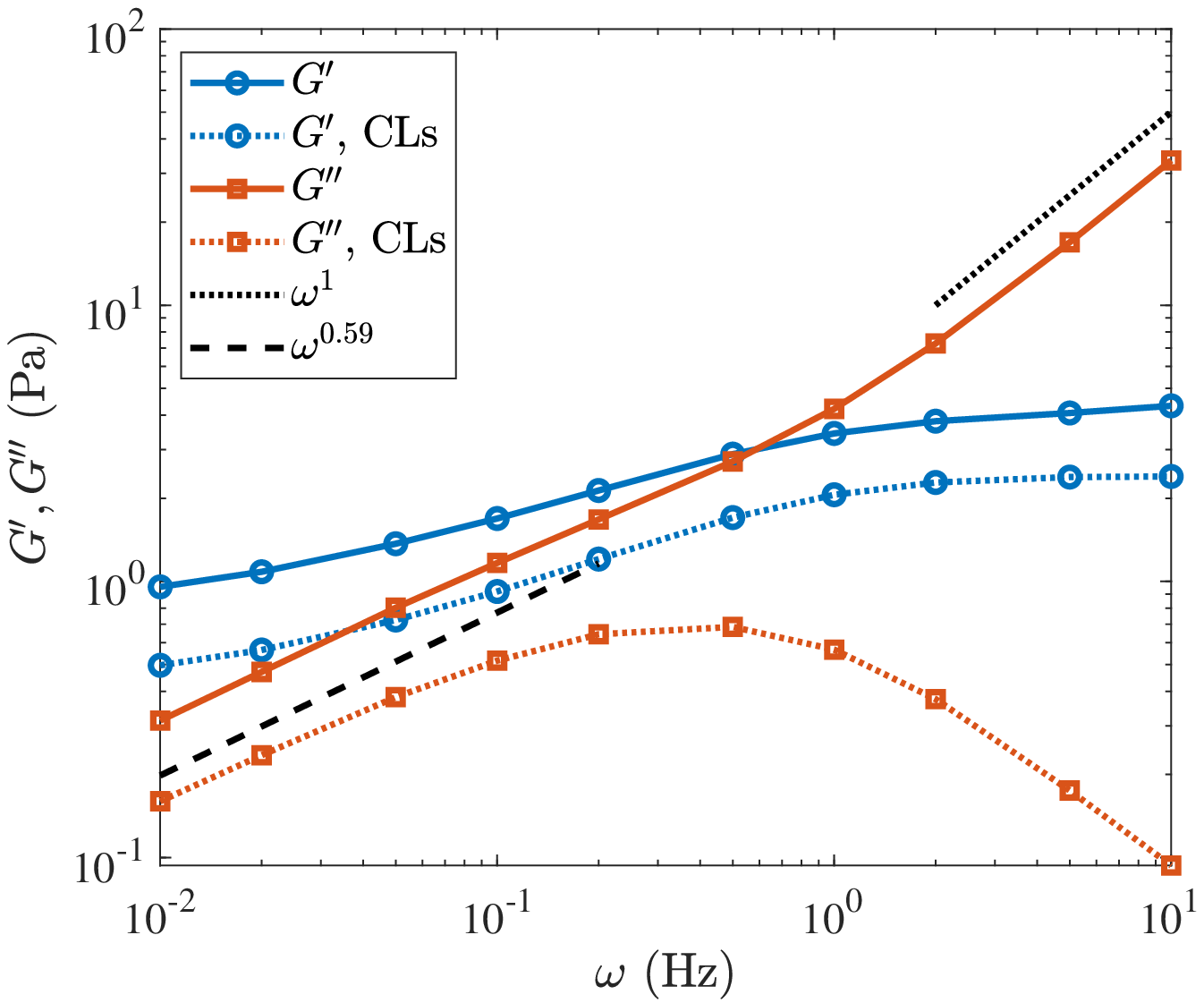}}
\subfigure[Changes in moduli without hydrodynamics]{\label{fig:ChangeModuli700} 
\includegraphics[width=0.48\textwidth]{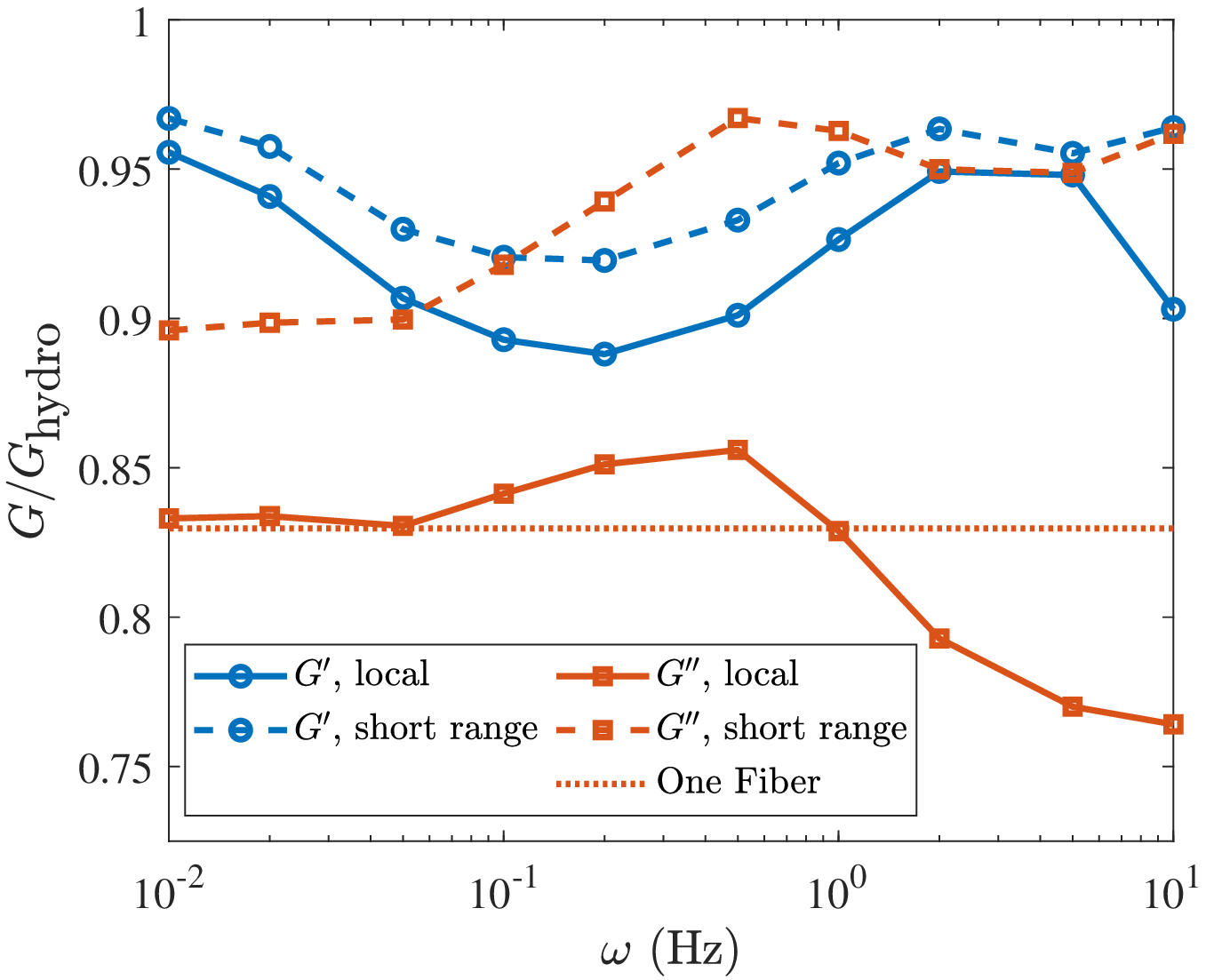}}
\caption{\rev{Elastic and viscous moduli for the network of 700 fibers and 8400 CLs. (a) Values of the elastic modulus $G'$ (blue circles) and viscous modulus $G''$ (orange squares, this excludes the piece due to the background fluid stress $\sigma_{21}^{(\mu)}$ defined in\ \eqref{eq:strdef}) for a gel of 700 fibers and 8400 CLs. We use fully nonlocal hydrodynamics to compute the dynamics and moduli. We show the contribution of the CL stress\ \eqref{eq:sigCL} to the moduli as dotted lines. (b) Fraction of $G'$ (blue circles) and $G''$ (orange squares) that can be recovered when the dynamics are computed using local drag (solid lines) and short-ranged hydrodynamics only (i.e., only intra-fiber but no inter-fiber hydrodynamics). We also show the scaling of the stress for a single fiber in a shear flow as a dotted orange line. }}
\end{figure}

Figure\ \ref{fig:Moduli700} shows the steady-state elastic and viscous moduli when the dynamics of the network \rev{are computed with fully nonlocal hydrodynamics}. \revp{Figure\ \ref{fig:Moduli700} also shows separately the contribution to the moduli of the CL stress\ \eqref{eq:sigCL}, which is $5-10$ times more elastic than viscous.} There is a clear transition in both of the moduli for $\omega \approx 0.5-1$ Hz which can be understood using the characteristic time scales in the problem. Since the characteristic relaxation time scale $\tau_c \approx 0.5-2$ s, the behavior of the moduli can be divided into three regimes: low frequency (background flow time scale $\tau_{\omega} \gg \tau_c$), medium frequency ($\tau_{\omega} \approx \tau_c$), and high frequency ($\tau_{\omega} \ll \tau_c$), where $\tau_\omega = \omega^{-1}$.

\begin{enumerate}
\item \rev{} In the low frequency regime ($\omega < 0.1$ Hz), $\tau_{\omega}$ is the longest time scale and the system is in a constant quasi-static state. If the frequency is small enough, the network has the opportunity to relax at every instant, and it therefore behaves more like an elastic solid where the links constrain the network. As in an elastic solid, the elastic modulus is unchanged with frequency and changes very little (less than $10\%$, as shown in Fig.\ \ref{fig:ChangeModuli700}) when nonlocal hydrodynamics is \rev{dropped}. In this regime, the viscous modulus scales like $G'' = \omega^{0.59}$; the reason for this particular scaling is not obvious to us.
\item In the mid frequency regime ($0.1 \leq \omega \leq 1$), fibers and CLs can deform and relax on the time scale of the background flow and the dynamics involve both an elastic and viscous response. In this regime, $G'' \approx G'$, as shown in Fig.\ \ref{fig:Moduli700}, and the change in the elastic modulus $G'$ due to both changes in frequency (Fig.\ \ref{fig:Moduli700}) and the inclusion of hydrodynamics (Fig.\ \ref{fig:ChangeModuli700}) is maximal. 
\item In the high frequency regime ($\omega > 1$), the background flow dominates the dynamics, the network is essentially fixed on the time scale of the shear flow, and it oscillates back and forth as a viscous fluid would. Figure\ \ref{fig:Moduli700} shows that for $\omega \gg 1/\tau_c$, the viscous modulus scales linearly with $\omega$, as would happen for a pure viscous fluid. In this regime, the viscous modulus \rev{decreases} by as much as \rev{25\%} when nonlocal hydrodynamics is not included, with higher frequencies giving larger decreases. The larger the frequency, the farther the network is from its quasi steady-state and the more important dynamics are for determining the viscous modulus. 
\end{enumerate}

\begin{figure}
\centering
\includegraphics[width=0.41\textwidth]{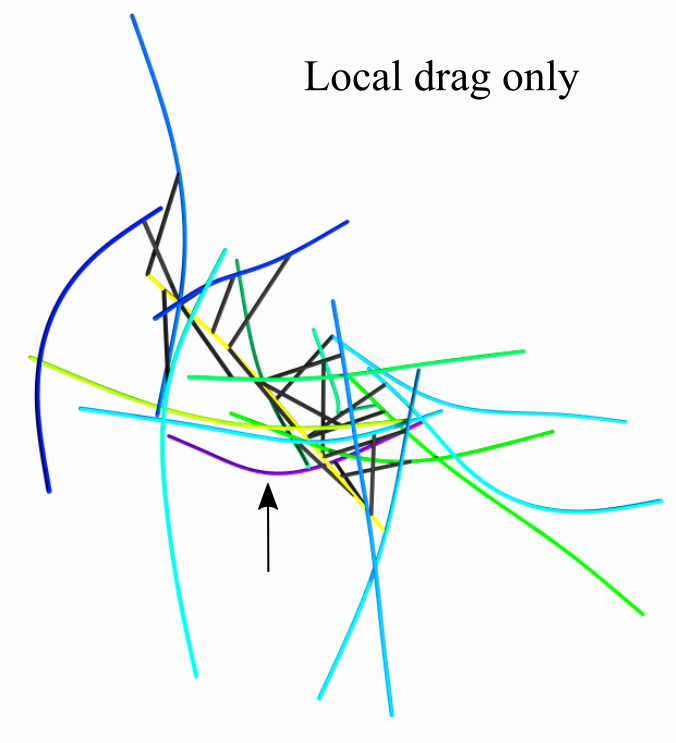}
\includegraphics[width=0.41\textwidth]{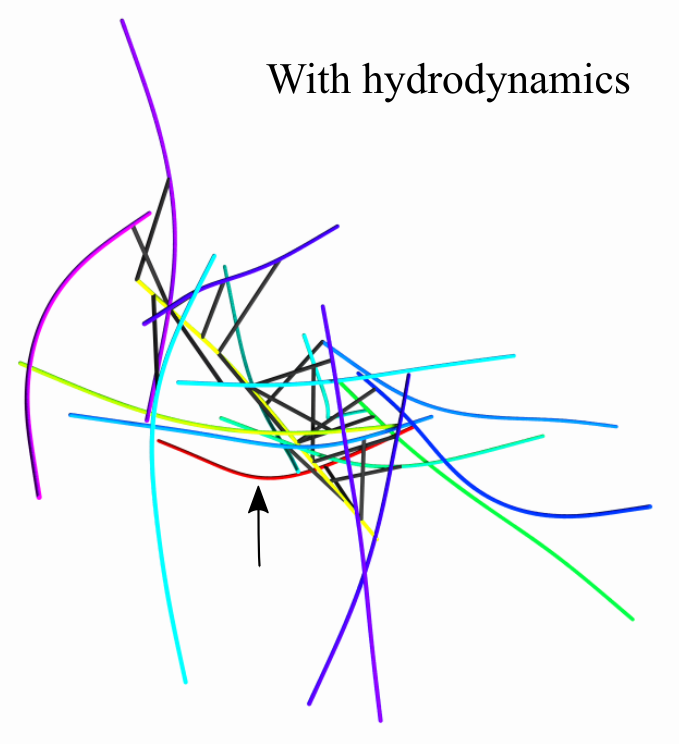}
\includegraphics[height=0.5\textwidth, trim = {13.5cm 0cm 0cm 0cm}, clip]{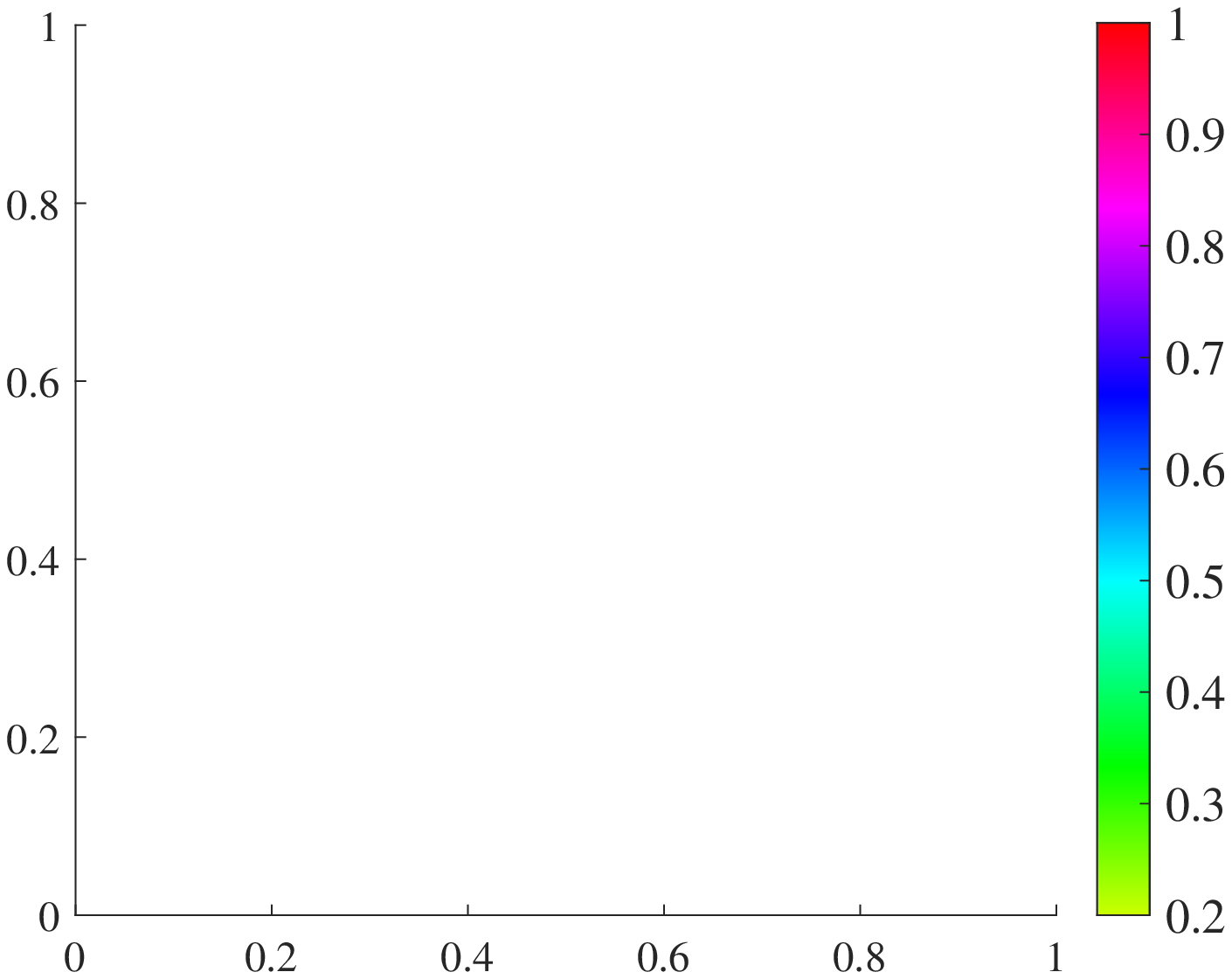}
\caption{\label{fig:lambdasLNL} Subset of fibers inside the unit cell, colored by the $L^2$ norm of the constraint forces $\ind{\V{\lambda}}{i}$ for $\omega=1$ Hz at $t=3.5$ s ($g=0$). Normalization is the maximum $L^2$ norm of $\ind{\V{\lambda}}{i}$ in the system, which is on the fiber indicated by the arrow. We show the case of local drag at left and nonlocal hydrodynamics at right, noting that the inclusion of nonlocal hydrodynamics increases the norm of the constraint forces (which is the dominant cause of the increase in stress), without changing the fiber positions significantly.}
\end{figure}

\revp{Generally speaking, the changes in the viscous modulus when hydrodynamics is switched off are attributable to a reduction in stress on each fiber individually. The changes in the viscous modulus at most frequencies can in fact be explained by considering an isolated straight fiber that makes an angle $\theta$ with the $x$ axis, $\Xs(s)=(\cos{\theta},\sin{\theta},0)$. We put a background shear flow on the fiber $\V{u}_0(x,y,z) = (y,0,0)$ and compute the resulting constraint force density $\V{\lambda}$ on the fiber using both local drag and nonlocal hydrodynamics. Averaging over $\theta$, we obtain a mean difference in the corresponding stress of $\sim17\%$, which is plotted as a dotted line in Fig.\ \ref{fig:ChangeModuli700}, and matches the change in the network's viscous modulus when nonlocal hydrodynamics is dropped, except at the largest frequencies. This implies that the change in the network's viscous modulus comes primarily from an increase in the stress on each fiber individually when the intra-fiber, but not necessarily inter-fiber, nonlocal hydrodynamics is included in the mobility calculation. We confirm this in Fig.\ \ref{fig:lambdasLNL}, which shows that including nonlocal hydrodynamics causes the magnitude of the constraint forces $\V{\lambda}$ on the fibers in the network to increase without changing their positions substantially. 

While Fig.\ \ref{fig:ChangeModuli700} shows that there is a change in the viscous modulus of as much as 25\% when the dynamics are computed by local drag, it also shows that including short-ranged hydrodynamics only, i.e., including intra-fiber but not inter-fiber hydrodynamics, gives a viscous modulus that is at least 90\% of the one computed with inter-fiber hydrodynamics. The changes in the elastic modulus when nonlocal hydrodynamics is dropped are smaller (at most 10\%), but not explainable simply by adding only intra-fiber hydrodynamics. Indeed, for $\omega=0.5$, when the change in the elastic modulus due to hydrodynamics is near maximal, we find that the short-ranged finite part integral\ \eqref{eq:Mfp} explains only one third of the change, with the rest coming from inter-fiber hydrodynamics. Physically, the elastic modulus is related to the interactions of the fibers with CLs, and the network is sufficiently connected that the CLs transmit stress across multiple links in the network, and so long-ranged hydrodynamics plays a role.   }

Our main findings here can be summarized as follows: there exists a critical time scale $\tau_c$, which is on the order of a second for our parameters. On time scales longer than $\tau_c$, the CLs are constantly in a steady state, as they would be in an elastic solid, and the network is more elastic than viscous ($G' > G''$). On time scales shorter than $\tau_c$, the network does not have time to respond to (penalize) strain deformations, and the network is more viscous than elastic ($G'' > G'$). On time scales comparable to $\tau_c$, the network is equally viscous and elastic. The more viscous the network behavior and the farther the network is from a quasi-steady state, the more nonlocal hydrodynamics impacts the moduli. Indeed, if a modulus $G'$ or $G''$ changes substantially with changes in frequency, then we expect dynamics, including whether they are computed with nonlocal or local hydrodynamics, to matter. \revp{Interestingly, we find that intra-fiber nonlocal hydrodynamics dominates the contributions of nonlocal hydrodynamics to the viscous, but not the elastic, modulus. }

\section{Conclusion \label{sec:conclusion}}
In this paper, we have developed a novel method for the simulation of slender filaments, such as actin filaments and microtubules, in the viscous environment of live cells. The key novelty is our reformulation of the continuum fiber centerline evolution in terms of tangent vector rotations, $\partial \Xs/\partial t = \V{\Omega} \times \Xs$, where $\V{\Omega}$ is an unknown rotation rate. We introduced constraint forces $\V{\lambda}$ and closed the system by requiring that the constraint forces perform no work. This virtual work constraint supplements the evolution equation to give a closed saddle-point system for the unknown rotation rates and constraint forces. Here we used slender body theory to obtain the mobility operator $\Lop{M}$ for cytoskeletal filaments in the zero Reynolds number regime; however, our continuum formulation could be used with other mobility relationships. In fact, we show in Appendix \ref{sec:rpyappen} that the fiber evolution equation takes the same form regardless of whether we represent the fiber as a cylinder using SBT or a continuum chain of regularized point forces (using the RPY tensor).

We have also shown that our formulation of inextensibility lends itself more naturally to numerical calculations than formulations involving an auxiliary line tension equation \cite{ts04}. The two main issues with the line tension equation are that it is highly nonlinear in $\V{X}$, which can create aliasing issues when using spectral methods, and that it does not naturally preserve inextensibility after discretization, thereby requiring additional penalty forces to do so. We were able to obtain the discrete tangent vector rotation rates $\V{\Omega}$ directly from the discretized saddle-point system, rotate the fiber tangent vectors, and obtain the positions by integration, thus preserving discrete inextensibility without penalty parameters. 

\rev{While we developed an efficient spectral discretization for each fiber's centerline, the numerical method  as described here breaks down when the fibers become cylindrical. In that case, which corresponds to $\delta \rightarrow 0$ in our local drag regularization\ \eqref{eq:creg}, the contribution of the finite part integral might exceed that of local drag at the fiber endpoints. This causes the SBT mobility matrix to become ill-conditioned, and the force required to produce even a uniform velocity becomes highly oscillatory near the fiber endpoints. The corresponding lack of smoothness in $\V{\Omega}$ and $\V{\lambda}$ at the endpoints makes the global Chebyshev method used here inaccurate and unstable. One way to address this might be to couple a spectral discretization, such as the one developed here, for the fiber interior with a lower-order finite difference discretization near the endpoints. This would allow the smooth part of the forcing in the fiber interior to be represented efficiently using a global interpolant, while the oscillatory parts of the forcing would be represented using a more robust local representation.}

In addition to our contributions for a single fiber, we have also contributed key numerical developments in nonlocal slender body hydrodynamics. By observing \rev{that slender body theory can be interpreted as asymptotic evaluation of a line integral of regularized singularities}, we reformulated the nonlocal hydrodynamics in terms of the Rotne-Prager-Yamakawa mobility tensor. We then used an Ewald splitting scheme for the RPY kernel to evaluate the hydrodynamic interactions on a triply periodic sheared domain in linear time with respect to the number of fibers. For nearby fibers, we supplemented Ewald splitting with a special quadrature scheme which is based on factoring out the near singularity and expanding what remains in a monomial series \cite{barLud}. This special quadrature method allowed us to develop an algorithm to efficiently compute inter-fiber interactions to 3 digits of accuracy most of the time, while avoiding calls to the special quadrature scheme when direct quadrature is sufficiently accurate.

Despite our improved treatment of near fiber interactions, there are still several issues with our treatment of nearly contacting fibers to be explored in future work. Foremost among these is our use of slender body theory itself. Indeed, SBT is designed to \textit{avoid} the length scale $\epsilon L$, which is the same length scale on which the fibers contact each other. While we made some adjustments to the nonlocal hydrodynamics that make it less likely for fibers to cross, these adjustments require too small of a time step to resolve the near-contacts. Nearly contacting fibers also reduce the effectiveness of our GMRES preconditioner for more concentrated suspensions, which is based on the assumption that the dynamics of a single filament are dominated by local drag. Our adjustments to nonlocal SBT also did not account for lubrication \revp{and friction (contact)} forces or steric repulsion between nearby fibers. \revp{These forces, which might be a function of the local orientation, twist, and microscopic shape of the filaments, play a role on length scales $\epsilon L$ and cannot be resolved by slenderness approximations. Our work has therefore left unanswered what contribution these localized forces make to the macroscopic rheology of a cross-linked fiber suspension.} A possible extension to this work is to use a collision tracking algorithm \cite{yan2019computing, yan2020scalable} to preserve larger time step sizes and prevent fibers from crossing, \revp{thereby enabling future studies to determine the role of lubrication, sterics, and friction.}

Since the nonlocal hydrodynamics is the most expensive to compute, we designed a temporal integrator that minimizes the number of nonlocal evaluations. For dilute suspensions, we split the velocity into a local and nonlocal part. The local part, which is the leading order term in $\epsilon$, was treated implicitly, while the sub-leading nonlocal part was treated explicitly. This scheme is stable for any value of the fiber stiffness $\kappa$, as long as the local dynamics dominate the nonlocal hydrodynamics. For semi-dilute and semi-concentrated fiber suspensions, our strategy was to use GMRES to solve for the second-order perturbations in the kinematic coefficients $\V{\alpha}$ and constraint forces $\V{\lambda}$ that give unconditional stability. Since the perturbations are needed for stability and not accuracy, we simply run GMRES for a fixed number of iterations until we get stability. Our tests showed that at most five total evaluations of the nonlocal hydrodynamics are needed per time step to ensure stability, even for semi-concentrated fiber suspensions.

Our preliminary results on cross-linked actin networks in Section\ \ref{sec:CLs} showed that the behavior of the network revolves around a timescale $\tau_c$ on which the network relaxes to a dynamic steady state. On timescales longer than $\tau_c$, the network is in a quasi-static state and is more elastic in nature, while on timescales shorter than $\tau_c$, the elastic part of the network has no time to respond to the deformations exerted on it by a background flow, and the network as a whole shows viscous behavior. On timescales shorter than $\tau_c$, including nonlocal hydrodynamics in the mobility changes the viscous modulus by as much as 25\% in the system we considered. \revp{However, the increases in the viscous modulus come primarily from an increase in stress on each filament separately rather than truly long-ranged hydrodynamic interactions (i.e., the increases due to hydrodynamics are reproducible when each fiber only interacts with itself through the fluid). This behavior differs from that of the elastic modulus, where the more modest increase of $10\%$ due to hydrodynamic interactions is only present with long-ranged hydrodynamics. While these numbers are within the error bounds of most rheology experiments, without our numerical method it would not have been possible to determine them in the first place, since we were able to perform reference simulations and evaluate the errors made by various approximations such as dropping inter-fiber or even intra-fiber nonlocal hydrodynamics (as is often done for simplicity). Omitting long-ranged inter-fiber hydrodynamic interactions, when sufficiently accurate, can speed up our calculations considerably since the more-expensive parts of our algorithm such as Ewald summation and nearly singular quadratures are no longer part of the mobility calculation\footnote{Our temporal integrator can also be simplified in this case and is stable without any additional GMRES iterations.} (this will allow for rapid scanning of parameter regimes in future work). At the same time, in a number of examples in biology, especially when activity is included, the flows generated by individual fibers add constructively to create large-scale macroscopic flows \cite{monteith2016mechanism, goldstein2016batchelor}.  In these systems, far-field hydrodynamics is necessary to capture large-scale fiber-generated flows, and the numerical method we developed here is a necessary tool for such simulations.}

Our findings for cross-linked actin networks are supported by a number of experimental studies. For example, Gardel et al. also observed $G'$ and $G''$ to be on the same order of magnitude for an actin gel of similar density \cite{gardel2004scaling}, and several studies have demonstrated a weak dependence of the elastic modulus $G'$ on frequency \cite{gardel2004scaling,janmey1994mechanical, kasza2010actin}. Janmey et al. found an elastic modulus $G' \approx 1$ Pa for an actin gel with filaments of mean length $L\approx2$ $\mu$m. They also obtained linear scaling of $G''$ with $\omega$ at high frequencies and sublinear scaling at low frequencies, with a transition occurring at a frequency near $1$ rad/s \cite{janmey1994mechanical}. Our agreement with only some of the existing experimental results is natural since we used only one set of parameters, while in vivo or in vitro experimental parameters can vary based on the system. For example, the CL $\alpha$-actinin is about an order of magnitude shorter ($\approx 40$ nm) \cite{ribeiro2014structure} than filamin ($150-200$ nm) \cite{weihing1985filamins}, and most cross-linking proteins have been estimated to be an order of magnitude stiffer \cite{claessens2006actin} than the ones we use here. Another parameter mismatch manifests itself in the viscous modulus, for which we obtained larger values than those reported experimentally \cite{gardel2004scaling,janmey1994mechanical, kasza2010actin} because our underlying fluid viscosity $\mu$ was larger than estimated for living cells \cite{luby1999cytoarchitecture}. While improvements and extensions to our numerical scheme are required to reach the entire span of biological parameters, the platform we have developed here can still be used to test the assumptions behind some of the prior physical theories \cite{gardel2004scaling, janmey1994mechanical, lieleg2009cytoskeletal}.

Our analysis here was also specific to a fixed $\tau_c$, since the CL stiffness $K_c$ and fiber bending constant $\kappa$ were fixed. This begs the question: how does $\tau_c$ change with changes in $\kappa$, $K_c$, and fluid viscosity $\mu$? How do the intrinsic timescales in the problem change when the CLs are dynamically binding and unbinding to the fibers? To fully explore these questions, we need to design an efficient temporal integrator for the CL forces. For stiff CLs, explicit treatment leads to a reduction in the stable time step size, so that $\mathcal{O}(10^5)$ time steps might be necessary to measure stress for frequencies as small as $0.01$ $s^{-1}$. To alleviate this restriction, we will develop a \mbox{(semi-)implicit} temporal integrator for cross-linked networks in future work. We will also investigate the utility of our spring model for shorter, stiffer CLs, which might be better modeled as actual fibers or as rigid connections between fibers.

In this work, we used Euler-Bernouilli beam theory to obtain the elastic force on a fiber due to its curvature. We neglect twist elasticity on the assumption that those deformations are in equilibrium on our timescales of interest \cite{powers2010dynamics}. While this assumption is reasonable for a system of fibers in shear flow, it precludes modeling flagellar beating \cite{lim2012fluid} and chains twisted by external forces such as magnetic fields \cite{sprinkle2020active}. An extension of this work is to account for twist by using the Kirchhoff rod model instead of the Euler beam. The Kirchhoff rod model has previously been used in combination with the immersed boundary method \cite{lim2008dynamics}, RPY tensor \cite{keavRPY}, and method of regularized Stokeslets \cite{olson2013modeling, inexRS} to model a bent twisted fiber interacting with a fluid, with the methods of \cite{inexRS, keavRPY} even preserving discrete inextensibility. As we've discussed at length, these regularization-based models become impractical as the fiber becomes slender since $1/\epsilon$ regularized points are required to properly resolve the fiber thickness. Slender body theory is once again a natural choice, but it is unclear how to account for rotation and torque on the fiber centerline. The most rigorous (and difficult) solution is to use an SBT that accounts for twist, for example the $\mathcal{O}(\epsilon^2)$ slender body fluid velocity of Johnson \cite{johnson} or the SBT of Keller and Rubinow \cite[~Section 10]{krub}. Another option is to take a continuum limit of the RPY tensor for translation-rotation and rotation-rotation coupling \cite{wajnryb2013generalization, keavRPY}, as we have done in Appendix\ \ref{sec:rpyappen} for translation-translation, and use these continuum limits in the grand mobility matrix. A still better approach would be to bridge the gap between the RPY tensor and SBT by deriving an RPY-type tensor for rings which gives the average linear and angular velocity on a ring (cross section of a fiber) due to the force/torque uniformly spread over another ring (cross section).

Because the persistence length of actin is $\mathcal{O}(10)$ $\mu$m and actin filaments in vivo are hundreds of nanometers of length, we have neglected thermal fluctuations in this work. For actin networks, some studies \cite{janmey1994mechanical, gardel2004scaling} have estimated large elastic moduli even in the absence of CLs, which indicates entropic effects. One of our goals for future work is to extend the method here to account for thermal fluctuations of the filaments. There are many challenges in doing this. Most notably, the loss of smoothness of the fiber centerline function $\V{X}(s)$ when thermal noise is present ($\Xs(s)$ has the same H$\ddot{\text{o}}$lder continuity as Brownian motion) makes it difficult to even write a well-posed stochastic equation of motion, let alone solve it accurately using a spectral method suited to smooth $\V{X}(s)$, as we assumed here. 

While there are many improvements awaiting attention, our work has provided solutions to several problems that have previously plagued numerical methods for inextensible slender fibers. From a general framework to view the dynamics of inextensible fibers to specific numerical methods for singular integrals, we have developed a platform for the efficient simulation of thousands of filaments that can be used to gain new insights into processes from sedimentation to cell division and motility.

\section*{Acknowledgments}
We thank those in the numerical analysis community who provided advice on some aspects of this work. Florencio Balboa Usabiaga provided results for the strong formulation that we use in Section\ \ref{sec:falling}. Anna-Karin Tornberg privately shared a preprint of Ref.\ \cite{tornquad} with us, and Ludvig af Klinteberg and Alex Barnett supplied us with code on the special quadrature schemes for nearly singular integration. Nick Trefethen and Nick Hale pointed us to rectangular spectral collocation for the calculation of bending forces. Yoichiro Mori and Laurel Ohm provided help with the theory of SBT. Thanks also to Mike Shelley, Leslie Greengard, and Wen Yan for helpful discussions on the numerical method and stress calculations. Ondrej Maxian is supported by the National Science Foundation via GRFP/DGE-1342536 and Alex Mogilner is supported by U.S. Army Research Office Grant W911NF-17-1-0417. This work was also supported by the NSF through Research Training Group in Modeling and Simulation under award RTG/DMS-1646339.

\appendix
\setcounter{equation}{0}
\renewcommand{\theequation}{\thesection.\arabic{equation}}
\section{Relationship between regularized singularity methods and SBT \label{sec:rpyappen}}
In this appendix, we compare our SBT-based fiber representation with a line of regularized point-like forces. Previous studies have used the method of regularized Stokeslets \cite{inexRS} or RPY tensor \cite{keavRPY} to represent exactly inextensible fibers. An unresolved question is how regularized singularity methods relate to an SBT-based approach for a single fiber. Bringley and Peskin \cite{ttbring08} partially answered this by numerically comparing results from a free space immersed boundary method to SBT. \rev{Cortez and Nicholas \cite{cortez2012slender} gave a more complete answer by asymptotically evaluating the velocity on a fiber made of regularized Stokeslets and doublets. They showed that the resulting velocity is a regularized form of SBT with an arbitrary constant regularization parameter. In this appendix, we extend these prior results on regularized singularity methods by considering the fiber to be a line of singularities regularized with the RPY tensor\ \eqref{eq:rpyknel}. We show that there is a unique choice of sphere radius $\aRPY$ that gives, to order $\epsilon$, the same formula as SBT for the centerline velocity. This is the sphere radius that we use in our inter-fiber mobility\ \eqref{eq:sbtother} to give a consistent formulation.}

Consider a line of regularized point forces, where each of the point forces is regularized over the surface of a sphere of radius $b$, and the centerline velocity at a point on the fiber is computed by averaging the fluid velocity over the same sphere of radius $b$. This means that force can be related to velocity via the RPY mobility tensor. The velocity at arclength coordinate $s$ on the fiber is then obtained by integrating the RPY kernel over the fiber length,
\begin{align}
\label{eq:rpyint}
8\pi \mu\V{U}\left(s\right) =  & \int_{R> 2b} \Knel{\V{X}(s),\V{X}\left(s'\right),\rev{2}b^2/3}\V{f}\left(s'\right) \, ds'\\[2 pt]
\nonumber + & \int_{R \leq 2b} \left(\left(\frac{4}{3b}-\frac{3R\left(s'\right)}{8b^2}\right)\M{I}+\frac{1}{8b^2R\left(s'\right)} \left(\V{R}\V{R}\right)\left(s'\right)\right)\V{f}\left(s'\right) \, ds'. 
\end{align}
Here $\V{R}\left(s'\right)=\V{X}\left(s'\right)-\V{X}(s)$ and $R=\norm{\V{R}}$. The separation of the integrals captures the change in the RPY tensor\ \eqref{eq:rpyknel} when $R < 2b$. Because the RPY tensor is nonsingular for $R=0$, it can be evaluated at any point on the centerline. This is in contrast to the asymptotics for SBT, which are based on evaluating the Stokeslet/doublet kernel on the fiber surface (i.e., $\epsilon L$ away from the centerline), and then assigning this result to be the velocity of the centerline \cite{gotz2001interactions}. Although the kernel\ \eqref{eq:rpyint} is nonsingular, it is still nearly singular for $s \approx s'$, and in the limit $b \ll L$ it is more efficient to evaluate\ \eqref{eq:rpyint} asymptotically. In this appendix we show that this results in an SBT-type formulation with a local drag term and finite part integral. \rev{The local drag terms can be matched with a specific choice of $\aRPY$, and the remainder is a finite part integral that is equivalent to that of SBT to $\mathcal{O}(\epsilon)$. } 

Our strategy is standard matched asymptotics and similar to the approach of Gotz \cite{gotz2001interactions} for SBT. We compute an outer expansion to the integral\ \eqref{eq:rpyint} by considering the region where $|s-s'|$ is $\mathcal{O}(1)$. We then construct an inner expansion in the region where $|s-s'|$ is $\mathcal{O}(b)$. This inner solution must be constructed in two parts for $|s-s'| > 2b$ and $|s-s'| \leq 2b$. We then add the inner and outer solutions together and subtract the common part to obtain the final solution. 

\subsection{Outer expansion}
In the outer expansion, we consider the part of the integral\ \eqref{eq:rpyint} where $|s-s'|$ is $\mathcal{O}(1)$. In this case, the doublet term in $\M{S}_D$ is insignificant and we obtain the outer velocity by integrating the Stokeslet over the fiber centerline, 
\begin{equation}
8\pi \mu \V{U}^{(\text{outer})}(s) = \int_{R > 2b} \Slet{\V{X}(s),\V{X}(s')}\V{f}\left(s'\right) \, ds'. 
\end{equation}
The part of the kernel\ \eqref{eq:rpyint} for $R \leq 2b$ makes no contribution to the outer expansion since $|s-s'|$ is $\mathcal{O}(b)$ there. 

\subsection{Inner expansion}
In the inner expansion, we consider the part of the integral\ \eqref{eq:rpyint} where $|s-s'|$ is $\mathcal{O}(b)$. In this case, we follow \cite{gotz2001interactions} and introduce the rescaled variable
\begin{equation}
\xi = \frac{s'-s}{b}, 
\end{equation}
so that $\xi$ is $\mathcal{O}(1)$. As in \cite{gotz2001interactions}, we will assume that the region $[-2,2]$ is contained in the domain of $\xi$, thereby ignoring the case when $s$ is $\mathcal{O}(b)$ away from the fiber endpoints. \rev{While it is feasible to directly evaluate the RPY integral\ \eqref{eq:rpyint} at the fiber endpoints, our goal here is to show equivalence with SBT, which is only valid away from the fiber endpoints.}

We will need the following asymptotics around $\V{X}(s)$, 
\begin{gather}
\V{R} = \V{X}\left(s'\right)-\V{X}(s) = \xi b \Xs(s) + \mathcal{O}\left(b^2\right), \quad 
\V{R}\V{R} = \xi^2 b^2 \Xs(s) \Xs(s) +\mathcal{O}\left(b^3\right),\\[2 pt] 
R^2 = \V{R} \cdot \V{R} = \xi^2 b^2 + \mathcal{O}\left(b^3\right), \quad R = |\xi|b+\mathcal{O}\left(b^2\right),\\[2 pt]
R^{-1} = \frac{1}{|\xi|b}+\mathcal{O}(1), \quad 
R^{-3} = \frac{1}{|\xi|^3 b^3} + \mathcal{O}\left(b^{-2}\right), \quad 
R^{-5} = \frac{1}{|\xi|^5 b^5} + \mathcal{O}\left(b^{-4}\right). \\[2 pt]
\V{f}\left(s'\right) = \V{f}(s) + \mathcal{O}(b)
\end{gather}

We begin with the part of the integral\ \eqref{eq:rpyint} that uses the kernel $\M{S}_D$ in the region $R> 2b$. For this we will need the expansion of the Stokeslet and doublet, 
\begin{gather}
\Slet{\V{X}(s),\V{X}\left(s'\right)} = \frac{\M{I}+\Xs(s) \Xs(s)}{|\xi| b} + \mathcal{O}(1)\\[2 pt]
\Dlet{\V{X}(s),\V{X}\left(s'\right)} = \frac{\M{I}-3\Xs(s) \Xs(s)}{|\xi|^3 b^3} + \mathcal{O}\left(b^{-2}\right). 
\end{gather}
We now integrate the Stokeslet along the centerline region 
\begin{align}
\int_{R > 2b} \Slet{\V{X}(s),\V{X}(s')}&\V{f}\left(s'\right) \, ds' = \int_{|\xi| > 2}  \left(\frac{\M{I}+\Xs(s) \Xs(s)}{|\xi| b}\right)\V{f}(s)  b \, d\xi +\mathcal{O}(b)\\[2 pt]
& = \left(\M{I}+\Xs(s)\Xs(s)\right)\V{f}(s)\left[\int_{-s/b}^{-2} -\frac{1}{\xi} \, d\xi + \int_2^{(L-s)/b} \frac{1}{\xi}\, d\xi\right]+\mathcal{O}(b)\\[2 pt]
& =\log{\left(\frac{(L-s)s}{4b^2}\right)} \left(\M{I}+\Xs(s)\Xs(s)\right)\V{f}(s)+\mathcal{O}(b). 
\end{align}
Likewise for the doublet, we have
\begin{align}
\frac{\rev{2}b^2}{3} & \int_{R > 2b} \Dlet{\V{X}(s),\V{X}(s')} \, ds' = \frac{\rev{2}b^2}{3} \int_{|\xi| > 2}  \left(\frac{\M{I}-3\Xs(s) \Xs(s)}{|\xi|^3 b^3}\right)\V{f}(s)  b \, d\xi + \mathcal{O}(b)\\[2 pt]
& = \frac{\rev{2}\left(\M{I}-3\Xs(s) \Xs(s)\right)\V{f}(s)}{3} \left[\int_{-s/b}^{-2} -\frac{1}{\xi^3} \, d\xi + \int_2^{(L-s)/b} \frac{1}{\xi^3}\, d\xi\right]+ \mathcal{O}(b)\\[2 pt]
& = \frac{\rev{2}}{12}\left(\M{I}-3\Xs(s)\Xs(s)\right)\V{f}(s)+ \mathcal{O}(b). 
\end{align}
Combining these results, we have, to $\mathcal{O}(b)$, 
\begin{gather}
\label{eq:rg2b}
\int_{R > 2b} \Knel{\V{X}(s),\V{X}\left(s'\right),\rev{2}b^2/3} \, ds' = \\[2 pt] \nonumber \left(\log{\left(\frac{(L-s)s}{4b^2}\right)}\left(\M{I}+\Xs(s)\Xs(s)\right) + \frac{1}{\rev{6}}\left(\M{I}-3\Xs(s)\Xs(s)\right)\right)\V{f}(s).
\end{gather}

It still remains to include in the inner expansion the term for $R \leq 2b$. For this we have the two terms
\begin{gather}
\int_{R < 2b} \left(\frac{4}{3b}-\frac{3R\left(s'\right)}{8b^2}\right)\V{f}\left(s'\right) \, ds' =\V{f}(s)\int_{-2}^2 \left(\frac{4}{3}-\frac{3|\xi|}{8}\right) \, d\xi +\mathcal{O}(b)= \frac{23}{6}\V{f}(s) + \mathcal{O}(b),\\[2 pt]
 \int_{R < 2b} \frac{1}{8b^2R\left(s'\right)} \left(\V{R}\V{R}\right)\left(s'\right)\V{f}\left(s'\right) \, ds' = \frac{1}{8}\int_{-2}^2 \Xs(s)\Xs(s)\V{f}\left(s'\right) |\xi| \, d\xi +\mathcal{O}(b)= \frac{1}{2}\Xs(s)\Xs(s)\V{f}(s)+\mathcal{O}(b)
\end{gather}
where we have used the fact that $\xi \in [-2,2]$ is on the fiber ($s$ is away from the endpoints). We therefore have, to $\mathcal{O}(b)$, 
\begin{gather}
\label{eq:rl2b} 
\int_{R \leq 2b} \left(\left(\frac{4}{3b}-\frac{3R\left(s'\right)}{8b^2}\right)\M{I}+\frac{1}{8b^2R\left(s'\right)} \left(\V{R}\V{R}\right)\left(s'\right)\right) \, ds' = \left(\frac{23}{6}\M{I} + \frac{1}{2}\Xs(s)\Xs(s)\right)\V{f}(s). 
\end{gather}
The inner expansion is therefore, adding the terms\ \eqref{eq:rg2b} and\ \eqref{eq:rl2b}, 
\begin{equation}
8 \pi \mu \V{U}^{(\text{inner})}(s) = \left(\log{\left(\frac{(L-s)s}{4b^2}\right)}\left(\M{I}+\Xs(s)\Xs(s)\right)+ \rev{4}\M{I}\right)\V{f}(s). 
\end{equation}
By adding and subtracting $\log(16)\left(\M{I}+\Xs(s)\Xs(s)\right)$, we obtain the same leading order coefficient as SBT, 
\begin{gather}
8\pi \mu \V{U}^{(\text{inner})}(s)=  \Bigg{(}\log{\left(\frac{4(L-s)s}{b^2}\right)}\left(\M{I}+\Xs(s)\Xs(s)\right) \\[2 pt]
\nonumber + \left(\rev{4}-\log{16}\right)\M{I} - \rev{\left(\log{16}\right)}\Xs(s)\Xs(s)\Bigg{)}\V{f}(s). 
\end{gather}

\subsection{Common part}
The common part is the outer velocity written in terms of the inner variables. That is, to $\mathcal{O}(b)$, 
\begin{equation}
8\pi \mu \V{U}^{(\text{common})}(s) = \int_{R > 2b} \left(\frac{\M{I}+\Xs(s)\Xs(s)}{|s-s'|}\right)\V{f}(s) \, ds'. 
\end{equation}

\subsection{Matched asymptotic expansion}
The total velocity is the sum of the inner and outer expansions, with the common part subtracted, 
\begin{gather}
\V{U}(s) = \V{U}^{(\text{inner})}(s) + \V{U}^{(\text{outer})}(s)-\V{U}^{(\text{common})}(s). 
\end{gather}
This can be written as 
\begin{gather}
\label{eq:totvelnofp}
8\pi\mu\V{U}(s) = \left(\log{\left(\frac{4(L-s)s}{b^2}\right)}\left(\M{I}+\Xs(s)\Xs(s)\right)
+ a_I \M{I} + a_\tau \Xs(s)\Xs(s)\right)\V{f}(s)\\[2 pt] \nonumber + \int_{R> 2b} \left(\Slet{\V{X}(s),\V{X}(s')} \V{f}\left(s'\right) -  \left(\frac{\M{I}+\Xs(s)\Xs(s)}{|s-s'|}\right)\V{f}(s)\right) \, ds',\\[2 pt]
\text{where} \quad \rev{a_I = 4-\log{16} \quad \text{and} \quad a_\tau = -\log{16}}. 
\end{gather}
We can now establish equivalence with SBT by observing that the integrand in\ \eqref{eq:totvelnofp} is $\mathcal{O}(\aRPY)$ when $R < 2\aRPY$, and so we can add that part of the integral back into the velocity without changing the asymptotic accuracy of the velocity\ \eqref{eq:totvelnofp}. This gives a velocity of the exact same form as SBT,
\begin{gather}
\label{eq:marpy}
8\pi\mu\V{U}(s) = \left(\log{\left(\frac{4(L-s)s}{b^2}\right)}\left(\M{I}+\Xs(s)\Xs(s)\right)
+ a_I \M{I} + a_\tau \Xs(s)\Xs(s)\right)\V{f}(s)\\[2 pt] \nonumber + \int_0^L \left(\Slet{\V{X}(s),\V{X}(s')} \V{f}\left(s'\right) -  \left(\frac{\M{I}+\Xs(s)\Xs(s)}{|s-s'|}\right)\V{f}(s) \, \right)ds'. 
\end{gather}
The velocity expression\ \eqref{eq:marpy} is the same as the SBT velocity\ \eqref{eq:onefib} when $a_I=1$ and $a_\tau=-3$. \rev{This equivalence is accomplished by the specific choice\footnote{Observe that we have two equations ($a_I$ and $a_\tau$) for one variable $\aRPY$, so the existence of a solution is surprising.} of $\aRPY$
\begin{equation}
\label{eq:rpyeps}
\aRPY = \frac{e^{3/2}}{4}\epsilon L,
\end{equation}
which gives the RPY doublet coefficient
\begin{equation}
\label{eq:dco}
\frac{2b^2}{3} = \frac{e^3}{24} \epsilon L, 
\end{equation}
that we use in Section\ \ref{sec:RPYSBT} for inter-fiber interactions. Because we used only the leading order terms in the asymptotics, we are guaranteed that\ \eqref{eq:marpy} approximates\ \eqref{eq:rpyint} to $\mathcal{O}(\aRPY)$, although we observe $\mathcal{O}(b^2)$ accuracy a distance $\mathcal{O}(1)$ from the endpoints.}

\section{Verification for sheared unit cell \label{sec:shearBC}}
\begin{figure}
\centering
\includegraphics[width=70mm]{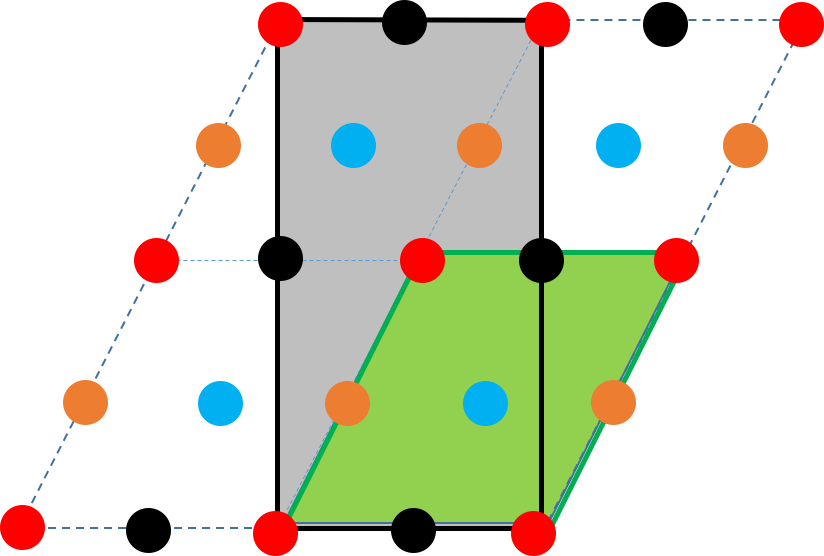}
\caption{\label{fig:hexlat} The lattice for the sheared unit cell test. The points are positioned on a lattice with a 1-particle (primitive) unit cell (shaded green) at $(0,0,0)$ (red), $(1,0,0)$ (black), $(0.5,1,0)$ (orange), and $(1.5,1,0)$ (blue). The lattice can also be viewed as periodic on the larger 2-particle rectangular unit cell shaded in gray. Colors indicate the magnitude of the force placed on each set of points: $+1$ in all three directions on the red points, $-1$ in each direction on the black points, $+2$ on the orange points, and $-2$ on the blue points.}
\end{figure}

\setcounter{equation}{0}
To test our implementation of sheared periodic boundary conditions, we consider a packing of points that is hexagonal in the $xy$ plane. As shown in Fig.\ \ref{fig:hexlat}, the points are positioned on a (green) periodic slanted cell at $(0,0,0)$ (red blob), $(1,0,0)$ (black), $(0.5,1,0)$ (orange), and $(1.5,1,0)$ (sky blue). To form a periodic hexagonal packing in the $xy$ plane, we set $g=0.5$ with periodic domain length $L_x=L_y=L_z=2$. Using a coloring scheme (see Fig.\ \ref{fig:hexlat}), it is easy to see that this arrangement is equivalent to the same set of points on a (gray) rectangular unit cell, with additional points at $(1,2,0)$, $(0,2,0)$, $(1.5,3,0)$, $(0.5,3,0)$, with the ordering of forces in the second set of points being the same as the first and periodic length $L_y=4$. We place a force of strength $+1$ in each direction (including $z$) on the first (red) pair of points, $-1$ in each direction on the second (black) pair, $+2$ on the third (orange), and $-2$ on the fourth (sky blue). Note that the $z$ direction is also periodic in all cases with length $L_z=2$, so that we are actually considering a set of stacked copies of Fig.\ \ref{fig:hexlat}. 

We solve for the RPY velocities induced by the forces at each point using the Ewald splitting technique described in Section\ \ref{sec:ewald}. We set $\xi=5$, sphere radius $b=10^{-2}$, and fluid viscosity $\mu=3$. The maximum relative 2-norm error in the velocity of the four points is less than $10^{-5}$ for all values of the NUFFT tolerance less than $10^{-2}$, with decay to $10^{-11}$ when the tolerance is $10^{-8}$. We conclude that our modified Ewald splitting scheme of Section\ \ref{sec:ewald} properly treats the strain in the periodic coordinate system. 

\section{Near fiber accuracy \label{sec:nearfibacc}}
\setcounter{equation}{0}
In this appendix, we test the accuracy of Algorithm\ \ref{fig:algflow} for computing the slender body interaction integrals\ \eqref{eq:sbtother}. We generate 100 smooth inextensible fibers by initializing an unnormalized tangent vector that is an exact Chebyshev series with 15 exponentially decaying terms. More precisely, the $k$th coefficient of the series is a Gaussian random variable with mean $0$ and standard deviation $e^{-10k/N}$, where $k=0, \dots 15$. We then normalize this tangent vector to obtain $\Xs(s)$ and integrate to obtain the fiber positions $\V{X}(s)$. To make sure the resulting fiber is smooth after tangent vector normalization, we compute the Chebyshev series of the fiber position $\V{X}(s)$. Denoting the coefficients of the position Chebyshev series by $\hat{a}_k$, we only accept fibers with Chebyshev series coefficients $|\hat{a}_k| \leq e^{-0.61k}$ for $k=2, \dots 15$ (the constant and linear modes play no role in the fiber smoothness). This means that the last coefficient $\hat{a}_{15}$ has value at most $10^{-4}$. 

We consider fibers with $L=2$ and $\epsilon=10^{-3}$ in a fluid of viscosity $\mu=1/8\pi$.  Our goal is to evaluate the velocity due to a fiber $\V{X}(s)$ at a target $\V{x}$,
\begin{equation}
\label{eq:vtarg}
\V{v}(\V{x}) = \int_0^L \Knels{\V{x}}{\V{X}(s)}{\rev{\frac{e^3}{24}}(\epsilon L)^2}\V{f}(s). 
\end{equation}
We choose $\V{f}(s)=\Xs(s)$, so that the force density is sufficiently smooth. 

To measure the accuracy of Algorithm\ \ref{fig:algflow}, we place 100 targets a distance $d$ away in a random normal direction from each fiber's centerline. To get a reference answer, we compute the integral\ \eqref{eq:vtarg} directly by upsampling the fiber to 6000 type 1 Chebyshev points. We then compute the integral using Algorithm\ \ref{fig:algflow}. We show the maximum relative error, $\V{E}_i /\norm{\V{v}(\V{x})}_\infty$, where the maximum is over the direction $i=1, 2, 3$ and $\bm{E}$ is the absolute difference between the approximate and reference values of the velocity\ \eqref{eq:vtarg}. 

\begin{figure}
\centering
\subfigure[Short distances]{\label{fig:testshort}
\includegraphics[width=0.45\textwidth]{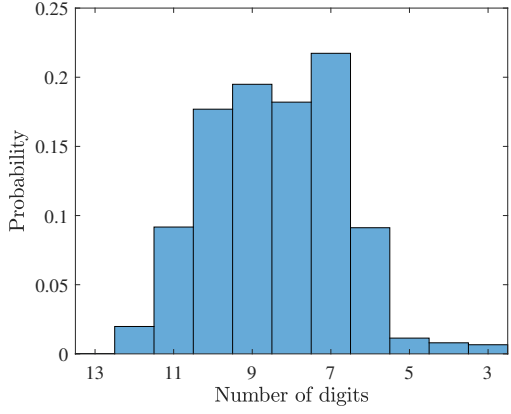}}
\subfigure[Long distances]{\label{fig:testlong}
\includegraphics[width=0.45\textwidth]{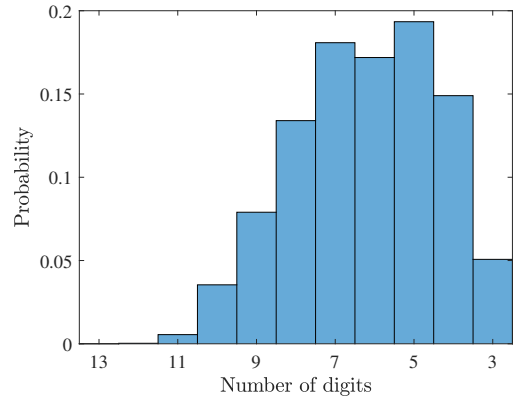}}
\caption{\label{fig:allds}Errors in Algorithm\ \ref{fig:algflow} for near singular quadratures. We show a histogram of the number of digits obtained in the velocity $\V{v}(\V{x})$, given in\ \eqref{eq:vtarg}, for 100 different fibers and 100 targets per fiber. Here $d$ is the distance from the target point to the fiber centerline and we show histograms of the number of digits obtained in the integral $\V{v}(\V{x})$. The number of digits is computed as $-\text{log}_{10}\left(\max_i \V{E}_i /\norm{\V{v}(\V{x})}_\infty \right)$, where the maximum is over the direction $i=1, 2, 3$ and $\bm{E}$ is the absolute difference between the approximate and reference values of the velocity\ \eqref{eq:vtarg}. (a) Errors from short distances $2 \leq d/(\epsilon L) \leq 10$. We see that we are over-working in most cases, since most of the time we obtain many more than 3 digits of accuracy, but there are some cases when we only obtain 3 digits. (b) Long distances $0.04 < d/L < 0.20$, where we obtain 4-7 digits most of the time.}
\end{figure}

We separate our results into short distances, $2\epsilon L < d < 10\epsilon L$, and long distances, $0.01L < d< 0.2L$ (there is overlap between the two regions since $10\epsilon L = 0.01L = 0.02$ with our parameters). Fig.\ \ref{fig:testshort} shows the errors at short distances. We see that we obtain many more digits than necessary in most cases. There are, however, a few cases where we obtain 3 digits. Since $d$ will rarely be $\mathcal{O}(\epsilon L)$, it is acceptable to expend extra computational effort to guarantee accuracy. In Fig.\ \ref{fig:testlong}, we show the errors for long distances. In particular, we see that we obtain 4-5 digits most of the time, and that $\approx 5\%$ of the time we obtain exactly 3 digits of accuracy.  

\bibliographystyle{plain}

\bibliography{PaperTex/SlenderBib}

\end{document}